\documentclass[11pt]{article}
\usepackage{amsmath}
\usepackage{amssymb}
\usepackage{amscd}
\usepackage{amsthm}

\def\endproof{\relax\ifmmode\expandafter\endproofmath\else
  \unskip\nobreak\hfil\penalty50\hskip.75em\hbox{}\nobreak\hfil\bull
  {\parfillskip=0pt \finalhyphendemerits=0 \bigbreak}\fi}
\def\endproofmath$${\eqno\bull$$\bigbreak}
\def\bull{\vbox{\hrule\hbox{\vrule\kern3pt\vbox{\kern6pt}\kern3pt\vrule}\hrule}}
\newtheorem{theorem}{Theorem}[subsection]
\newtheorem*{main}{Theorem}
\newtheorem*{mainc}{Corollary}
\newtheorem{proposition}[theorem]{Proposition}
\newtheorem{lemma}[theorem]{Lemma}
\newtheorem{claim}[theorem]{Claim}

\newtheorem{corollary}[theorem]{Corollary}

\newtheorem{D}[theorem]{Definition}
\newenvironment{defn}{\begin{D} \rm }{\end{D}}

\newtheorem{R}[theorem]{Remark}
\newenvironment{remark}{\begin{R}\rm }{\end{R}}

\def\Zee{\mathbb{Z}}
\def\Q{\mathbb{Q}}
\def\Ar{\mathbb{R}}
\def\Cee{\mathbb{C}}
\def\Pee{\mathbb{P}}
\def\WP{\mathbb{WP}}

\def\scrO{\mathcal{O}}
\def\ov{\overline}
\def\spcheck{^{\vee}}
\def\frak{\mathfrak}
\def\Pic{\operatorname{Pic}}
\def\Ker{\operatorname{Ker}}

\def\Sym{\operatorname{Sym}}
\def\Hom{\operatorname{Hom}}
\def\Aut{\operatorname{Aut}}

\def\Id{\operatorname{Id}}

\def\Spec{\operatorname{Spec}}
\def\ad{\operatorname{ad}}

\def\Lie{\operatorname{Lie}}

\title{Holomorphic Principal Bundles Over Elliptic Curves II:\\
The Parabolic Construction}
\author{Robert Friedman\thanks{The first author was partially
    supported by NSF grant DMS-99-70437.}
\  and John W. Morgan\thanks{The second author was partially supported
   by NSF grant DMS-97-04507.}}

\begin{document}

\maketitle

\section*{Introduction}

Let $E$ be a smooth elliptic curve and let $G$ be a simple  complex algebraic
group of rank $r$. We shall always assume that $\pi_1(G)$ is cyclic and
that $c$ is a generator. We shall freely identify $c$ with the corresponding
element in the center of the universal cover $\widetilde G$ of $G$. A
$C^\infty$
$G$-bundle
$\xi_0$ over
$E$ has a characteristic class $c_1(\xi_0)\in H^2(E;\pi_1(G)) \cong \pi_1(G)$ which
determines $\xi_0$ up to $C^\infty$ isomorphism. The goal of this paper is to
continue the study, begun in
\cite{FM}, of the moduli space $\mathcal{M}(G,c)$ of semistable holomorphic
$G$-bundles $\xi$ with $c_1(\xi) = c$. In
\cite{FM}, this space was studied from the transcendental viewpoint of
$(0,1)$-connections using the results of Narasimhan-Seshadri and
Ramanathan that in every S-equivalence class there is a unique
representative whose holomorphic structure is given by a flat connection.
This viewpoint, however, is not suitable for many questions,
such as finding universal bundles,  studying singular elliptic curves, or
generalizing to families of elliptic curves. In this paper, which is
largely independent of \cite{FM}, we describe
$\mathcal{M}(G,c)$ from an algebraic point of view.
As we shall show in later papers, this construction
is much better adapted for dealing with the questions described above.

Our motivation comes from the theory of deformations of singularities with a
$\Cee^*$-action. In this theory, the deformations corresponding to
nonnegative weights are topologically equisingular. Thus, from the point of
view of smoothings, it is the negative weight deformation space which is
interesting. This space can be studied infinitesimally, by looking at the
$\Cee^*$-action on the Zariski tangent space to the deformation functor.
There is also a globalization  of this local description. For example, in
the case of hypersurface simple elliptic singularities which are given as
weighted cones over an elliptic curve $E$, the global moduli space, modulo the
action of
$\Cee^*$, can (essentially) be identified with pairs $(X,D)$, where $X$ is a
smooth del Pezzo surface or a del Pezzo surface with rational double points,
and $D$ is a hyperplane section of $X$ isomorphic to $E$.

By analogy, our method is to describe $\mathcal{M}(G,c)$ as a certain set of
deformations of a ``singular" bundle $\eta$, where in this context singular
means unstable. Of course, it is natural to require that $\eta$ be
minimally unstable, in the sense that all deformations which are not roughly
speaking topologically equivalent to $\eta$ should be semistable. As
we showed in \cite{FMAB},
such minimally unstable bundles always exist, and, in most cases,
are unique once we fix the determinant. There is a distinguished
subgroup $\Cee^*$ in the
automorphism group of $\eta$, and (with our conventions) it acts with
nonnegative weights. The positive weight deformations of $\eta$
correspond to semistable bundles.

More precisely, to every unstable bundle $\eta$ there is associated a
conjugacy class of parabolic subgroups of $G$.
The parabolic subgroup corresponding to a minimally unstable $\eta$ is
a maximal parabolic subgroup.
The conjugacy classes of maximal parabolic
subgroups are indexed by the simple roots of $G$. Let $P$ be such a maximal
parabolic subgroup, with unipotent radical $U$ and Levi factor $L$. Then $L$
has a one-dimensional center, and hence there is a natural inclusion of
$\Cee^*$ in $L$. Correspondingly there is a unique primitive dominant
character $\chi_0$ of $L$.  We consider semistable
$L$-bundles
$\eta$ whose degree is $-1$ with respect to $\chi_0$. In particular, this
implies that $P_-$, the parabolic subgroup opposite to $P$, is a
Harder-Narasimhan parabolic subgroup for the unstable $G$-bundle
$\eta\times _LG$. The central
subgroup $\Cee^*$ of $L$ acts on the deformation space of $\eta\times _LG$,
and we shall be interested in the positive weight deformations. As in the
singularities case, there is a global interpretation of this infinitesimal
picture: it is the space of isomorphisms classes of pairs
$(\xi,\varphi)$ consisting of a $P$-bundle $\xi$ and an isomorphism
$\varphi\colon \xi/U \to \eta$.
The isomorphism classes of such pairs $(\xi, \varphi)$  are
classified by the cohomology group $H^1(E; U(\eta))$, where $U(\eta)$ is the sheaf
of unipotent groups $\eta\times_LU$. In general, $U(\eta)$ is not a sheaf of
abelian groups, and so $H^1(E; U(\eta))$ is \textit{a priori} simply a set.
However, the fact that $U(\eta)$ has a filtration whose successive quotients
are vector groups implies, in our situation, that $H^1(E; U(\eta))$ carries
in a functorial sense the structure of an affine space. The group $\Cee^*$
acts on this space, fixing the origin and such that the differential of the
action has positive weights at the origin. We show that this action can
be linearized. From this point of view, it is the non-abelian nature of
$U$ which allows there to be different weights for the $\Cee^*$-action, so
that the quotient $(H^1(E; U(\eta))-\{0\})/\Cee^*$ is a weighted projective space
$\WP(\eta)$, typically with distinct weights. In fact,
choosing $\eta$ to be
minimally unstable, the weights are given as follows in  case
$G$ is simply connected.
 Let $\widetilde \Delta$ be the
extended set of simple roots for $G$, and for each $\alpha \in \widetilde\Delta$,
let
$\alpha\spcheck$ be the corresponding dual coroot. Then there is a unique
linear relation
$$\sum _{\alpha \in \widetilde \Delta}g_{\alpha}\alpha\spcheck = 0,$$
provided that we require that the coefficient of the coroot dual to the negative
of the  highest root is $1$. Up to multiplying by a common factor, the weights of
$\WP(\eta)$ are then the integers
$g_\alpha$, counted with multiplicity. In particular, the dimension of $\WP(\eta)$
is equal to the rank of $G$, which we know to be the dimension of
$\mathcal{M}(G,1)$.  When $G$ is non-simply connected and we 
consider  $G$-bundles whose first Chern class is a generator $c\in
\pi_1(G)$, there is a similar result where, again up to a common factor, the
weights are given by the coroot integers $g_{\ov \beta}$ on the quotient of
the extended Dynkin diagram of
$G$ by the action of the central element $c$ in the universal covering
group $\widetilde G$ of $G$ as described in \cite{BFM}.

Assuming that the bundle $\eta\times _LG$ described above is minimally
unstable, all points of $\WP(\eta)$ correspond to semistable bundles.
Thus, there is an induced morphism $\Psi\colon \WP(\eta)\to
\mathcal{M}(G,c)$. Now every morphism whose
domain is a weighted projective space is either constant or finite, and in
our case it is easy to see that the morphism is non-constant. This already
yields some information about $L$ and
$\eta$: assuming
that the action of $\Aut _L\eta$  on $\WP(\eta)$ is essentially effective, since
$\Psi$ is constant on the $\Aut _L\eta/\Cee^*$-orbits, it follows that
$\Cee^*$ is the identity component of $\Aut_L\eta$. To go further, we need to
use the fact that $E$ is an elliptic curve, which implies that the map
$\Psi\colon \WP(\eta)\to  \mathcal{M}(G,c)$ is dominant, and more generally that
the rational map from  $\WP(\eta)$ to $\mathcal{M}(G,c)$ is dominant for every
maximal parabolic. This implies that $\dim \WP(\eta) \geq \dim \mathcal{M}(G,c)$,
with equality for the minimally unstable case. By contrast, for curves of genus at
least
$2$, even for the case of $SL_2(\Cee)$, the map $\WP(\eta) \to \mathcal{M}(G,c)$ is
typically not dominant, but rather maps $\WP(\eta)$ to a proper subvariety of
$\mathcal{M}(G,c)$.

To sum up, then, in the minimally unstable case, we have a weighted
projective space $\WP(\eta)$ and a finite morphism $\Psi\colon \WP(\eta)\to
\mathcal{M}(G,c)$. In fact, we show the following theorem:

\begin{main} Suppose that the $L$-bundle $\eta$ is minimally unstable.
Then the map $\Psi\colon \WP(\eta)\to \mathcal{M}(G,c)$ is an isomorphism.
\end{main}

\begin{mainc} The moduli space $\mathcal{M}(G,c)$ is a weighted projective
space, with weights $g_{\ov \beta}/n_0$, where $n_0$ is the $\gcd$ of the
$g_{\ov \beta}$. In particular, if $G$ is simply connected, $\mathcal{M}(G,c)$ 
is a weighted projective space, with weights $g_\alpha, \alpha\in \widetilde
\Delta$.
\end{mainc}

In the simply connected case, the corollary is due to Looijenga \cite{Loo}
(see also \cite{BS}). Note however
that the theorem goes beyond an abstract description of $\mathcal{M}(G,c)$ as 
a weighted projective space: it identifies an algebraically defined moduli
space, $\WP(\eta)$, with a transcendentally defined  moduli space which in some
sense is obtained by taking the periods of a flat connection.
We view this as  a theorem of Torelli type in a non-linear context.

The bundles produced by the parabolic construction,
in addition to being semistable, are regular in the sense that their
automorphism groups have minimal possible dimension. This is
reminiscent of the Steinberg cross-section of regular elements for the
map of $G$ to its adjoint quotient
and of the Kostant section of regular elements
for the adjoint quotient of the Lie algebra $\frak g$ of $G$.
In fact, as we shall show in a future paper, the parabolic construction
extends to the case of nodal curves of genus one and to cuspidal curves of
genus one when $G\not= E_8$. For nodal
curves, the parabolic construction
produces a weighted projective space and an open subset which is
identified with a Steinberg-like cross-section  of regular elements in
each conjugacy class. For cuspidal curves (and $G\not= E_8$)  the
weighted projective space contains an open subset producing a
Kostant-like section of the adjoint quotient of $\frak g$.
Thus, in both cases, the parabolic construction yields a new approach to the
proof of the existence of  a section of regular elements for the adjoint
quotient, and produces  a natural compactification of the adjoint quotient which
is a weighted projective space.

 In case $G=E_8, E_7, E_6$ (as well as $D_5$ and $A_4$), the
 relationship  between deformations of minimally unstable $G$-bundles
 and deformations of
simple elliptic singularities goes far beyond a formal
 analogy. Indeed, this observation, which is connected to what is
 called in the physics literature F-theory, was a major motivation for
 us to study $G$-bundles over elliptic curves. This connection will be
 described elsewhere.

This paper is related to \cite{FMAB}, where we enumerate the  minimally
unstable strata in the space of
$(0,1)$-connections on a $G$-bundle. While we make no use of these results, at
least in the simply connected case, that paper helps explain the
characterizing properties that a minimally unstable bundle
$\eta$ has to satisfy: it predicts, for example, that  its Harder-Narasimhan
parabolic is the maximal parabolic associated with what we call a
special root and that the $L$-bundle has degree $-1$ with respect to the
dominant character. But even without knowing that these bundles lie in
minimally unstable strata in the Atiyah-Bott formalism,  one can
establish the isomorphism given in the main theorem above.

The contents of this paper are as follows. In Section 1, we collect together
preliminary technical results. Many of these
results concern numerical facts related to irreducible representations
of the Levi factor of a maximal parabolic on the unipotent radical of
that parabolic.
These are used in Section 2 to compute the dimensions of various
cohomology groups related to
bundles over maximal parabolic subgroups of $G$, as well as to understand the
weights of the $\Cee^*$-action. These dimensions and weights could be
computed by case-by-case checking of the root tables. We have tried instead
to find classification-free arguments wherever possible. In calculating the
$\Cee^*$-weights, we make use of a property we call circular symmetry. This
property was introduced by Witten and established for the coroot integers in
\cite{BFM}. Its name derives from the relation of this property to a
symmetry statement for points placed on a circle according to these numbers,
as described in  \cite[\S 3.8]{BFM}. Here, we do not need this geometric
interpretation. Rather we need to know only that (as was proved in \cite{BFM})
the coroot integers and the coroot integers on the quotient diagram by
the action of the center satisfy circular symmetry,
and that numbers satisfying circular symmetry are completely
determined by three pieces of information, which in our context are the
dimension of the weighted projective space, the highest weight
appearing in the $\Cee^*$-action and the dual Coxeter number of the
group. In the minimally unstable case, we are able to show that these
invariants agree with the corresponding ones for the coroot integers and this
agreement is what allows us to identify the $\Cee^*$-weights with the coroot
integers. We emphasize, however, that circular symmetry holds for all
maximal parabolic subgroups, not just those which correspond to minimally
unstable bundles. Unfortunately, our proof of this resorts in the end to
case-by-case checking. It would be extremely illuminating to have a more
conceptual understanding of the meaning of circular symmetry. In
\S\ref{s1.5}, we discuss the volume of the moduli space of flat connections on
$E$ and again relate it to the coroot integers. These two pieces of numerical
information, the
$\Cee^*$-weights for
$\WP(\eta)$ and the volume of the space of flat connections, will turn out to be
crucial in Section 5 for the proof that $\deg \Psi =1$.

In Section 2, we give a general description, for every maximal parabolic
subgroup $P$, of the bundles
$\eta$ over the Levi factor $L$. We then compute the dimensions of the
nonabelian cohomology space (or rather its tangent space at the origin) and
the $\Cee^*$-weights in terms of the numbers introduced in Section 1. In
Section 3, we study the minimally unstable case in detail. As we have
mentioned above, we expect from general principles that $\dim \Aut_L\eta = 1$
and  that the cohomology dimension must be $\dim \mathcal{M}(G,c)+1$. We
identify the simple roots in the minimally unstable case, verifying the above
facts in a classification-free way in the simply connected case, and identify
the $\Cee^*$-weights with the coroot integers via circular symmetry. The
description of the minimally unstable case is also given in \cite{FMAB},
although the discussion here in the simply connected case is independent of
that paper. In \S\ref{ns}, we consider the non-simply
connected case. Here we use the results of \cite{FMAB} as well as a
case-by-case analysis to identify the minimally unstable bundles and to
identify the cohomology dimension, the dimension  of the automorphism group
of the bundle, and the $\Cee^*$-weights.

Section 4 is concerned with the nonabelian cohomology space, i.e.\ the affine
space $H^1(E; U(\eta))$. We show that the $\Cee^*$-action can be linearized
and discuss the obstructions to the existence of a universal bundle. In
\S\ref{s4.3}, we show that, in the minimally unstable case, the points of
$H^1(E; U(\eta)) -\{0\}$ correspond to \textbf{regular} semistable
bundles. This
means that the algebraic families provided by the parabolic
construction are different from the families provided by the space of
flat connections. Every S-equivalence class of bundles has two extreme
representatives which are unique up to isomorphism, the flat
representative and the regular representative. On an open dense set of
$\mathcal{M}(G,c)$ these representatives agree and, when they do, all bundles
of the given S-equivalence class are isomorphic. But, along a
codimension one subvariety of $\mathcal{M}(G,c)$, the regular representative does
not have a flat connection. It turns out that, because the dimension of the
automorphism group of regular representatives is constant, the regular
representatives behave better in families.

Finally, in Section 5, we prove the main theorem by calculating the degree
of $\Psi$. The crux of the argument is to study the determinant line bundle
on $\mathcal{M}(G,c)$, which pulls back via $\Psi^*$ to the determinant line bundle
on
$\WP(\eta)$. Thus, once we show that both determinant line bundles have the same
top self-intersection, then it follows that the degree of $\Psi$ is one.
Some of the technical results concerning the nonabelian cohomology
space and its interpretation are deferred to the appendix. We prove that
the cohomology space naturally has the structure of an affine space and
represents an appropriate functor.

The parabolic construction of semistable $G$-bundles was originally introduced
and explained, for the simply connected case, in \cite{FMW}, in a paper
written for an audience of physicists, as well as in the announcement
\cite{FMW1}. It is a pleasure to thank Ed Witten for originally
raising  the questions which led to this work and for the insights he
shared with us during the course of our joint work on these subjects. We
would also like to thank A. Borel, P. Deligne, and W. Schmid for various
helpful conversations and correspondence. Finally, during the preparation of
this paper, S. Helmke and P. Slodowy sent us their preprint \cite{HS}, which
has a considerable overlap with the first part of this paper and which also
analyzes the case where the automorphism group of $\eta\times_LG$ is just
slightly larger than in the minimally unstable case.

\section{Preliminaries}

\subsection{Notation}

Throughout this paper, $R$ denotes a reduced and irreducible root
system of rank $r$ in a real vector space $V$, with Weyl group $W=W(R)$, and
$\Delta$ is a set of simple roots for
$R$. Let
$R^+$ be the set of positive roots corresponding to the choice of
$\Delta$. There exists a
$W$-invariant positive definite inner product $\langle \cdot, \cdot \rangle$ on $V$
and it is unique up to scalars. Given a root
$\alpha$, there is an associated coroot $\alpha \spcheck\in V^*$.
Using the inner product to identify $V$ with $V^*$, we have $\alpha
\spcheck = 2\alpha/\langle \alpha, \alpha\rangle$. As usual, we denote the
Cartan integer  $\alpha( \beta \spcheck)$ by $n(\alpha,
\beta)$. The
{\sl coroot lattice\/}
$\Lambda$ is the lattice inside
$V^*$ spanned by the coroots. Given
  $\alpha \in \Delta$, we have the {\sl fundamental weight\/}
$\varpi _{\alpha}\in V$, which satisfies $\varpi _{\alpha}(\beta \spcheck) =
\delta_{\alpha\beta}$ for all $\alpha, \beta \in \Delta$.  The fundamental
coweights
$\varpi_{\alpha}\spcheck\in V^*$ are defined similarly. As usual, let
$\rho$ be the
sum of the fundamental weights, so that
$\rho =\sum _{\alpha \in \Delta}\varpi _{\alpha} = \frac12\sum _{\beta \in
R^+}\beta$.

Let $\widetilde \alpha$ be the highest root of $R^+$. We have
$\widetilde \alpha =\sum _{\beta \in \Delta}h_\beta\beta$, with
$h_\beta > 0$.
We set $\alpha_0 =-\widetilde \alpha$, $\widetilde
\Delta =\Delta \cup\{\alpha_0\}$, and
$h_{\alpha_0} =1$ so that $\sum _{\beta \in \widetilde
  \Delta}h_\beta\beta =0$. The
number
$h = 1+\sum _{\beta \in \Delta}h_\beta$ is the {\sl Coxeter number} of
$R$. Similarly, we have
$\widetilde\alpha\spcheck = \sum _{\beta \in
  \Delta}g_\beta\beta\spcheck$ with $g_\beta>0$. We set
$g_{\alpha_0}=1$ so that $\sum_{\beta \in \widetilde
  \Delta}g_\beta\beta\spcheck =0$. We call
$g = 1+ \sum _{\beta \in \Delta}g_\beta = \sum _{\beta \in
\widetilde\Delta}g_\beta$
the {\sl dual Coxeter number} of $R$. An easy calculation then shows:

\begin{lemma} In the above notation, we have
$$g_\alpha = \frac{h_\alpha\langle \alpha , \alpha \rangle}{\langle
\widetilde\alpha,
\widetilde\alpha\rangle}.$$ Thus $g_\alpha| h_\alpha$, and $g_\alpha= h_\alpha$ if
and only if
$\alpha$ is a long root of $R$.
\endproof
\end{lemma}

Let $Q\in \Sym^2\Lambda^*$ be the quadratic form  defined by
$$Q= \sum _{\alpha \in R}\langle \alpha, \cdot\rangle \langle \alpha,
\cdot\rangle.$$ This form was first introduced by  Looijenga in~\cite{Loo}, where
he showed:

\begin{lemma}\label{Looform} Let $g$ be the dual Coxeter number. Then
$$Q = (2g)I_0,$$ where $I_0$ is the unique $W$-invariant quadratic form on
$\Lambda$ such that
$I_0(\widetilde \alpha \spcheck) =2$. For example, if $R$ is
simply laced, then $I_0$ is the usual intersection form on $\Lambda$.
\endproof
\end{lemma}

Throughout this paper, we use the inner product on $V^*$ defined by
$I_0$ and the
corresponding dual inner product on $V$. It has the property that all long
roots have length $2$.

\begin{lemma}\label{next} Let
$\{\alpha^* \}_{\alpha \in \Delta}$ be the dual basis to $\{\alpha \spcheck
\}_{\alpha \in \Delta}$ with respect to $I_0$. Then
$\alpha^* =g_{\alpha}\varpi_{\alpha}\spcheck/h_{\alpha}$.
\end{lemma}
\begin{proof} This is just the statement that
$\displaystyle\langle \alpha  \spcheck, \varpi_\beta\spcheck\rangle =
\frac{2}{\langle
\alpha
 , \alpha \rangle}\alpha( \varpi_\beta\spcheck)=
\frac{h_{\alpha}}{g_{\alpha}}\delta_{\alpha\beta}$.
\end{proof}

\subsection{Structure of maximal parabolic subgroups}\label{parabolic}

Let
$H\subseteq G$ be a Cartan subgroup, and let $R\subseteq \frak h^*$ be the
set of roots for the pair
$(G,H)$. We denote by
${\frak g}$ the Lie algebra of $G$ and by ${\frak h}\subseteq {\frak g}$ the
Cartan subalgebra of $\frak g$ corresponding to $H$. Let $R^+$ be a set of
positive roots, and let $\Delta$ denote
the corresponding set of simple roots. We shall also view the roots
$\alpha\in R$ as characters
$\alpha\colon H\to  \Cee^*$ on $H$.  There is a
unique dual  coroot
$\alpha\spcheck\in {\frak h}$.
We can view $\alpha\spcheck$ as defining a linear map $\Cee\to {\frak h}$.
Exponentiating this map gives us a cocharacter of $H$, i.e. a one-parameter
subgroup $\ell_\alpha\colon \Cee^*\to H$.

If $P$ is a parabolic subgroup of $G$, then the unipotent radical $U$ of $P$ is a
normal subgroup, and there is a semidirect product $P= U\cdot L$, where $L$ is a
reductive subgroup of $P$, unique up to conjugation, called the Levi factor of
$P$. If $\widetilde G$ is the universal cover of $G$, then there is a
one-to-one correspondence between parabolic subgroups of $G$ and those of
$\widetilde G$, which associates to the subgroup $P$ of $G$ its preimage
$\widetilde P\subseteq \widetilde G$. Since a unipotent group is torsion-free, the
unipotent radicals of
$P$ and $\widetilde P$ are isomorphic, and the Levi factor $\widetilde L$ of
$\widetilde P$ is the preimage of the Levi factor $L$ of $P$. {\bf For the
remainder of \S\ref{parabolic}, unless otherwise stated, we assume that
$G$ is simply connected.}

For $\alpha \in \Delta$, let $P^\alpha$ be the connected subgroup of $G$ whose Lie
algebra is spanned by $\frak h$ and the root spaces $\frak g^\beta$, where either
$\beta \in R^+$ or $\beta$ lies in the span of $\Delta -\{\alpha\}$. Then
$P^\alpha$ is a maximal parabolic subgroup of $G$, and every maximal parabolic
subgroup $P$ is conjugate to exactly one $P^\alpha$, $\alpha \in
\Delta$.   Thus there are exactly $r$ maximal parabolic subgroups of
$G$ up to conjugation.
We denote the unipotent radical of $P^\alpha$ by $U^\alpha$ and its
Levi factor by $L^\alpha$.
The torus $H$ is a maximal torus of $L^{\alpha}$.
The semisimple part $S^{\alpha}$ of
$L^{\alpha}$ (or equivalently the derived subgroup of $L^{\alpha}$) has
Lie algebra spanned by ${\frak h}'$ and by the root spaces
corresponding to the  set of roots in the linear span of
$\Delta-\{\alpha\}$, where ${\frak h}'\subseteq {\frak h}$ is the subspace
spanned by the coroots $\beta\spcheck\in {\frak h}$ dual to the simple roots
$\beta\in\Delta-\{\alpha\}$. A maximal torus $H'$ of
$S^{\alpha}$ is given by  the subtorus which is the image under
exponentiation of ${\frak h}'\subseteq {\frak h}$, and $H' =H\cap
S^{\alpha}$. The Dynkin diagram of
$S^{\alpha}$ is the  subdiagram of the Dynkin diagram of $G$ obtained
by deleting
the vertex corresponding to
$\alpha$ and the edges incident to this vertex.
Since $\{\alpha\spcheck\}_{\alpha \in \Delta}$ is
a basis for the coroot lattice $\Lambda$ of $G$, the intersection of
the coroot lattice
$\Lambda$ for $G$ with ${\frak h}'$ is exactly the coroot
lattice for $S^{\alpha}$.
Since $G$ is simply connected,   $S^{\alpha}$ is also simply
connected.
Note that $S^{\alpha}$ is a semisimple group of rank $r-1$. Clearly we have:

\begin{lemma} Let $\Delta -\{\alpha\} = \coprod_{i=1}^t\Delta_i$, where each
$\Delta_i$ defines a connected component of the Dynkin diagram of $\Delta
-\{\alpha\}$. Then $S^\alpha =
\prod_{i=1}^tS_i$, where $S_i$ is the simple and simply connected
group corresponding to $\Delta_i$. \endproof
\end{lemma}

\begin{defn}\label{zetaalpha} Let $K(\alpha)=\bigcap _{\beta \in
    \Delta-\{\alpha\}}\Ker \beta \subseteq
\Lambda$. Then $K(\alpha)$ is an infinite cyclic group. Let
$\zeta_\alpha=\sum_{\beta \in \Delta}m_\beta\beta\spcheck$ be the generator of 
$K(\alpha)$   such that  $m_\alpha >0$. It then follows that $m_\beta > 0$
for all $\beta \in \Delta$.  Define $n_\alpha=\alpha(\zeta_\alpha)$, so that
$\zeta_\alpha= n_\alpha\varpi_\alpha\spcheck$. Note that $m_\alpha =
\varpi_\alpha(\zeta_\alpha)$.
\end{defn}

\begin{lemma}\label{Lalpha}
Define the map
$\varphi_\alpha\colon \Cee^*\to H$ by
$$\varphi_\alpha(\lambda)= \prod_{\beta \in \Delta} \ell
_\beta\left(\lambda^{m_\beta}\right).$$ Then $\varphi_\alpha$ is an
isomorphism from
$\Cee^*$ to the identity component of the center of $L^{\alpha}$. Moreover,
$$L^{\alpha} = S^{\alpha}\times_{\Zee/m_\alpha\Zee}\Cee^*,$$
where  $1\in \Zee/m_\alpha\Zee$ maps to $e^{2\pi i/m_\alpha}\in
\Cee^*$ and to the central element $\prod _{i=1}^t\gamma_i^{-n(\beta_i,\alpha)} \in
S^{\alpha} =\prod_{i=1}^tS_i$, where $\beta_i\in \Delta_i$ is the
unique element for which $n(\beta_i,\alpha)\not=0$
and $\gamma _i=\exp(2\pi\sqrt{-1}\varpi_{\beta_i}\spcheck)$.
\end{lemma}
\begin{proof}
With $\varphi_\alpha$ defined as above, $\varphi_\alpha(\Cee^*)$ is in
the kernel
of all simple roots
$\beta$ distinct from $\alpha$, and thus, since $\zeta_\alpha$ is primitive,
$\varphi_\alpha$ is an embedding of $\Cee^*$ into the center of $L^{\alpha}$.
Also, since $S^{\alpha}\cap H = H'$, if $\lambda = e^{2\pi
 it}$, then $\varphi _\alpha(\lambda) \in S^{\alpha}$ if and only if
$\sum _{\beta \in \Delta}tm_\beta\beta\spcheck\equiv 0 \bmod \frak h'+\Lambda$, if
and only if
$m_\alpha t\in
\Zee$.  Thus
$\varphi_\alpha(\Cee^*)
\cap S^{\alpha}$ is the cyclic subgroup of order
$m_\alpha$ in $\Cee^*$, and so
$$L^{\alpha} = S^{\alpha}\times_{\Zee/m_\alpha\Zee}\Cee^*,$$
where the image of $1\in \Zee/m_\alpha\Zee$ in the first factor lies in the center
of
$S^{\alpha}$, and corresponds to the element
$c=\zeta_\alpha/m_\alpha-\alpha\spcheck$. To describe this central element, let
$\beta\in
\Delta$ be a root of $S_i$. Then $\beta(c) = 0$ if $\beta \neq
\beta_i$ and $\beta_i(c) = -n(\beta_i,\alpha)$. Thus $c$ is the
central element of $S_i$ given by
$\exp(-2n(\beta_i,\alpha)\pi \sqrt{-1}\varpi_{\beta_i}\spcheck)$.
\end{proof}

The next lemma gives a more precise description of the center of $L^{\alpha}$.

\begin{lemma}
The center of $L^{\alpha}$ is generated by
$\varphi_\alpha\left(\Cee^*\right)$ and
the center of $G$. The intersection of
$\varphi_\alpha\left(\Cee^*\right)$ with
the center of $G$ is a cyclic group of order $n_\alpha$.
\end{lemma}
\begin{proof}
We have
$$Z(L^{\alpha})=Z(S^{\alpha})\times_{\Zee/m_\alpha\Zee}\Cee^*.$$
Then  $Z(G)$ is the subgroup of $Z(L^{\alpha})$ which
is in the kernel of the character $\alpha$.
The restriction of $\alpha$ to $\Cee^*\subseteq Z(L^{\alpha})$
is non-trivial, and hence surjective, and for an element of
$Z(S^{\alpha})\times_{\Zee/m_\alpha\Zee}\Cee^*$, written as $[z, \zeta]$, we
clearly have $\alpha ([z, \zeta]) =\alpha (z)\alpha(\zeta)$. Thus, for each
$z\in Z(S^{\alpha})$ there is an element
$u\in {\Ker}(\alpha)$ of the form $u=[z,\zeta]\in
Z(S^{\alpha})\times_{\Zee/m_\alpha\Zee}\Cee^*$. This element
$u$ is in the center of $G$, since it is in the kernel of $\beta$ for all
$\beta\neq \alpha$, as well as in $\Ker \alpha$. It follows that an arbitrary
$[z,
\lambda]\in Z(L^{\alpha})$ is of the form $u\cdot \mu$, where $\mu \in \Cee^*$, as
claimed.

To see the second statement, note that $\alpha \circ
\varphi_\alpha(\lambda) = \lambda ^{n_\alpha}$. Since
$\varphi_\alpha(\lambda)$ is in the kernel of all of
the remaining roots,   $\varphi_\alpha(\lambda)$ lies in the center of
$G$ if and only if $\lambda^{n_\alpha} = 1$.
\end{proof}

\begin{remark} Along similar lines, one can show that there is an exact sequence
$$\{1\} \to \Zee/n_\alpha\Zee \to Z(G) \to Z(S^\alpha)/(\Zee/m_\alpha\Zee) \to
\{1\}.$$
In particular, $m_\alpha/n_\alpha = \#Z(S^\alpha)/\#Z(G)$.
\end{remark}

To state the next result, recall that a character $\chi\colon H \to
\Cee^*$ is {\sl dominant\/}  if the corresponding linear function
$\frak h \to \Cee$ is nonnegative
on the positive coroots. A character of any subgroup containing $H$
will be called dominant if its restriction to $H$ is dominant.

\begin{lemma}\label{dominant}
The group of characters of $L^{\alpha}$ is isomorphic to $\Zee$. There
  is a   unique primitive
dominant character  $\chi_0\colon L^{\alpha}\to \Cee^*$, and $\chi_0\circ
\varphi_\alpha (\lambda) = \lambda^{m_\alpha}$.
\end{lemma}
\begin{proof} Let $\varpi_\alpha$ be the fundamental weight
  corresponding to the
simple root $\alpha$. The unique primitive dominant character of
$L^{\alpha}$ is
$\varpi_\alpha$, viewed as a character on $H$, and $\varpi_\alpha\circ
\varphi_\alpha (\lambda) = \lambda^{m_\alpha}$.
\end{proof}

If $G$ is not simply connected, we continue to denote the parabolic subgroup 
associated to the simple root $\alpha$ by $P^\alpha$, the unipotent radical of
$P^\alpha$ by 
$U^\alpha$,  and the  Levi factor of $P^\alpha$ by $L^\alpha$.  
The map $\varphi_\alpha\colon \Cee^*\to L^\alpha$ is defined to be
embedding of $\Cee^*$ into to the center of $L^\alpha$ so that the
composition of $\varphi_\alpha$ followed by the primitive dominant
character of $L^\alpha$ is a positive character on $\Cee^*$.
 Let $\widetilde
\varphi_\alpha\colon \Cee^*\to \widetilde L^\alpha$ be as given
in  Definition~\ref{zetaalpha} for the simply connected  group
$\widetilde G$ and $\alpha\in \Delta$. Clearly,  $\varphi _\alpha(\Cee^*)$ is
the quotient of $\widetilde
\varphi_\alpha( \Cee^*)$ by the finite subgroup $\widetilde
\varphi_\alpha( \Cee^*)\cap \langle c \rangle$.  It is still the case, of course,
that the center of
$L^\alpha$ is generated by the center of $G$ and $\varphi_\alpha(\Cee^*)$. The
following lemma is then clear:

\begin{lemma}\label{ncalpha} Let $n_{c, \alpha}$ be the order of 
$\varphi_\alpha(\Cee^*)\cap Z(G)$. Then the order of $\widetilde
\varphi_\alpha( \Cee^*)\cap \langle c \rangle$ is $n_\alpha/n_{c, \alpha}$, and
the induced map from $\widetilde
\varphi_\alpha( \Cee^*) \cong \Cee^*$ to $\varphi _\alpha(\Cee^*) \cong \Cee^*$ is
given by raising to the power $n_\alpha/n_{c, \alpha}$. \endproof
\end{lemma}

\subsection{The unipotent radical of a maximal parabolic
  subalgebra}\label{para2}

{\bf In \S\ref{para2}, we assume that $G$ is simply connected.}
Fix the simple root $\alpha$ and consider the maximal parabolic
subgroup  $P^{\alpha}$. Here we will describe
the Lie algebra
$\frak u = \frak u (\alpha)$ of
$U^{\alpha}$. It is spanned by the root spaces $\frak g^\delta$   such
that, if $\delta =
\sum _{\beta\in \Delta}x_\beta\beta$, then the  coefficient of $\alpha$ in the
sum is positive.
The action of
$\varphi_\alpha(\Cee^*)$  on $\frak u $ is given as follows.
Let $\frak g^\delta$ be a  root space and let $X\in \frak g^\delta$.
Then
$\varphi_\alpha(\lambda)(X) = \lambda^{\delta(\zeta_\alpha)}\cdot X$.
Of course, since $\zeta_\alpha$ is in the kernel of all the simple roots
except $\alpha$,
$\delta(\zeta_\alpha)=\delta(\varpi_\alpha\spcheck)\alpha(\zeta_\alpha)$ where
$\delta(\varpi_\alpha\spcheck)$ is the coefficient of
$\alpha$ in the expression for $\delta$ as a sum of simple roots. In
particular, using Lemma~\ref{unbroken} below, the weights of the action of $\Cee^*$
on
$\frak u$ are
$n_\alpha,2n_\alpha,\ldots,h_\alpha\cdot n_\alpha$ where 
$n_\alpha=\alpha(\zeta_\alpha)$.

Let $\frak u ^k$ be the sum of all such
root spaces $\frak g^\beta$ where the coefficient of $\alpha$ in
$\beta$ is exactly $k$. The Lie
algebra $\frak u $ is a direct sum of the spaces $\frak u^k$ for
$k>0$.  Each space
$\frak u^k$ is an $L^\alpha$-module, and as we shall see below it is in fact
irreducible.

The top exterior power  $\bigwedge^{\text{\rm top}}{\frak u}$ is a
one-dimensional
$L^\alpha$-module, and as such it is given by a character $\chi_+$ of
$L$. The next
lemma identifies this character:

\begin{lemma}\label{defd1} Let $d_1(\alpha) = 2\rho(\zeta_\alpha)/m_\alpha =
2\rho(\varpi_\alpha\spcheck)/\varpi_\alpha(\varpi_\alpha\spcheck)$. In
the above notation, $\chi_+ = \chi_0^{d_1(\alpha)}$ where $\chi_0$ is
the primitive dominant character.
\end{lemma}
\begin{proof} Since $\chi_0$ is primitive
$\chi_+ = \chi_0^{N_+}$ for some integer $N_+$.
Since the character
$\chi_0\circ\varphi_{\alpha}$ of $\Cee^*$ is  given by $\lambda
\mapsto \lambda^{m_\alpha}$, the character
$\chi_+\circ\varphi_\alpha$ of $\Cee^*$ is given by $\lambda
\mapsto
\lambda^{m_\alpha\cdot N_+}$.
The action of
$\varphi_\alpha(\Cee^*)$ is diagonal with respect to the
decomposition of ${\frak u}$ as a sum of root spaces,
and  the character on the one-dimensional subspace spanned
by a root
$\delta$ is simply the restriction of $\delta$ to
$\varphi_\alpha(\Cee^*)$. The space ${\frak u}$  is the
subspace of ${\frak g}$ spanned by the  set of  all roots
$\delta = \sum _{\beta \in \Delta}x_\beta\beta$ with the property that $x_\alpha
>0$. Recall that
$\varphi_\alpha(\Cee^*)$ is in the kernel of all  the
simple roots except $\alpha$. Thus, to compute the
character of $\Cee^*$ given by the product of $\varphi
_\alpha^\delta$ over all $\delta$ such that the coefficient of
$\alpha$ in $\delta$ is positive, we may as well take the
product over all of the positive roots $\delta$. In other
words, the character of $\Cee^*$ which gives the degree of the top
exterior power is simply the character
$\varphi_\alpha^{\sum_{\delta\in R^+}\delta}$.  We can rewrite this
expression as $\varphi_\alpha^{ 2\rho}$.   Thus, the character
that we are
computing is
$\varphi_\alpha^{2\sum_{\beta \in \Delta}\varpi_\beta}$. Recalling that
the embedding of $\varphi_\alpha\colon\Cee^*\to H$ is given
by
$\prod_{\beta \in \Delta}\ell_\beta^{m_\beta} $,   we
see that the character
$\chi_+\circ \varphi_\alpha$ is given by raising to the power
$2\rho(\zeta_\alpha)$. Thus $N_+ =
2\rho(\zeta_\alpha)/m_\alpha=d_1(\alpha)$.
\end{proof}

\begin{remark} The proof above also shows that the integer
  $d_1(\alpha)$ is the degree of divisibility of the canonical bundle
  of the homogeneous space $G/P^{\alpha}$.
\end{remark}

Now we  give a purely root theoretic formula for $d_1(\alpha)$.

\begin{lemma}\label{d1lemma} $d_1(\alpha) = \sum _  { \beta
(\varpi_\alpha\spcheck) >0}  n(\beta,
\alpha)$, where the $\beta$ in the sum range  over $R$.
\end{lemma}
\begin{proof} By definition, $\chi_+= \sum_{\beta (\varpi_\alpha\spcheck)
>0}\beta$ as additive  characters on $\frak h$. By the previous lemma,
$\chi_+(\zeta_\alpha/m_\alpha) = d_1(\alpha)\chi_0(\zeta_\alpha/m_\alpha) =
d_1(\alpha)$. On the other hand, $\zeta_\alpha/m_\alpha =
\alpha\spcheck + \nu$,
where $\nu$ lies in the $\Q$-span of the simple coroots
$\beta\spcheck$, $\beta\neq
\alpha$, and hence in the Lie algebra of the derived group $S^\alpha$. Thus
$$\chi_+(\zeta_\alpha/m_\alpha) = \chi_+(\alpha\spcheck)=
\sum_{\substack {\beta\in
R\\
 \beta (\varpi_\alpha\spcheck) >0}}\beta(\alpha\spcheck),$$
and so $d_1(\alpha) = \sum _{\beta (\varpi_\alpha\spcheck) >0}n(\beta,
\alpha)$.
\end{proof}

Similarly, we compute the character $\chi_k$ of $L^\alpha$ corresponding to
$\bigwedge^{\text{\rm top}}{\frak u}^k$:

\begin{lemma}\label{deguk}  For $k>0$, let $c(\alpha, k)$ be the dimension of
$\frak u^k$, in other words the number of roots $\beta$ such that the
coefficient of $\alpha$ in
$\beta$ is $k$. Let $i(\alpha, k) = kc(\alpha,
k)/\varpi_\alpha(\varpi_\alpha\spcheck)= kn_\alpha c(\alpha,
k)/m_\alpha$. Then the character $\chi_k$ of
$L^\alpha$ corresponding to
$\bigwedge^{\text{\rm top}}{\frak u}^k$ is given by
$$\chi_k = \chi_0^{i(\alpha, k)}.$$
\end{lemma}
\begin{proof} As before, $\chi_k = \chi_0^{N_k}$ for some integer
  $N_k$. To compute
$N_k$, note that, if $\frak g^\beta \subseteq \frak u^k$, then
$\varphi_\alpha$ acts on
$\frak g^\beta$ via raising to the power $\beta(\zeta_\alpha) =
kn_\alpha$. Since $\chi_0\circ
\varphi_\alpha (\lambda) =
\lambda^{m_\alpha}$, we must have
$$m_\alpha N_k = kn_\alpha\dim {\frak u}^k.$$
Hence, $N_k =  kn_\alpha
c(\alpha,k)/m_\alpha = i(\alpha, k)$.
\end{proof}

\begin{remark} An argument very similar to the proof of
  Lemma~\ref{d1lemma} shows that
$$i(\alpha, k) = \sum_{ \beta(\varpi_\alpha\spcheck)
=  k}n(\beta,\alpha).$$
\end{remark}

\subsection{Some lemmas on root systems}

Fix $\alpha \in \Delta$. Our goal now is to analyze the
$L^\alpha$-representations ${\frak u}^k$, and in particular the numbers
$d_1(\alpha)$ and $i(\alpha, k)$  introduced above.

\begin{defn} Fix a positive integer $k$, and
consider the set
$$S(\alpha ,k)= \{\beta \in R: \beta(\varpi_\alpha\spcheck) =k\}.$$
For $k=0$, we define similarly
$$S^+(\alpha ,0)= \{\beta \in R^+: \beta(\varpi_\alpha\spcheck) =0\}.$$
The latter is a set of positive roots for the root system $R'=\Delta
-\{\alpha\}$, i.e.
for $S^\alpha$. We define $S^-(\alpha ,0)$ similarly.
  A {\sl lowest root\/}
$\sigma_k(\alpha)$ for $S(\alpha ,k)$ is a root $\sigma_k(\alpha) \in
S(\alpha ,k)$ such that, for every $\beta \in S(\alpha ,k)$, $\beta -
\sigma_k(\alpha)$ is a sum (possibly empty) of simple roots. For example,
$\sigma_1(\alpha) =\alpha$. A lowest root in
$S(\alpha ,k)$ is clearly unique, if it exists. A highest root
$\lambda_k(\alpha)\in S(\alpha ,k)$ is defined similarly, and is also clearly
unique if it exists.
\end{defn}

The following is related to results of Borel-Tits
(unpublished) as well as Azad-Barry-Seitz~\cite{ABS}.

\begin{proposition}\label{highest} In the above notation, if $S(\alpha
  ,k) \neq
\emptyset$, then lowest roots and highest roots always exist and are
unique.
\end{proposition}
\begin{proof} Let $R(\alpha ,k)$ be the subset of $R$ consisting of
  roots $\beta$
such that $k$ divides $\beta(\varpi_\alpha\spcheck)$.  Clearly
$R(\alpha ,k)$ is again a root system. Let $V$ be the real span of $R$
and let $V'$ be the subspace of $V$  spanned by $\Delta
-\{\alpha\}$. Then clearly $V'\cap R = V'\cap R(\alpha ,k) = R'$ is
the set of all roots which are linear combinations of
elements of $\Delta -\{\alpha\}$. Thus $R'$ is a root system with simple roots
$\Delta -\{\alpha\}$. By~\cite[VI \S 1, Proposition 24]{Bour}, since
$R'\subseteq R(\alpha, k)$, there exists a set of simple roots for
$R(\alpha ,k)$ containing
$\Delta -\{\alpha\}$. In fact, the proof of this proposition shows that there
are at least two different sets of simple roots, each of the form  $(\Delta
-\{\alpha\})\cup
\{\beta\}$. Then
$k|\beta(\varpi_\alpha\spcheck)$, and since $S(\alpha ,k) \neq
\emptyset$, in fact
$\beta(\varpi_\alpha\spcheck) =
\pm k$. Suppose for example that $\beta(\varpi_\alpha\spcheck) = k$. If
$\gamma \in S(\alpha ,k)$, then   $\gamma = \sum
_{\delta\neq \alpha}m_\delta\delta + m\beta$ and
$\gamma(\varpi_\alpha\spcheck)=\beta(\varpi_\alpha\spcheck) = k$. Thus
$m =1$ and
$m_\delta\geq 0$ for all
$\delta\in \Delta-\{\alpha\}$, so that $\beta$ is a lowest root for
$S(\alpha ,k)$.
Suppose now that
$(\Delta -\{\alpha\})\cup \{\beta'\}$ is also a set of simple roots
for $R(\alpha ,k)$, where
$\beta'\neq \beta$. It follows that $\varpi_\alpha\spcheck(\beta') =
-k$. In this case,
it is easy to check that $-\beta'$ is a highest root for $S(\alpha
,k)$.
\end{proof}

\begin{corollary} The $L^\alpha$-modules $\frak u^k$ are irreducible.
\endproof
\end{corollary}

\begin{remark}  Using the Borel-de Siebenthal theorem~\cite[p.\ 229,
  ex.\ 4]{Bour}, the Dynkin
diagram of
$R(\alpha, k)$ is given abstractly as follows. Begin with the extended Dynkin
diagram of $R$. There exists a root $\beta$ such that $h_\beta = k$,
and such that
the Dynkin diagram for $\widetilde \Delta -\{\beta\}$ contains a root $\gamma$
such that the diagram for $\widetilde \Delta -\{\beta, \gamma\}$ is the same as
the diagram for $\Delta -\{\alpha\}$. In practice, these properties uniquely
determine $\beta$ and $\gamma$. The Dynkin diagram for $R(\alpha, k)$
is then the
Dynkin diagram of $\widetilde \Delta -\{\beta\}$, and $\gamma$ corresponds to
$\lambda_k(\alpha)$.
\end{remark}

\begin{lemma}\label{unbroken}
Let $k$ be a positive integer.
Then $S(\alpha,k)\neq \emptyset$ if and only if $1\leq k\leq h_\alpha$.
\end{lemma}

\begin{proof}
Since $\alpha \in S(\alpha, 1)$, 
$S(\alpha,1)\neq \emptyset$. Suppose inductively we have shown that
$S(\alpha,k)\neq \emptyset$ for all $k,1\le k\le \ell$. If $\beta\in
\Delta-\alpha$, then 
$\lambda_\ell(\alpha)+\beta$ is not a root. If
$\lambda_\ell(\alpha)+\alpha$ is not a root, then
$\lambda_\ell(\alpha)$ is the highest root and
$\ell=h_\alpha$. Otherwise, $\lambda_\ell(\alpha)+\alpha\in
S(\alpha,\ell+1)$.
\end{proof}

For future reference, we record the following properties of
$\sigma_k(\alpha)$ and $\lambda_k(\alpha)$:

\begin{lemma}\label{lambdaprops} In the above notation, let
  $R'\subseteq R$ be the
subroot system with simple roots $\Delta -\{\alpha\}$. Suppose that
$S(\alpha,k) \neq \emptyset$, and let
$w_0'\in W(R')\subseteq W(R)$ be the unique element such that $w_0'(\Delta
-\{\alpha\}) = -(\Delta -\{\alpha\})$. Suppose that $\tau $ is the permutation
of $\Delta -\{\alpha\}$ induced by $-w_0'$, which we can also view as
a permutation of
$\Delta$ fixing $\alpha$. Let $\tau$ act on $V$ and $V^*$ in the natural
way.  Then:
\begin{enumerate}
\item[\rm (i)] $w_0'\sigma_k(\alpha) = \lambda_k(\alpha)$;
\item[\rm (ii)] $\lambda _1(\alpha\spcheck) = (\lambda_1(\alpha))\spcheck$;
\item[\rm (iii)] $\sigma_k(\alpha) = k\lambda_1(\alpha)-\tau(\lambda_k(\alpha))
+k\alpha$. Likewise $\sigma_k(\alpha)\spcheck =
k'\lambda_1(\alpha\spcheck)-\tau(\lambda_k(\alpha)\spcheck) +k'\alpha\spcheck$,
where $k'$ is the coefficient of $\alpha\spcheck$ in
$\lambda_k(\alpha)\spcheck$.
\end{enumerate}
\end{lemma}
\begin{proof} Clearly $w_0'\sigma_k(\alpha) \in S(\alpha,k)$ has the property
that $(\Delta -\{\alpha\}) \cup \{-w_0'\sigma_k(\alpha)\}$ is a set of simple
roots for $R(\alpha,k)$. By the proof of Proposition~\ref{highest},
$w_0'\sigma_k(\alpha) = \lambda_k(\alpha)$, proving (i). To see (ii), note that
$\alpha\spcheck = \sigma _1(\alpha\spcheck)$, and thus
$\lambda _1(\alpha\spcheck) = w_0'\sigma _1(\alpha\spcheck) =
w_0'(\alpha\spcheck)
= (w_0'\alpha)\spcheck = (\lambda_1(\alpha))\spcheck$. To see (iii), write
$\lambda_k(\alpha) = k\alpha + \lambda '$, where $\lambda' \in V'$. Then
$$\sigma _k(\alpha) = w_0'\lambda_k(\alpha) = kw_0'\alpha - \tau (\lambda') =
k\lambda_1(\alpha) - \tau (\lambda').$$
On the other hand, $-\tau (\lambda') = k\alpha-\tau(\lambda_k(\alpha))$, and
plugging this back in gives the first part of (iii). The second part is proved
in a very similar way.
\end{proof}

Recall that, for $k> 0$,   $c(\alpha ,k)$ is the cardinality of
$S(\alpha ,k)$, in other words the number of $\beta \in R$ such that
$\beta(\varpi_\alpha\spcheck) = k$, and that $i(\alpha, k)  =
kc(\alpha, k)/\varpi
_{\alpha}( \varpi_{\alpha}^\vee )$.

\begin{defn}\label{defdk} Define
$$d_k(\alpha) =\sum _{k|x}i(\alpha, x)= \sum _{\ell > 0}i(\alpha, \ell   k).$$
In particular $d_1(\alpha)$ agrees with the previous definition.
Of course, the $d_k(\alpha)$ determine $i(\alpha, k)$
via Moebius inversion. Since the $i(\alpha, k)$  are all positive
integers, the $d_k(\alpha)$ are positive integers as well.
\end{defn}

\begin{lemma}\label{formulas} With notation as above,
\begin{align*}
i(\alpha, k) &=  \frac{h_\alpha g kc(\alpha, k)}{g_\alpha\sum
_{\beta
\in R^+}(\beta(\varpi_\alpha\spcheck))^2};\\
d_1(\alpha) &=  2 \rho(\varpi_\alpha\spcheck)/\varpi
_{\alpha}( \varpi_{\alpha}^\vee ) = \frac{h_\alpha g\sum _{\beta \in
R^+}\beta(\varpi_\alpha\spcheck)}{g_\alpha \sum _{\beta \in
R^+}(\beta(\varpi_\alpha\spcheck))^2}.
\end{align*}
Moreover $\sum _{k>0}\varphi(k)d_k(\alpha) = {h_\alpha} g/g_\alpha$.
\end{lemma}
\begin{proof} The first equality follows since, by Lemma~\ref{Looform} and
Lemma~\ref{next},
$$2g\left(\frac{h_\alpha}{g_\alpha}\right)\varpi
_{\alpha}( \varpi_{\alpha}^\vee ) =2g\langle
\varpi _{\alpha}^\vee,
\varpi_{\alpha}^\vee \rangle = 2gI_0(\varpi _{\alpha}^\vee)  =2\sum _{\beta \in
R^+}(\beta(\varpi_\alpha\spcheck))^2.$$ The second follows since $\sum _{k>0}k
c(\alpha, k) =
\sum _{\beta \in R^+}\beta(\varpi_\alpha\spcheck)$. The third is an easy
consequence of the fact that $\sum _{k|x}\varphi (k) = x$, using $\sum
_{k>0}k^2c(\alpha, k) =
\sum _{\beta \in R^+}(\beta(\varpi_\alpha\spcheck))^2$.
\end{proof}

\begin{proposition}\label{d1} We have $d_1(\alpha) = 
\rho(\lambda_1(\alpha\spcheck)) + 1$, where $\lambda_1(\alpha\spcheck)$
is the highest coroot such that the coefficient of $\alpha\spcheck$ in
$\lambda_1(\alpha\spcheck)$ is $1$.
\end{proposition}
\begin{proof} By Lemma~\ref{d1lemma}, it suffices to prove that
$\rho(\lambda_1(\alpha\spcheck)) + 1 = \sum _{\beta
  (\varpi_\alpha\spcheck) >0}n(\beta,
\alpha)$. By Lemma~\ref{lambdaprops}, since $\sigma_1(\alpha\spcheck)
=\alpha\spcheck$,
$w_0'(\alpha\spcheck) = \lambda_1(\alpha\spcheck)$. Since $w_0'$ exchanges
positive and negative roots in $R'$, it follows that $-w_0'$ is a
  permutation of
$S^+(\alpha, 0)$. Clearly, given $\beta \in S^+(\alpha, 0)$,
$$\beta(\lambda_1(\alpha\spcheck)) = -(-w_0'(\beta)(\alpha\spcheck)).$$
Thus, since $-w_0'$ permutes
$S^+(\alpha, 0)$, we have
$$\sum _{ \substack{ \beta \in S^+(\alpha, 0)\\
    \beta(\lambda_1(\alpha\spcheck))
>0}}\beta(\lambda_1(\alpha\spcheck)) + \sum _{ \substack{ \beta \in
S^+(\alpha, 0)\\ \beta( \alpha\spcheck) <0}}\beta( \alpha\spcheck ) =
0.$$
Next we claim:

\begin{lemma}\label{sum1}
$$-\sum _{ \substack{ \beta \in S^+(\alpha,
0)\\ \beta( \alpha\spcheck) <0}}\beta( \alpha\spcheck )\  \ + \ \  2 = \sum
_{\beta (\varpi_\alpha\spcheck) >0}\beta(\alpha\spcheck )= \sum
_{\beta (\varpi_\alpha\spcheck) >0}n(\beta, \alpha).$$
\end{lemma}

\begin{proof} Suppose that $\beta \neq \pm \alpha\in R$, and consider the
$\alpha$-string defined by $\beta$, say $\beta - q\alpha, \dots, \beta
+ p\alpha$.
Since $p-q+1=-n(\beta, \alpha) +1$, it is easy to see that $\sum
_{i=-q}^p(\beta +
i\alpha)(\alpha\spcheck) = 0$. If
$\beta$ is a negative root, then every root in the $\alpha$-string is negative
and thus none of them appears in the right hand side of the above
equality.   After
reindexing, we can assume that
$\beta$ is the origin of the
$\alpha$-string. If $\beta( \alpha\spcheck)>0$, then $(\beta+i\alpha)(
\alpha\spcheck)>0$ for all $i>0$, and so the total contribution to the right
hand side from the sum over the $\alpha$-string is zero. If $\beta \in S^+(\alpha,
0)$, then the contribution
to the sum on the right hand side is $\sum _{i\geq 1}(\beta +
i\alpha)(\alpha\spcheck) = -\beta(\alpha\spcheck)$. The remaining
possibility for
the right hand side is $\beta =\alpha$, and in this case
$\beta(\alpha\spcheck) =
2$. Thus we see that the right hand side in Lemma~\ref{sum1} is equal
to the left
hand side.
\end{proof}

\begin{lemma}\label{sum2}
$$\sum _{ \substack{ \beta \in S^+(\alpha, 0)\\
    \beta(\lambda_1(\alpha\spcheck))
>0}}\beta(\lambda_1(\alpha\spcheck)) =   \rho(\lambda_1(\alpha\spcheck)) -
1.$$
\end{lemma}
\begin{proof} First note that
$$2 \rho( \lambda_1(\alpha\spcheck)) = \sum _{\beta \in
R^+}\beta(\lambda_1(\alpha\spcheck)).$$
We consider as before the
$\lambda_1(\alpha)$-strings defined by a root $\beta\neq
\lambda_1(\alpha)$ which
lie in $R^+$. If the origin of such a string lies in
$R^+$, then so does every $\gamma$ lying in the string, and the sum
over all such
$\gamma$ of $n(\gamma, \lambda_1(\alpha))$ is zero. Next we claim:

\begin{claim} If a $\lambda_1(\alpha)$-string meets $R^+$ but is not
  contained in $R^+$, then either:
\begin{enumerate}
\item[\rm (i)] The origin of the string lies in $S^-(\alpha, 0)$ and all other
elements of the string lie in $R^+$.
\item[\rm (ii)] The extremity of the string lies in $S^+(\alpha, 0)$, and all
other elements lie in $R^-$.
\end{enumerate}
Moreover, there is a length-preserving bijection between strings of types {\rm
(i)} and {\rm (ii)} above.
\end{claim}
\begin{proof} First note that, if $\beta \in S(\alpha, 1)$, then
$\beta -2\lambda_1(\alpha)$ cannot be a root. For then $\beta -
\lambda_1(\alpha)$
would also be a root, necessarily negative, and then
$2\lambda_1(\alpha) -\beta$
would be an element of $S(\alpha, 1)$ higher than
$\lambda_1(\alpha)$. Thus, every
$\lambda_1(\alpha)$-string meeting $R^+$ but not contained in it must either begin
or terminate in
$S^+(\alpha, 0)\cup S^-(\alpha, 0)$. Clearly, the only possibilities
are (i) and
(ii) above, and the bijection is  given by sending the origin of a
string of type
(i) to its negative, which is the extremity of a string of type (ii).
\end{proof}

Returning to the proof of Lemma~\ref{sum2}, the only nonzero
contributions to the
sum $\sum _{\beta \in
R^+}\beta(\lambda_1(\alpha\spcheck))$ come from
\begin{enumerate}
\item[\rm (i)] $\lambda_1(\alpha)$-strings whose origin is $\gamma =-\beta
\in S^-(\alpha,0)$, and these contribute $-\gamma(\lambda_1(\alpha\spcheck)) =
\beta(\lambda_1(\alpha\spcheck))$;
\item[\rm (ii)] $\lambda_1(\alpha)$-strings whose extremity is $\beta
\in S^+(\alpha,0)$, and these contribute $\beta(\lambda_1(\alpha\spcheck))$;
\item[\rm (iii)] the root $\lambda_1(\alpha)$ and this contributes
$\lambda_1(\alpha)(\lambda_1(\alpha\spcheck)) =2$.
\end{enumerate}
Summing these up, we see that
$$2 \rho(\lambda_1(\alpha\spcheck)) = \sum _{\beta \in
R^+}\beta(\lambda_1(\alpha\spcheck)) = 2\sum _{ \substack{ \beta \in
  S^+(\alpha,
0)\\ \beta(\lambda_1(\alpha\spcheck)) >0}}\beta(\lambda_1(\alpha\spcheck))\ \
+\ \  2.$$ Dividing by $2$ gives the final formula of Lemma~\ref{sum2}.
\end{proof}

To complete the proof of Proposition~\ref{d1}, we have
$$\sum _{\beta (\varpi_\alpha\spcheck) >0}n(\beta, \alpha)=  -\sum _{
  \substack{
\beta \in S^+(\alpha, 0)\\ \beta( \alpha\spcheck) <0}}\beta(
\alpha\spcheck )\
+\   2 =
\sum _{ \substack{ \beta \in S^+(\alpha, 0)\\
    \beta(\lambda_1(\alpha\spcheck))
>0}}\beta(\lambda_1(\alpha\spcheck))\ +\ 2 = 
\rho(\lambda_1(\alpha\spcheck)) + 1,$$ as claimed.
\end{proof}

There is a generalization of the previous proposition to the computation of
$d_k(\alpha)$ for every $k>0$:

\begin{proposition}
Let $k' =k\langle \alpha,
\alpha\rangle/\langle \lambda_k(\alpha), \lambda_k(\alpha)\rangle$ be the
coefficient of $\alpha\spcheck$ in
$\lambda_k(\alpha)\spcheck$. Then  $\displaystyle d_1(\alpha) + d_k(\alpha) =
\frac{2} {k'}( \rho(\lambda_k(\alpha)\spcheck) + 1)$.
\end{proposition}
\begin{proof} In the notation of Proposition~\ref{highest},
set $\alpha_k = -\lambda_k(\alpha)$ and let $R_k=R(\alpha ,k)$, with simple
roots
$(\Delta -\{\alpha\})\cup \{\alpha_k\}$. Although $R_k$ need not be irreducible,
we can still define the integer
$d_1^{R_k}(\alpha_k)$ with respect to the root system $R_k$ as in
Lemma~\ref{defd1}. By Proposition~\ref{d1}, which holds even if $R_k$ is reducible,
$d_1^{R_k}(\alpha_k) =
\rho_k(\lambda_1((\alpha_k)\spcheck)) + 1$, where $\rho_k$ is the sum of the
fundamental weights of the root system $R_k$. Applying Lemma~\ref{lambdaprops},
and using the notation introduced in its statement, we have
\begin{align*}\lambda_1((\alpha_k)\spcheck) &=  \lambda_1(\alpha_k)\spcheck =
w_0'(\alpha_k)\spcheck = (w_0'\alpha_k)\spcheck =
-(\sigma_k(\alpha))\spcheck\\
&=  -k'\lambda_1(\alpha\spcheck) +\tau(\lambda_k(\alpha))\spcheck)
-k'\alpha\spcheck,
\end{align*}
where $k' =k\langle \alpha,
\alpha\rangle/\langle \lambda_k(\alpha), \lambda_k(\alpha)\rangle$ is
the coefficient of $\alpha\spcheck$ in
$\lambda_k(\alpha)\spcheck$. Next we compute $\rho_k$. Write
$(\alpha_k)\spcheck =
\sum _{\beta \in \Delta -\{\alpha\}}c_\beta\beta\spcheck -
k'\alpha\spcheck = \sum
_{\beta \in \Delta}c_\beta\beta\spcheck$, where $c_\alpha = -
k'$. Denote the fundamental weights of
$R_k$ by $\varpi_\beta^k$,
$\beta \neq \alpha$, and $\varpi_{\alpha_k}$. Then it is easy to check that
\begin{eqnarray*} \varpi_{\alpha_k} &=& -\frac{1}{k'}\varpi_{\alpha};\\
\varpi_\beta^k &=& \varpi_\beta + \frac{c_\beta}{k'}\varpi_\alpha
\mbox{\qquad  for $\beta \neq \alpha$.}
\end{eqnarray*}
Thus
\begin{multline*}
\rho_k = \sum _{\beta \neq \alpha}\varpi_\beta^k + \varpi_{\alpha_k} =
\sum _{\beta
\neq \alpha}\varpi_\beta +\frac{1}{k'}\left(\sum _{\beta \neq \alpha} c_\beta
-1\right)\varpi_{\alpha} = \sum _\beta \varpi_\beta +\frac{1}{k'}\left(\sum
_{\beta } c_\beta -1\right)\varpi_{\alpha} \\=\rho
-\frac{1}{k'}\left( \rho(\lambda_k(\alpha)\spcheck)
+1\right)\varpi_{\alpha}.
\end{multline*}
Since $\varpi_\alpha( -k'\lambda_1(\alpha\spcheck) +
\tau(\lambda_k(\alpha))\spcheck) -k'\alpha\spcheck) =-k'$ and $ \rho(
\tau (\gamma)) =  \rho(\gamma)$ for all $\gamma$, we see that
\begin{align*}
d_1^{R_k}(\alpha_k) &=  -k'( \rho(\lambda_1(\alpha\spcheck) +1) +
 \rho(\lambda_k(\alpha)\spcheck)) +1 +
\rho(\lambda_k(\alpha)\spcheck)) +1\\
&= -k'd_1(\alpha) + 2( \rho(\lambda_k(\alpha)\spcheck)) +1).
\end{align*}
Now an argument similar to the calculation of $\rho_k$ above shows that
$\varpi_{\alpha_k}\spcheck$, the fundamental coweight for $R_k$ dual
to $\alpha
_k$, is given by $-(1/k)\varpi_\alpha\spcheck$. Also, $c(\alpha, nk) =
c(\alpha_k,
n)$ for all $n$. Thus $i(\alpha_k, n) = k'i(\alpha, nk)$ for all
positive integers
$n$. It follows that $d_1^{R_k}(\alpha_k) = k'd_k(\alpha)$. Putting
this together
with the above gives $k'd_1(\alpha) + k'd_k(\alpha) = 2(
\rho(\lambda_k(\alpha)\spcheck)) +1)$, which is the statement of the
proposition.
\end{proof}

\begin{corollary}\label{firstcase} $\displaystyle d_1(\alpha) +
d_{h_\alpha}(\alpha) =
\frac{2g}{g_\alpha}$.
\end{corollary}
\begin{proof}  For $k=h_\alpha$, we have $\lambda_k(\alpha)=\widetilde
  \alpha$ and
$k' = kg_\alpha/h_\alpha=  g_\alpha$, and the corollary is clear.
\end{proof}

To put the corollary in a more general context, we  have the following
definition
which is taken from~\cite{BFM}:

\begin{defn} For a sequence  $d_1, \dots, d_N$  of positive integers, let $M
= \sum _{k>0}\varphi (k)d_k$. Let the Farey sequence  ${\mathcal F}_N
=\{0/1, 1/N,
1/(N-1), \dots\}$ be the sequence of rational numbers between $0$ and $1$, written
in lowest terms, whose denominator is at most $N$, written in increasing order. We
say that the sequence $\{d_k\}$ has the {\sl circular symmetry property\/} with
respect to
$N$ and $M$ if,  for all  consecutive terms $r/x$ and
$s/y$ in ${\mathcal F}_N$,
$$d_x+d_y = \frac{2M }{xy} .$$
The geometric meaning of this property is explained in~\cite{BFM}.
\end{defn}

Since every integer $x$, $1\leq k\leq N$, appears as a denominator of some 
element of ${\mathcal F}_N$, the following is clear:

\begin{lemma}\label{circdetermines} Suppose that $d_1, \dots, d_N$ and $d_1',
\dots, d_N'$ are two sequences of positive integers which both satisfy the
circular symmetry property with respect to $N$ and $M$. If $d_1=d_1'$, then
$d_x=d_x'$ for all $x$, $1\leq k\leq N$.
\endproof
\end{lemma}

For the sequence $\{d_k(\alpha)\}$, the largest integer $N$ such that
$d_N(\alpha)\neq 0$ is $h_\alpha$, and $\sum _{k>0}\varphi (k)d_k(\alpha) =
gh_\alpha/g_\alpha$, by Lemma~\ref{formulas}. Moreover
$d_1(\alpha)$ is given by Proposition~\ref{d1}.  Then Corollary~\ref{firstcase} is
the first case
$x=1, y=h_\alpha$ of the following:

\begin{proposition}\label{circular} With notation as above,  the sequence
$\{d_k(\alpha)\}$ has the circular symmetry property with respect to
$h_\alpha$ and $gh_\alpha/g_\alpha$.
\end{proposition}
\begin{proof} If $h_\alpha \leq 3$, it is easy to check that the above
  conditions
follow from Corollary~\ref{firstcase} and the fact that $\sum
_{k>0}\varphi (k)d_k(\alpha)
= gh_\alpha/g_\alpha$. The remaining cases are: $F_4$ with $\alpha$
the root such
that $h_\alpha = 4$,  $E_7$ with $\alpha$ the root such that $h_\alpha =
4$, or $E_8$ with $\alpha$ a root such that $h_\alpha =
4,5,6$. These cases may be checked by hand.
\end{proof}

It would be very interesting to find a more conceptual proof of
Proposition~\ref{circular}.

\subsection{The moduli space of semistable $G$-bundles}\label{s1.5}

 If
$\xi_0$ is a $C^\infty$ principal $G$-bundle over $E$, then there is a
characteristic class $c_1(\xi_0) \in H^2(E; \pi_1( G))\cong
\pi_1(G)=\langle c\rangle$.
Let $\mathcal{M}(G,c)$ be the set of S-equivalence classes of holomorphic
semistable $G$-bundles $\xi$ over $E$ with $c_1(\xi) = c$. Suppose that $K$ is
the compact form of $G$. Then $\langle c\rangle= \pi_1(K)$.
Let $\widetilde K$ be the universal covering group of $K$.
By the theorem of
Narasimhan-Seshadri
and Ramanathan and~\cite[5.8(i)]{FM}, $\mathcal{M}(G,c)$ is
homeomorphic to the the
space
$$\{(x,y)\in \widetilde K\times \widetilde K: xyx^{-1}y^{-1} =c\}/\widetilde K,$$
where the action of $\widetilde K$ is by simultaneous
conjugation. By~\cite[3.2]{BFM}, corresponding to $c\in Z(\widetilde G)$
there is an element $w_c\in W$ and an affine
isomorphism
$\varphi_c$ of $V^*$ which permutes the set $\widetilde \Delta$ and
thus acts on
the fundamental alcove $A$ corresponding to the Weyl chamber defined
by $\Delta$.

\begin{defn}\label{defrc} The Weyl element $w_c$ induces a   permutation of the
set
$\widetilde\Delta$ which induces an automorphism of the extended
Dynkin diagram of $G$.
Thus
$g_{w_c(\alpha)} = g_\alpha$ for all $\alpha \in \widetilde\Delta$.  Let
$\ov\alpha$ denote the
$w_c$-orbit of $\alpha$ and let $n_{\ov \alpha}$ be the cardinality of $\ov
\alpha$. Set
$g_{\ov \alpha} = n_{\ov \alpha}g_\alpha$, for any choice of $\alpha \in \ov
\alpha$. Let $n_0$ be the gcd of the integers $g_{\ov \alpha}$. Let 
$\widetilde\Lambda$ be the free abelian group  with basis
$\widetilde \Delta$. Then   $w_c$
acts on $\widetilde\Lambda$,
preserving the relation $\sum _{\beta \in \widetilde
\Delta}g_\beta\beta\spcheck =0$, and so  $w_c$ acts on the quotient which is
$\Lambda$. Define $r_c$
by: $r_c+1$ is the cardinality of
$\widetilde \Delta /w_c$.
\end{defn}

Let $A^c$ be the fixed subspace of $\varphi_c$, acting on $A$. With
$V^* =\Lambda \otimes
\Ar$,   let $T=V^*/\Lambda$, so that $T$ is a real torus of dimension $r$. Then
$W$ and $w_c$ act on $T$. Let $T_0 = (T^{w_c})^0$ be the identity
component of the group $T^{w_c}$. Thus $r_c =\dim T_0$.

For an abelian group $\mathcal A$ and an automorphism $\sigma$ of
$\mathcal A$, we
denote as usual by $\mathcal A^\sigma$ the subgroup of invariants of
$\sigma$ and
by $\mathcal A_\sigma$ the group $\mathcal A/\operatorname{Im}(\Id -\sigma)$ of
coinvariants of $\sigma$.  The automorphism   $w_c$ acts on $T=V^*/\Lambda$,
and we can  define $T^{w_c}$ and $T_{w_c}$ as before. Let $T_0 = (T^{w_c})^0$
be the identity component of the group $T$. There is an induced map
from $T_0$ to
$T_{w_c}$ and it is finite. We then have:

\begin{theorem}\label{degree}
Fix $x_0\in A^c$ and $y_0$ an element in the normalizer $N_T(\widetilde K)$ of
$T$ in $\widetilde K$ which projects to $w_c$ in $W=N_T(\widetilde K)/T$.
Then for $s,t\in T_0$ the pair $x=sx_0,y=ty_0$ satisfies $[x,y]=c$.
We define a map
$T_0\times T_0 \to \mathcal{M}(G,c)$
by $(s,t)\mapsto [(sx_0,ty_0)]$. This map
 is finite and surjective. Its degree is
$$(r_c)!\frac{\det (I_0|\Lambda^{w_c})}{n_0}\prod _{\ov \alpha}g_{\ov
  \alpha},$$
where the product is over the $w_c$-orbits of $\widetilde \Delta$.
\end{theorem}
\begin{proof} By \cite[Lemma 6.1.7]{BFM}, with $\widetilde K$ and $c$
  as above,   every pair
$(x,y)$ such that $xyx^{-1}y^{-1} =c$ is conjugate to such a pair with
$x\in A^c$
and $y\in T_0\cdot w_c$. This proves that the map is surjective.
Clearly, it is finite-to-one.
If $x$ is in the interior of $A^c$, then it
is regular,
and the only further possible conjugation is via an element $t\in T$,
which acts
on $y$ via $t-w_c(t)$. Thus, a fundamental domain for the map
$T_0\times T_0\to
\mathcal{M}(G,c)$ is given by $A^c\times S$, where $S$ is a
fundamental domain for
the quotient map $T_0 \to T_{w_c}$. It follows that the degree of the
map $T_0\times T_0\to
\mathcal{M}(G,c)$ is the product of the degree of the map from $T_0$
to $T_{w_c}$
with the ratio $\operatorname{vol}(T_0)/\operatorname{vol}(A^c)$,
where volume is
computed with respect to any Weyl invariant metric. We consider these
two integers separately.

\begin{lemma}\label{orders} Let $\ov \alpha \spcheck = \sum _{\alpha \in \ov
\alpha}\alpha\spcheck$.
Then
$$\Lambda^{w_c} \cong \bigoplus_{\ov \alpha}\Zee\cdot \ov \alpha
\spcheck\Bigg/\sum_{\ov \alpha}g_{ \alpha}{\ov \alpha}\spcheck.$$
Moreover, the set $\{\ov \alpha\spcheck: \ov \alpha \neq \ov \alpha_0\}$ is an
integral basis for
$\Lambda^{w_c}$. Finally, for each orbit $\ov \alpha$, choose $\alpha \in \ov
\alpha$ and let $e_{\ov \alpha}$ be the image of $\alpha$ in  $\Lambda_{w_c}$.
Then
$$\Lambda_{w_c} \cong \bigoplus_{\ov \alpha}\Zee\cdot e_{\ov
\alpha}\Bigg/\sum_{\ov \alpha}g_{\ov \alpha}e_{\ov \alpha}.$$
\end{lemma}
\noindent {\bf Proof.} There is
an exact sequence
$$0 \to \Zee \to \widetilde \Lambda \to \Lambda \to 0,$$
where $1\in \Zee \mapsto \sum _{\beta \in \widetilde
  \Delta}g_\beta\beta\spcheck$.
The homomorphisms in this sequence are equivariant with respect
to the action of ${w_c}$. Moreover, ${w_c}$ acts on $\widetilde \Lambda$ by a
permutation of the basis. The proof of the lemma follows easily by
considering the
associated long exact sequence
$$0\to \Zee^{w_c} \to (\widetilde \Lambda)^{w_c} \to \Lambda ^{w_c} \to \Zee_{w_c}
\to (\widetilde \Lambda)_{w_c} \to \Lambda _{w_c} \to 0. \endproof$$

\begin{corollary}\label{torsorder} The torsion subgroup
  $\operatorname{Tor}\Lambda _{w_c} \cong \Zee/n_0\Zee$, and
$$\Lambda_{w_c}/\operatorname{Tor}\Lambda _{w_c}  \cong \bigoplus_{\ov
\alpha}\Zee\cdot e_{\ov
\alpha}\Bigg/\sum_{\ov \alpha}\frac{g_{\ov \alpha}}{n_0}e_{\ov \alpha}.\endproof$$
\end{corollary}

\begin{lemma} The order of $T^{w_c}/T_0$ is $n_0$. The natural map $T_0 \to
T_{w_c}$ is finite and surjective of degree $\prod _{\ov \alpha}n_{\ov
\alpha}\Big/n_0$.
\end{lemma}
\begin{proof} Beginning with the short exact sequence of ${w_c}$-modules
$$0 \to \Lambda \to V \to T \to 0,$$
we get a long exact sequence
$$0 \to \Lambda^{w_c} \to V^{w_c} \to T^{w_c} \to  \Lambda_{w_c} \to V_{w_c}
\to T_{w_c} \to 0.$$
The quotient $V^{w_c} /\Lambda^{w_c}= T_0$. Since $T^{w_c}/T_0$ is a finite
group and $V_{w_c}$ is torsion free, the induced map $T^{w_c}/T_0 \to
\Lambda_{w_c}$ is an isomorphism from $T^{w_c}/T_0$ to
$\operatorname{Tor}\Lambda _{w_c} $, and hence $T^{w_c}/T_0 \cong
\Zee/n_0\Zee$. Moreover, it is clear that the degree of the map from $T_0$ to
$T_{w_c}$ is the index of the image of $\Lambda^{w_c}$ in
$\Lambda_{w_c}/\operatorname{Tor}\Lambda _{w_c}$. By Lemma~\ref{orders} and
Corollary~\ref{torsorder}, it suffices to compute the order of the quotient of
$\bigoplus_{\ov \alpha}\Zee\cdot e_{\ov \alpha}$ by the relations $n_{\ov
\alpha}e_{\ov \alpha}$ and $\sum _{\ov \alpha}g_{\ov \alpha}e_{\ov \alpha}$. Since
$g_{\ov \alpha_0}=1$, it is clear that the quotient has order
$\prod _{\ov \alpha}n_{\ov \alpha}\Big/n_0$.
\end{proof}

Now we compute $\operatorname{vol}(T_0)/\operatorname{vol}(A^c)$ using
the volume
determined by the inner product $I_0$. The alcove $A$ has vertices
equal to $0$ and
$\varpi_\alpha\spcheck/h_\alpha$, $\alpha \in \Delta$. By
Lemma~\ref{next}, the vertices are $0$ and $\alpha^*/g_\alpha$, where
the $\alpha^*$ are the dual basis
with respect to $I_0$  to the basis $\alpha\spcheck$ of $\Lambda$. Now
we have the following elementary lemma, whose proof is left to the reader:

\begin{lemma}  Let $A$ be a simplex in $\Ar ^n$ with vertices $0=e_0,
e_1, \dots, e_n$. Let $\varphi$ be an affine linear transformation of
$\Ar ^n$ which acts via a permutation of the vertices of $A$.
Suppose that the orbits of $w$ on the vertices are $\mathbf{o}_0,
\dots,\mathbf{o}_s$, with $0\in \mathbf{o}_0$. If $n_{\mathbf{o}}$ is the order
of the orbit $\mathbf{o}$, set
$$v_{\mathbf{o}} =\frac{1}{n_{\mathbf{o}}}\sum _{e_i\in \mathbf{o}}e_i.$$
Then the fixed set of
$A^\varphi$ for the action of $\varphi$ on $A$ is a simplex with
vertices
$$v_{\mathbf{o}_0}, v_{\mathbf{o}_1} +v_{\mathbf{o}_0}, \dots, v_{\mathbf{o}_s}
+v_{\mathbf{o}_0}. \endproof$$
\end{lemma}

Applying the lemma to $A^c$, we see that $A^c$ is a translate of the simplex in
$(V^*)^{w_c}$ spanned by $0$ and $1/g_{\alpha}v_{\ov \alpha}$, where $v_{\ov
\alpha} = 1/n_{\ov \alpha}\sum _{\alpha \in \ov \alpha}\alpha^*$ and
$\alpha$ is
any representative for $\ov \alpha$.  It follows by Lemma~\ref{orders} that
$\{\ov\alpha\spcheck; \ov \alpha \neq \ov \alpha_0\}$ is an integral basis for
$\Lambda^{w_c}$, where $\ov\alpha\spcheck =
\sum _{\alpha \in \ov \alpha}\alpha\spcheck$. Since
$$I_0(v_{\ov \alpha}, \ov \beta\spcheck) =
\begin{cases}
0, &\text{if $\ov \alpha \neq \ov \beta$;}\\
1, &\text{if $\ov \alpha = \ov \beta$,}
\end{cases}$$
we see that $\{\ov\alpha\spcheck\}$ and $\{v_{\ov \alpha}\}$ are dual bases for
the restriction of $I_0$ to $(V^*)^{w_c}$. Now
$$\operatorname{vol}(A^c) = \frac{\operatorname{vol}(C_1)}{(r_c)!\prod _{\ov
\alpha}g_\alpha},$$
where $C_1$ is the parallelepiped spanned by the basis $\{v_{\ov
  \alpha}\}$ and as
usual $\alpha$ is any representative for $\ov \alpha$. On the other hand,
$\operatorname{vol}(T_0) = \operatorname{vol}(C_2)$, where $C_2$ is the
parallelepiped spanned by the dual basis $\{\ov\alpha\spcheck\}$. Thus
$$\frac{\operatorname{vol}(T_0)}{\operatorname{vol}(A^c)} =
(r_c)!\left(\prod _{\ov
\alpha}g_\alpha\right)\frac{\operatorname{vol}(C_2)}{\operatorname{vol}(C_1)} =
(r_c)!\left(\prod _{\ov
\alpha}g_\alpha\right)\det (I_0|\Lambda^{w_c}).$$
To complete the proof of Theorem~\ref{degree}, the degree in question is the
product
$$(r_c)!\left(\prod _{\ov
\alpha}g_\alpha\right)\det (I_0|\Lambda^{w_c})\cdot \left(\frac{\prod _{\ov
\alpha}n_{\ov \alpha}}{n_0}\right) = (r_c)!\frac{\det
(I_0|\Lambda^{w_c})}{n_0}\prod _{\ov \alpha}g_{\ov \alpha},$$
as claimed.
\end{proof}

\section{Bundles over maximal parabolic subgroups}

\subsection{Description of bundles and their automorphisms}

Fix $\alpha\in \Delta$. We consider $L^{\alpha}$-bundles
$\eta$ over $E$
 such that $c_1(\eta\times _{L^\alpha} G) = c$, and will refer to such
 a bundle as a {\sl
unliftable bundle of type $c$}.
The primitive dominant character
$\chi_0$ of $P^\alpha$ lifts to a character on
$\widetilde P^\alpha$ which is a positive power of the primitive dominant
character $\varpi_\alpha$ of $\widetilde P^\alpha$. We denote this power
by $o_{c,\alpha}$.
 Note that $\varpi_\alpha(c)$ is well-defined  as an element of $\Q/\Zee$ and
$o_{c,
\alpha}$ is its order. In fact, we have:

\begin{lemma} Let $\beta$ be a root such that $c\equiv
  \varpi_\beta\spcheck \pmod \Lambda$. Then
$o_{c, \alpha}$, the order of
$\varpi_\alpha(c)$, is the order of
$\varpi_\alpha(\varpi_\beta\spcheck) \mod\Zee$.
\endproof
\end{lemma}

In the notation of~\cite[\S 3.4]{BFM}, $o_{c, \alpha} =1$ if and only
if $\alpha
\notin \Delta(c)$. The Dynkin diagram for $\Delta(c)$ is a union of diagrams of
$A$-type  and is described in the tables at the end of~\cite{BFM}.

If $\eta$ is an $L^\alpha$-bundle, then the character
$\chi_0$  defines an associated line bundle $\eta \times_{
  L^{\alpha}}\Cee$. This line
bundle is the {\sl determinant\/} of $\eta$, which we write as $\det \eta$. Its
degree is called the {\sl degree\/}
of $\eta$ and is denoted $\deg \eta$.

 Let
$\eta\to E$ be a principal
$ L^{\alpha}$-bundle whose degree $d$ is negative. We shall study
 the corresponding bundles
$\eta\times_{ L^{\alpha}} P^{\alpha}$ and
$\xi=\eta\times_{ L^{\alpha}} G$.
Associated to $\eta$ and the action of
$ L^{\alpha}$ on the Lie algebra
$\frak{g}$ there is the vector bundle
$\eta\times_{ L^{\alpha}}{\frak g}=\ad \xi$. The Lie algebra ${\frak
g}$ decomposes under $ L^{\alpha}$ as ${\frak u} \oplus
\frak l\oplus {\frak u}_-$ where $\frak l$ is the Lie
algebra of $ L^{\alpha}$,
${\frak u}$ is the subspace of ${\frak g}$ on which
$\varphi_\alpha(\Cee^*)\subseteq  L^{\alpha}$ acts with  positive weights,
and ${\frak u}_-$ is the subspace of
${\frak g}$ on which
$ \varphi_\alpha(\Cee^*)\subseteq  L^{\alpha}$ acts with
negative weights. Since the coefficients of $\zeta_\alpha$ are
nonnegative,
${\frak u}$ is the Lie algebra of $U^{\alpha}$, $\frak p  = \frak l
\oplus {\frak u}$ is the Lie algebra of $ P^{\alpha}$, and
$\frak u_-$ is the orthogonal space to $\frak u$ under the Killing form.
Clearly:

\begin{lemma}\label{adG} There is a direct sum decomposition
$$\ad \xi = \eta\times
_{ L^{\alpha}}{\frak g} = \ad_{ L^{\alpha}}\eta\oplus {\frak
u}(\eta)\oplus{\frak u}_-(\eta).$$
The action of the
$\Cee^*\subseteq  L^{\alpha}$ is trivial on
$\ad_{ L^{\alpha}}\eta$ and has positive {\rm{(}}resp. negative{\rm{)}} weights on
${\frak u}(\eta)$ {\rm{(}}resp. ${\frak u}_-(\eta)${\rm{)}}. \endproof
\end{lemma}

We now have:

\begin{lemma}\label{list} For  every negative integer $d$, there exists a 
semistable  $L^{\alpha}$-bundle
$\eta$ over $E$ of degree $d$. There is an unliftable semistable
$L^{\alpha}$-bundle $\eta$ of type $c$ if and only if  $\deg \eta/o_{c, \alpha}
\equiv \varpi_\alpha(c) \mod \Zee$.  In particular, if $\eta$ is unliftable of
type $c$ and degree $-1$,
then $-1/o_{c, \alpha} \equiv \varpi_\alpha(c) \mod \Zee$.
For a   semistable $L^{\alpha}$-bundle
$\eta$ of degree $d<0$, we have:
\begin{enumerate}
\item[\rm (i)] The bundle $  \eta  \times _{ L^{\alpha}} G$ is unstable;
\item[\rm (ii)] The parabolic $ P^{\alpha}_-$
opposite to $ P^{\alpha}$  is a Harder-Narasimhan parabolic of
$\xi$;
\item[\rm (iii)] $\frak u(\eta)$ is a direct sum of semistable vector
  bundles of
strictly negative degrees.
\item[\rm (iv)] The Atiyah-Bott point of
$\eta$ as defined in \cite{AtBo, FM, FMAB} is given by
$$\mu(\eta) = \frac{d\zeta_\alpha}{o_{c, \alpha}m_\alpha}=
\frac{dn_\alpha}{o_{c, \alpha}m_\alpha}\varpi_\alpha\spcheck.$$
\end{enumerate}
\end{lemma}
\begin{proof} The dominant character $\chi_0$ lifts to the character
$\varpi_\alpha^{o_{c,\alpha}}$ on $\widetilde L^\alpha$.
By  \cite[Definition 2.1.1]{FMAB}, $\mu(\eta)$ is the
  unique point
$\mu$ in the center of   $\frak l$ such that $o_{c,
  \alpha}\varpi_\alpha(\mu) = d$. Thus $\mu =  d\zeta_\alpha/o_{c,
  \alpha}m_\alpha$, showing (iv). The
congruence condition $\deg \eta/o_{c, \alpha} \equiv \varpi_\alpha(c) \mod
\Zee$ follows
from \cite[Lemma 2.1.2 (ii)]{FMAB}. Statement (iii) follows from \cite[Lemma
2.2.1]{FMAB}. The remaining statements are clear.
\end{proof}

Our goal now is to study the spaces $H^1(E; \ad \xi)$ and $H^0(E; \ad \xi)$. We
shall primarily be interested in the case where $\eta$ is
semistable. It is enough
to study the spaces $H^i(E; \ad_{ L^{\alpha}} \eta)$, $H^i(E; \frak
u(\eta))$, and $H^i(E; \frak u_-(\eta))$. Since
$L^{\alpha}$ has a one-dimensional center, regardless of the
choice of $\eta$ we must always have $\dim H^0(\ad_{ L^{\alpha}} \eta)
\geq 1$, and thus, since $\deg \ad_{ L^{\alpha}} \eta = 0$, $\dim
H^1(\ad_{ L^{\alpha}}
\eta) \geq 1$ as well by Riemann-Roch on $E$. More precisely, we have:

\begin{lemma}\label{adL}  Let $\widehat S^{\alpha}$ be the quotient of
  $L^{\alpha}$ by  its center. Let $\eta$ be a semistable
  $L^{\alpha}$-bundle and  let
$\widehat \eta$ be the induced
$\widehat S^{\alpha}$-bundle. Let $r(\widehat \eta)$ be the dimension of $\Aut
_{\widehat S^{\alpha}}(\widehat \eta)$. Then $\dim
  H^0(E;\ad_{L^{\alpha}} \eta) =
\dim H^1(E;\ad_{L^{\alpha}} \eta) =1+r(\widehat \eta)$.
\end{lemma}
\begin{proof} On the level of Lie algebras, there is a direct sum
decomposition $\frak l = \Cee \oplus \Lie (\widehat S^{\alpha})$, where $\Cee =
\Lie (Z(L^{\alpha}))$, and the proof follows.
\end{proof}
 We turn next to the groups $H^i(E; \frak u (\eta))$ and $H^i(E; \frak u_-
(\eta))$:

\begin{lemma}\label{uplus} Let $\eta$ be a principal $L^{\alpha}$-bundle
of negative degree. Then
$$\dim  H^1(E;{\frak u}(\eta)) \geq
-\deg\left({\frak u}(\eta)\right),$$ with equality
holding if and only if $H^0(E;{\frak u}(\eta))=0$.
Likewise,
$$\dim H^0(E;{\frak
u}_-(\eta)) \geq\deg\left({\frak u}_-(\eta)\right),$$
with equality holding if and only if  $H^1(E;{\frak
u}_-(\eta))=0$.  Finally, if $\eta$ is semistable, then $H^0(E;{\frak
u}(\eta))= H^1(E;{\frak
u}_-(\eta))=0$.
\end{lemma}
\begin{proof} The first two statements are immediate from Riemann-Roch on
$E$. The final one follows from   Statement (iii) of Lemma~\ref{list}.
\end{proof}

{}From the decomposition
$\frak p  =\frak l \oplus {\frak u}$, it follows
that $H^0(E;{\frak u}(\eta))=0$ if and only if
$$\dim \Aut_{L^{\alpha}}\eta = \dim \Aut _{P^{\alpha}}(\eta
\times _{L^{\alpha}}P^{\alpha}).$$
The vanishing of $H^1(E;{\frak u}_-(\eta))$ says on the other hand that the map
$$H^1(E; \ad _{ P^{\alpha}}(\eta\times _{ L^{\alpha}} P^{\alpha}))\to
 H^1(E; \ad _{G}(\eta\times _{ L^{\alpha}} G))$$
is an isomorphism. In particular, every small deformation of the $ G$-bundle
$\eta\times _{ L^{\alpha}} G$ arises  from a small $
P^{\alpha}$-deformation of $\eta\times _{L^{\alpha}}P^{\alpha}$.

To complete the determination of $H^1(E;{\frak u}(\eta))$, we must compute the
degree of ${\frak u}(\eta)$.

\begin{proposition}\label{degrees} We have $\deg {\frak u}(\eta)
= (\deg \eta) d_1(\alpha)/o_{c, \alpha}$, where $ d_1(\alpha)$ is the
integer defined
in Lemma~\ref{defd1} for $\widetilde G$. Likewise,  $\deg {\frak u}_-(\eta)
= -(\deg \eta) d_1(\alpha)/o_{c, \alpha}$.
\end{proposition}
\begin{proof} We compute the degree  of the line
bundle  $\bigwedge^{\text{\rm top}}{\frak u}(\eta)$.
Let $\widetilde\chi_0$ be the dominant character for $\widetilde P^\alpha$.
The line bundle $\bigwedge^{\text{\rm top}}{\frak u}(\eta)$ is
associated to $\eta$ by the character
$\chi_+\colon L^{\alpha}\to
\Cee^*$. By Lemma~\ref{defd1}, $\chi_+$ lifts to the character
$\widetilde\chi_+ =
\widetilde\chi_0^{ d_1(\alpha)}$ of $\widetilde L^\alpha$. Since $\chi_0$ lifts
to  $\widetilde
\chi_0^{o_{c,\alpha}}$ on $\widetilde L^\alpha$ and
since the line bundle associated to
$\eta$ by the character $\chi_0$ has degree $\deg \eta$, the
degree of ${\frak u}(\eta)$ is $(\deg \eta)\cdot  d_1(\alpha)/o_{c, \alpha}$.
A similar argument (or duality) handles the case of
${\frak u}_-(\eta)$.
\end{proof}

Combining Lemmas~\ref{adL} and \ref{uplus} with
Proposition~\ref{degrees}, we have:

\begin{corollary}\label{dimad} If $\eta$ is semistable of negative
  degree, then
$$\dim
H^1(E;{\frak u}(\eta)) =\dim
H^0(E;{\frak u}_-(\eta)) = -(\deg \eta) d_1(\alpha)/o_{c, \alpha}.$$
Thus
$$\dim
H^0(E;\ad (\eta\times _{L^{\alpha}} G)) =\dim
H^1(E;\ad (\eta\times _{L^{\alpha}} G)) =1 +r(\widehat \eta) -(\deg
\eta) d_1(\alpha)/o_{c, \alpha}. \endproof$$
\end{corollary}

\subsection{The $\Cee^*$-action in cohomology}

The Lie
algebra $\frak u $ is a direct sum of the subspaces $\frak u ^k$, $k>0$, where
$\frak u ^k$ is the sum of all the root spaces $\frak g^\beta$ where the
coefficient of
$\beta$ is exactly
$k$. By Proposition~\ref{highest}, $\frak u ^k$ is an irreducible
$L^{\alpha}$-module. Thus ${\frak u}(\eta)$ is the direct sum of
vector bundles
${\frak u}^k (\eta)$ associated to irreducible
representations of $L^{\alpha}$. By~\cite{Ra}, ${\frak u}^k (\eta)$
is  semistable. Our goal now is to study the action of
$\Cee^*=\varphi_\alpha(\Cee^*)$ on ${\frak
u}(\eta)$ and on $H^1(E;{\frak u}(\eta))$.
As we saw in Section~\ref{para2},
$\widetilde\varphi_\alpha(\Cee^*)\subseteq \widetilde L^\alpha$ acts on ${\frak
u}^k$ with weight $kn_\alpha$. Thus, by Lemma~\ref{ncalpha}, 
$\varphi_\alpha(\Cee^*)$ acts on  ${\frak u}^k $ with weight $kn_{c,\alpha}$. 

\begin{lemma}\label{weights} The degree $\deg {\frak u}^k (\eta)$ is equal to
$(\deg \eta)\cdot i(\alpha, k)/o_{c, \alpha}$, where
$$i(\alpha, k) = kc(\alpha,
k)/\varpi_\alpha(\varpi_\alpha\spcheck)= kn_\alpha c(\alpha, k)/m_\alpha.$$
 The slope of ${\frak u}^k
(\eta)$ is $k\varpi_\alpha(\mu(\eta))/\varpi_\alpha(\varpi_\alpha\spcheck)$, where
$\mu(\eta)$ is the Atiyah-Bott point of
$\eta$.  Thus, if
$\eta$ is semistable of negative degree, then
$\dim H^1(E;{\frak u}^k (\eta)) = - (\deg \eta)\cdot
i(\alpha, k)/o_{c, \alpha}$.
\end{lemma}
\begin{proof} By Lemma~\ref{deguk} applied to $\widetilde G$,   the
  character of
$L^{\alpha}$ defined by the determinant on ${\frak u}^k $ lifts to the
character
$\chi_0^{i(\alpha, k)}$ on $\widetilde L^\alpha$.  The degree of ${\frak
  u}^k (\eta)$ is thus
$(\deg \eta)\cdot i(\alpha, k)/o_{c, \alpha}$. It follows that
the slope of ${\frak u}^k
(\eta)$ is  $kn_\alpha\deg \eta/o_{c, \alpha}m_\alpha=
k\varpi_\alpha(\mu)/\varpi_\alpha(\varpi_\alpha\spcheck)$.
\end{proof}

\section{Special roots and the associated bundles}

In \S\ref{3.1}--\ref{3.3}, we assume that $G$ is
simply connected. We will defer the
discussion of the non-simply connected case until \S\ref{ns}.

\subsection{Definition of special roots}\label{3.1}

\begin{defn} A simple root $\alpha$ is {\sl special\/} if
\begin{enumerate}
\item[(i)] The Dynkin diagram associated to $\Delta -\{\alpha\}$ is a union of
diagrams of $A$-type;
\item[(ii)] The simple root $\alpha$ meets each component of the Dynkin diagram
associated to $\Delta -\{\alpha\}$ at an end of the component;
\item[(iii)] The root $\alpha$ is a long root.
\end{enumerate}
If $R$ is of type $A_n$, then every simple root  is special. All other
irreducible
root systems have a unique special simple root. It corresponds to the unique
trivalent vertex if the Dynkin diagram is of type
$D_n, n\ge 4$ or $E_n, n=6,7,8$. For  $R=C_n, n\ge 2$ and $G_2$, it is the
long simple root.
For $R=B_n, n\ge 2$,  and $F_4$ it is the  unique
long simple root which is not orthogonal to a short simple root.
\end{defn}

We shall investigate the structure of the group
$L^\alpha$ and the space
$H^1(E; \frak u(\eta))$ more closely in case $\alpha$ is special.

Let $\alpha$ be special and let $\eta$ be a semistable bundle over
$P^\alpha$ of
degree $-1$. By the results of
\cite{FMAB}, the unstable bundles
$\eta\times_{L^\alpha}G$ are minimally unstable bundles, in the sense
that every
small deformation of such a bundle is either of the same type or semistable.
Moreover, if $G$ is not of $A$-type, then for every unstable
$G$-bundle $\xi$, there is a small deformation of $\xi$ to a bundle of
the form
$\eta\times_{L^\alpha}G$. Now for every unstable $G$-bundle $\xi$, there is the
Harder-Narasimhan reduction to a  parabolic subgroup, not necessarily maximal,
and in fact $\xi$ reduces to a bundle $\eta$ the Levi factor $L$. Let
$P$ be the
opposite parabolic to the Harder-Narasimhan parabolic and let $\frak
u$ be the Lie
algebra of the unipotent radical of $P$. It is easy to see that the function
$\xi\mapsto \dim H^1(E; \frak u(\eta))$ is strictly decreasing for the
Atiyah-Bott
ordering, and hence attains its minimum in the case where $P$ is a maximal
parabolic corresponding to a special root and $-\deg \eta$ is
minimal. We shall see this directly below.

\subsection{Bundles associated to special roots}\label{3.2}

Our first lemma determines the structure of $S^{\alpha}$ and
$L^{\alpha}$ in case $\alpha$ is special:

\begin{lemma}\label{Lalpha2} Suppose that $\alpha$ is special, and let
  $t$ be the
number of components in the Dynkin diagram of $S^\alpha$. Then there exist
integers $n_i\geq 2$ such that
$S^{\alpha}\cong
\prod _{i=1}^tSL_{n_i}(\Cee)$ and  $m_\alpha = \operatorname{lcm}
(n_i)$. Moreover,
$$L^{\alpha}\cong \left\{(A_1, \dots, A_t) \in \prod
  _{i=1}^tGL_{n_i}(\Cee): \det A_1= \cdots = \det A_t\right\},$$
in such a way that the primitive dominant character of $L^{\alpha}$
corresponds to the common value of the determinant.
\end{lemma}
\begin{proof} It follows from the definition of a special root that
  $S^{\alpha}\cong \prod
_{i=1}^tSL_{n_i}(\Cee)$. By
Lemma~\ref{Lalpha}, there is an isomorphism
$$L^{\alpha}\cong  S^{\alpha}\times_{\Zee/m_\alpha\Zee}\Cee^*,$$
where the image of $1\in \Zee/m_\alpha\Zee$ is mapped to  $e^{2\pi
  i/m_\alpha}\in
\Cee^*$ and to $e^{-2\pi i/n_i}\Id \in SL_{n_i}(\Cee)$. From this, we must have
$m_\alpha = \operatorname{lcm} (n_i)$. The map from
$S^{\alpha}\times\Cee^*$ to $\prod _{i=1}^tGL_{n_i}(\Cee)$ which is the natural
inclusion   $\prod _{i=1}^tSL_{n_i}(\Cee)\subseteq \prod
_{i=1}^tGL_{n_i}(\Cee)$
and which maps $\lambda \in \Cee^*$ to $(\lambda ^{m_\alpha/n_1}\Id,
\dots, \lambda
^{m_\alpha/n_t}\Id)$ then factors to give an induced homomorphism
$S^{\alpha}\times_{\Zee/m_\alpha\Zee}\Cee^* \to \prod
_{i=1}^tGL_{n_i}(\Cee)$. It
is clear from the construction that this induced homomorphism is injective and
that its image is the subgroup of matrices of equal
determinant. Let
$\det\: L^{\alpha} \to \Cee^*$ denote the value of any of these
determinants under
the inverse isomorphism. For
$\lambda \in
\Cee^*$, we see that
$\det \circ \varphi_\alpha = \lambda ^{m_\alpha}$, and hence $\det = \chi_0$.
\end{proof}

\begin{lemma} If $\alpha$ is special, the positive integers $n_\alpha$
  are equal to
$1$ except in the following cases:
\begin{enumerate}
\item[\rm{(a)}] If $G = SL_n(\Cee)$ and $\alpha$ corresponds to the
$k^{\text{\rm th}}$ vertex in the usual ordering of the simple roots, then
$n_\alpha = n/\gcd(k,n)$.
\item[\rm{(b)}] If $G$ is of type $B_n$ and $n$ is even, then $n_\alpha=2$.
\item[\rm{(c)}] If $G$ is of type $C_n$, then $n_\alpha=2$.
\item[\rm{(d)}] If $G$ is of type $D_n$ and $n$ is odd, then $n_\alpha=2$.
\end{enumerate}
\end{lemma}
\begin{proof} If the center of $G$ is trivial, then
  $\varpi_\alpha\spcheck$ is a
primitive element of $\Lambda$, and hence $n_\alpha = 1$. This handles
the cases
$E_8, F_4, G_2$. Next suppose that $R$ is simply laced and not of type $A_n$,
so that the Dynkin diagram of $R$ is a $T_{p,q,r}$ diagram, with $(p,q,r) =
(2,2,n)$ or $(2,3,s)$ with
$s=3,4,5$. Let $N = (1/p+ 1/q+1/r-1)^{-1}$, so that $N = n$, $6,12,
30$ in the
respective cases above. In particular $N\in \Zee$. There exists a
labeling of the
roots as
$\{\alpha,
\beta_1, \dots, \beta_{p-1}, \gamma_1, \dots, \gamma_{q-1}, \delta_1, \dots,
\delta_{r-1}\}$, where $\beta_1, \gamma_1, \delta_1$ are ends of the diagram,
$\langle \beta_i, \beta_{i+1}\rangle =-1$, $1\leq i\leq p-2$ and
similarly for the $\gamma_j$ and $\delta_k$, and
$\beta_{p-1}, \gamma_{q-1}, \delta_{r-1}$ meet $\alpha$, such that
$$\varpi_\alpha\spcheck = \frac{N}{p}\sum _{i=1}^{p-1}i\beta_i\spcheck +
\frac{N}{q}\sum _{j=1}^{q-1}j\gamma_j\spcheck + \frac{N}{r}\sum
_{k=1}^{r-1}k\delta_k\spcheck + N\alpha.$$
It follows that $\varpi_\alpha\spcheck$ is integral, and hence $n_\alpha = 1$,
unless
$(p,q,r) = (2,2,n)$ and
$n$ is odd, in which case $n_\alpha = 2$. A similar argument handles
the case of
$A_n$. In case $B_n$, if we number the roots as in \cite{Bour} beginning at the
long end of the Dynkin diagram, then $\alpha = \alpha_{n-1}$ and
$$\varpi_{\alpha_{n-1}}\spcheck = \sum _{i=1}^{n-1}i\alpha_i\spcheck +
\frac{n-1}2\alpha_n\spcheck.$$
Thus $n_\alpha = 1$ if $n$ is odd and $2$ if $n$ is even. Finally, for
the case of
$C_n$, again numbering the roots in order as in \cite{Bour} beginning
at a short
root, so that $\alpha = \alpha_n$, we have
$$\varpi_{\alpha_n}\spcheck = \frac12\sum _{i=1}^ni\alpha_i\spcheck .$$
Thus $n_\alpha = 2$.
\end{proof}

We turn now to the existence of special  bundles over $L^{\alpha}$.

\begin{proposition}\label{eta0}  Suppose that  $\alpha$
is special. Then there is a unique
principal $L^{\alpha}$-bundle $\eta_0$ over $E$ with the following
properties:
\begin{enumerate}
\item[\rm{(i)}] $\det \eta_0 =\scrO_E(-p_0)$
\item[\rm{(ii)}] For $1\leq i\leq t$, if $V_i$ is  the
vector bundle associated to the principal $GL_{n_i}(\Cee)$-bundle
obtained from the composition of the inclusion  $L^{\alpha}\subseteq
\prod_{i=1}^tGL_{n_i}(\Cee)$ followed by projection onto the $i^{\rm th}$
factor, then each
$V_i$ is a stable vector bundle.
\end{enumerate}
 The automorphism group of  $\eta_0$, as an
$L^{\alpha}$-bundle, is identified with  the center of
$L^{\alpha}$ which is of the form
$Z(S^{\alpha})\times_{\Zee/m_\alpha\Zee}\varphi_\alpha(\Cee^*)$,
acting by multiplication.
\end{proposition}
\begin{proof} Recall that, for every $d\geq 1$, there is a unique
stable vector bundle
$W_d$ of rank
$d$ over $E$ such that $\det W_d =\scrO_E(p_0)$. Given the
structure of $L^{\alpha}$ as in Lemma~\ref{Lalpha2}, it is clear that there
is a unique principal
$L^{\alpha}$-bundle, up to isomorphism, satisfying (i) and (ii)
above, with $V_i = W_{n_i}^*$ for every
$i$.     Since the vector bundles $V_i$ in (i) are simple,
the automorphism group of each of these is isomorphic to
the center of
$GL_{n_i}(\Cee)$ acting by multiplication. It then follows
that the
$L^{\alpha}$-automorphisms of $\eta_0$ are given by the action of the
center of $L^{\alpha}$ acting by multiplication.
\end{proof}

\begin{defn}\label{translate} If $\eta$ is an $L^{\alpha}$-bundle
  which satisfies
(ii) of the Proposition, together with the condition that $\deg \eta = -1$,
then $\eta$ is the pullback of $\eta_0$ under a translation map $E\to E$. We
call $\eta$ a {\sl translate\/} of
$\eta_0$.
\end{defn}

Let us describe the unstable bundle $\eta_0 \times _{L^{\alpha}}G$ for
the classical groups. First, we need the following
notation. As above, let $W_d$ be the unique stable vector
bundle over $E$ of rank $d$ and such that $\det W_d
=\scrO_E(p_0)$. Let $\theta_i$,
$i=1,2,3$ be the three nontrivial line bundles of order two
on $E$. Let
$$Q_3 = \theta _1\oplus\theta _2\oplus \theta _3$$ be the
corresponding rank three vector bundle. Fix isomorphisms
$\theta _i\otimes
\theta_i \to \scrO_E$ and give $Q_3$ the corresponding
diagonal symmetric bilinear form. Define similarly
$$Q_4 = \scrO_E\oplus \theta _1\oplus\theta _2\oplus \theta _3,$$
together with a similar choice of a diagonal symmetric
bilinear form. We then have:

\begin{proposition}\label{classical}  With notation as above, and supposing
that $\alpha$ is special, let $\eta_0$
be the principal $L^{\alpha}$-bundle constructed in Proposition~\ref{eta0}.
Then the
vector bundle associated to $\eta_0
\times _{L^{\alpha}}G$ under the standard representation of $G$ is:
\begin{enumerate}
\item[\rm{(i)}] $W_k^*\oplus W_{n-k}$, if $G= SL_n(\Cee)$ and
$\alpha$ is the root corresponding to the
$k^{\text{th}}$ vertex in the Dynkin diagram, ordered in the
usual way.
\item[\rm{(ii)}] $W_n^*\oplus W_n$, if $G= Sp(2n)$, where each
factor is isotropic and we choose an isomorphism $W_n^*\to
W_n^*$, unique up to a scalar, to define the alternating
form on the direct sum.
\item[\rm{(iii)}] $W_{n-2}^*\oplus Q_4 \oplus W_{n-2}$, if $G=
Spin(2n)$, where $Q_4$ is given the form described above,
$W_{n-2}^*$ and $W_{n-2}$ are isotropic, we choose an
isomorphism
$W_{n-2}^*\to W_{n-2}^*$, unique up to a scalar, to define
the symmetric form on the direct sum $W_{n-2}^*\oplus
W_{n-2}$, and $Q_4$ is orthogonal to this direct sum.
\item[\rm{(iv)}] $W_{n-1}^*\oplus Q_3 \oplus W_{n-1}$, if $G=
Spin(2n+1)$, where $Q_3$ is given the form described above,
$W_{n-1}^*$ and $W_{n-1}$ are isotropic, we choose an
isomorphism
$W_{n-1}^*\to W_{n-1}^*$, unique up to a scalar, to define
the symmetric form on the direct sum $W_{n-1}^*\oplus
W_{n-1}$, and $Q_3$ is orthogonal to this direct sum.
\end{enumerate}
\end{proposition}
\begin{proof} The cases $G=SL_n(\Cee)$ and $G= Sp(2n)$ follow
easily from the explicit descriptions of the maximal
parabolic subgroups and are left to the reader. In case $G
= Spin(2n)$, the corresponding maximal parabolic of
$SO(2n)$ is the set of
$g\in SO(2n)$ preserving an isotropic subspace of
$\Cee^{2n}$ of dimension $n-2$. The corresponding Levi
factor $L$ is the subgroup of matrices in
$GL_{n-2}(\Cee)\times GL_2(\Cee) \times GL_2(\Cee)$ with
equal determinant. If
$\rho_1$ is the representation of $L$ induced by the
standard representation of
$GL_{n-2}(\Cee)$ on $\Cee^{n-2}$ and
$\rho_2, \rho_3$ are the two representations of $L$ induced
by the standard representations of the second and third
factors of $L$ on $\Cee^2$, then it is easy to check that
the representation of $L$ on $\Cee^{2n}$ which is the
restriction of the standard representation of $Spin (2n)$
is just
$$\rho _1\oplus \rho_1^* \oplus (\rho_2\otimes \rho_3^*).$$
The vector bundle associated to
$\eta_0$ is thus
$$W_{n-2}^*\oplus W_{n-2} \oplus (W_2\otimes W_2^*).$$
Moreover this is an orthogonal direct sum with respect to
the induced form and
$W_{n-2}^*$ and $W_{n-2}$ are isotropic subspaces.
Furthermore, by a result of Atiyah \cite{Atiyah},
$W_2\otimes W_2^* \cong Q_4$, and since the each line
bundle summand of
$Q_4$ is not isomorphic to the dual of any other summand,
the direct sum decomposition of $Q_4$ must be orthogonal
with respect to the quadratic form and thus as described
above.

The case of $Spin(2n+1)$ is similar.
\end{proof}

\subsection{Cohomology dimensions and weight spaces}\label{3.3}

We have seen  that the bundles $\eta_0$ are minimally
unstable in the sense of deformation theory. Here we begin by showing
that their
deformation spaces have minimal dimension among all unstable bundles.

\begin{theorem}\label{mindim} Let $\eta_0$ be the bundle described in
Proposition~\ref{eta0}. Then
$$\dim H^0(E; \ad (\eta_0\times_{L^{\alpha}}G)) = \dim H^1(E; \ad
(\eta_0\times_{L^{\alpha}}G)) = r+2.$$
If $\xi$ is any unstable $G$-bundle, then $\dim H^1(E; \ad \xi) \geq
r+2$, with
equality if and only if $\xi$ is isomorphic to $\eta\times
_{L^{\alpha}}G$, where $\eta$ is a translate of $\eta_0$.
\end{theorem}
\begin{proof} First, by Corollary~\ref{dimad},
$$\dim H^0(E; \ad (\eta_0\times_{L^{\alpha}}G)) = \dim H^1(E; \ad
(\eta_0\times_{L^{\alpha}}G)) = 1+ r(\widehat \eta_0) +
d_1(\alpha)=1+d_1(\alpha),$$
since $\dim \Aut_{\widehat S^\alpha}\widehat \eta_0 =0$.

Next we show that $d_1(\alpha)$ has the following minimality property:

\begin{lemma}\label{d2} If $\alpha$ is special, then $d_1(\alpha) = r+1$. If
$\beta$ is not special, then $d_1(\beta) > r+1$.
\end{lemma}
\begin{proof} To see that $d_1(\alpha) = r+1$, it suffices by
  Proposition~\ref{d1}
to show that $\lambda_1(\alpha\spcheck) = \sum _{\beta\in \Delta
  }\beta\spcheck$.
First, by~\cite[cor.\ 3, p.\ 160]{Bour}, $\lambda = \sum _{\beta\in \Delta
}\beta\spcheck$ is always a coroot. Clearly, $n(\lambda\spcheck,
\beta\spcheck) \geq 0$ for all
$\beta \in \Delta -\{\alpha\}$, $n(\lambda\spcheck, \beta\spcheck) > 0$ if
$\beta\neq
\alpha$ is an end of the Dynkin diagram and
$\lambda - \beta\spcheck$ can only be a coroot if $\beta\neq \alpha$ is an
end of the Dynkin diagram. These properties say that, if $\beta\neq
\alpha$, then
$\lambda+\beta\spcheck$ is not a coroot. Thus $\lambda =
\lambda_1(\alpha\spcheck)$.

Now suppose that $R$ is not of $A$-type, and hence that $R\spcheck$ is
not of $A$-type. Then
$\lambda_1(\alpha\spcheck)$ is not the highest coroot of
$R\spcheck$. Thus there
exists a simple root $\beta$ such that $\lambda_1(\alpha\spcheck) +
\beta\spcheck$
is again a coroot. By what we have just seen, we must have $\beta =
\alpha$. It
follows that, for $\alpha \neq \beta$, $\lambda_1(\beta\spcheck)$ is
equal to
$\lambda_1(\alpha\spcheck)+\alpha\spcheck$ plus a sum of simple
coroots. Hence
$d_1(\beta) =  \rho( \lambda_1(\beta\spcheck) ) + 1 \geq
\rho( \lambda_1(\alpha\spcheck) ) + 2  = r+2$.
\end{proof}

Thus, we have proved the first statement in the theorem. To see the
second, first
assume that the Harder-Narasimhan parabolic subgroup for $\xi$ is maximal. In this
case, we can assume that the Harder-Narasimhan parabolic for $\xi$ is
$P^\beta_-$ for some
$\beta$. Thus
$\xi$ is isomorphic to $\eta\times _{L^\beta}G$, where $\eta$ is a
semistable bundle of negative degree $-n$ on $L^\beta$. Now
$$\dim H^1(E; \ad
(\eta\times_{L^\beta}G)) = 1+ r(\widehat \eta ) + nd_1(\beta) \geq r+2,$$
with equality holding if and only if $r(\widehat \eta) =0$, $n=1$, and
$\beta$ is
special. In this last case, it follows from Definition~\ref{translate} that
$\eta$ is a translate of $\eta_0$.

Now suppose that the Harder-Narasimhan parabolic   for $\xi$ is not
maximal. There exists a maximal parabolic
subgroup $P^\beta$ such that $\xi$ has a reduction to an $L^\beta$-bundle $\eta$,
where $\eta$ has degree $-n<0$.  By Lemma~\ref{list},
$\eta$ is unstable, for otherwise the Harder-Narasimhan parabolic for
$\xi$ would be
$P^\beta_-$, which is maximal. Hence, the associated $\widehat S^\beta$-bundle
$\widehat
\eta$ is also unstable. Thus the vector bundle $\ad_{\widehat
S^\beta}\widehat
\eta$ is unstable of degree zero, and hence contains a semistable summand of
negative degree. It follows that
$\dim H^1(E;
\ad_{\widehat S^\beta}\widehat
\eta)
\geq 1$. Applying Lemmas~\ref{adL} and \ref{uplus} and
Proposition~\ref{degrees},
we see that
\begin{align*}
h^1(E; \ad \xi) &= 1+ h^1(E; \ad_{\widehat S^\beta}\widehat
\eta)+h^1(E; \frak u (\eta)) + h^1(E; \frak u_- (\eta))\\
&\geq 2 + nd_1(\beta) \geq r+3.
\end{align*}
This completes the proof in case the Harder-Narasimhan parabolic   for
$\xi$ is not  maximal.
\end{proof}

We turn now to the $\Cee^*$-weights for the action of $\Cee^*$ on $\frak
u (\eta_0)$.

\begin{proposition}\label{1weights} Suppose that $\alpha$ is special. Then the
$\Cee^*$-weights for the action of the center of $L^{\alpha}$ on
$H^1(E; \frak u
(\eta_0))$, with multiplicity, are the integers $n_\alpha g_\beta, \beta \in
\widetilde \Delta$.
\end{proposition}
\begin{proof} The group $\Cee^*$ acts on $H^1(E;{\frak
u}^k (\eta))$ with weight $kn_\alpha$.  By Lemma~\ref{weights},
$$\dim H^1(E;{\frak
u}^k (\eta)) = - (\deg
\eta_0)\cdot i(\alpha, k)= i(\alpha, k).$$ Thus, it suffices to show that
$i(\alpha, k)=\#\{\beta \in \widetilde\Delta: g_\beta = k\}$. If we define
\begin{align*}
i(k) &= \#\{\beta \in
\widetilde\Delta: g_\beta = k\};\\
d(k) &= \sum _{k|x}i(x) = \sum _{\ell \geq 1}i(\ell k),
\end{align*}
then it clearly suffices to show that, for all $k$,  $d(k) = d_k(\alpha)$, in
the notation of Definition~\ref{defdk}. By Proposition~\ref{circular},
the integers
$d_k(\alpha)$ have the circular symmetry property with respect to
$g_\alpha$ and
$g$, since $\alpha$ is a long root, and by Lemma~\ref{d2}, $d_1(\alpha)=r+1$.
By the proof of Theorem 3.8.7 in~\cite{BFM}, the integers $d(k)$ have
the circular
symmetry property with respect to $N$ and $g$, where $N =
\max\{g_\beta:\beta \in
\Delta\}$. By Corollary 6.2.5 of \cite{FMAB}, $N=g_\alpha$. (We will give
another proof of this fact in Part III.) Clearly
$d(1)
=\#\widetilde \Delta = r+1$. By Lemma~\ref{circdetermines},  $d(k) = d_k(\alpha)$
for all
$k$.
\end{proof}

\subsection{The non-simply connected case}\label{ns}

We now establish the analogues in the non-simply connected
case of the previous results. While we
believe there should be classification-free arguments for these results,
we argue here in a case-by-case analysis.

\begin{defn}
Let $o(c)$ denote the order of $c\in \pi_1(G)$. A
root $\alpha\in \Delta$ is $c$-special if there exists an integer $d<0$
such that
\begin{enumerate}
\item[\rm (i)] $d/ o_{c, \alpha} \equiv \varpi_\alpha(c) \mod \Zee$;
\item[\rm (ii)] The Weyl orbit of the point $\mu_{c, \alpha} =
  d\zeta_\alpha/o_{c,\alpha}m_\alpha$ is minimal in the
Atiyah-Bott ordering \cite{FMAB} among all Weyl orbits of nonzero
points of Atiyah-Bott type for $c$.
\end{enumerate}
\end{defn}

\begin{remark}
a) The first condition means that that there is  a holomorphic
semistable
$ L^\alpha$-bundle $\eta$  with $c_1(\eta\times_{ L^\alpha}
  G)=c$ and Atiyah-Bott point equal to $\mu_{c,\alpha}$, and the
topological
type of the $ L^\alpha$-bundle
$\eta$ is uniquely determined by $\mu_{c,\alpha}$. This topological type
is
specified by an element $c_1(\eta) =\gamma \in \pi_1( L^{\alpha})$. The
second condition means that the point
$\mu_{c,
\alpha}$, or the corresponding stratum of
$(0,1)$-connections,
is   minimally unstable in the sense of \cite[Definition 6.1.1]{FMAB}.

b)
It is not always true that there is a unique $c$-special root.
Uniqueness fails exactly when  $ G=SL_n(\Cee)/\langle c\rangle$ and $c$
does not
generate the center, where there are $n/o(c)$ special roots, and for
$ G=SO(2n)$, where there are two special roots.

c) As defined here, $c$-special roots are certain simple roots.
In \cite{FMAB} we   used the simple roots to index
strata of the space of $(0,1)$-connections, by associating to
$\alpha$ the stratum  lying in the Lie algebra
of the center of $L^\alpha$ and having the smallest possible {\bf
  positive} value under the dominant character.
The convention here differs by a sign from the one of \cite{FMAB}. This
means that the image under the automorphism of the Dynkin diagram
induced by $-w_0$ of the roots which are $c$-special as defined here
 correspond to the roots indexing the minimally unstable strata in
\cite{FMAB}. This automorphism sends $c$ to $c^{-1}$ and thus fixes the
set of
$c$-special roots if and only if $c$ is of order $1$ or $2$.
Otherwise, the $c$-special roots correspond to the roots indexing the
minimally unstable strata for $c^{-1}$ in \cite{FMAB}.
\end{remark}

\begin{theorem}\label{cthm}
Let $\alpha$ be a $c$-special root for $G$, and let
$\gamma\in
\pi_1( L^\alpha)$ be the first Chern class of the $
L^\alpha$-bundle $\eta$ corresponding to $\mu_{c,\alpha}$. Then:
\begin{enumerate}
\item[\rm (i)] The integer $d=-1$ and $o_{c, \alpha} =o(c)$.
\item[\rm (ii)] The adjoint quotient $\ad( L^\alpha)=  L^\alpha/Z(
 L^\alpha)$ is a product $
\prod_{i=1}^k\widehat S_i$,
where the $\widehat S_i$ are simple groups of $A$-type.
\item[\rm (iii)]
Let $\widehat \gamma$ be the image of $\gamma\in \pi_1( L^\alpha)$
under
the projection
$\pi_1( L^\alpha)\to \pi_1(\prod_{i=1}^k\widehat
S_i)=\prod_{i=1}^k\pi_1(\widehat S_i)$. For $i=1,\ldots, k$, the image
of
$\widehat \gamma$ in $\pi_1(\widehat S_i)$ generates the cyclic group
$\pi_1(\widehat S_i)$.
\item[\rm (iv)] $d_1(\alpha)/o_{c, \alpha}=r_c+1$.
\end{enumerate}
\end{theorem}

\begin{proof} The minimally unstable points $\mu_{c, \alpha}$ are listed
in \S
6.3 of \cite{FMAB}. From this list, it is easy to check that $d=-1$ and
$o_{c, \alpha} =o(c)$. To prove the remaining statements, we make a
case-by-case analysis.

\medskip
\noindent $\widetilde G=SL_n(\Cee)$:

We choose an identification of $\Lambda$ with $\{(x_1,\ldots,x_n)\in
\Zee^n\bigl|\bigr. \sum_{i=1}^nx_i=0\}$ in such a way that
$\Delta=\{\alpha_1,\ldots,\alpha_{n-1}\}$ with $\alpha_i=e_i-e_{i+1}$
with $e_i$ being the standard unit vector in the $i^{\rm th}$-coordinate
direction.
We write $c$ as the image of an element of the form
  $(m/n)(\alpha_1\spcheck+2\alpha_2\spcheck+\cdots
  +(n-1)\alpha_{n-1}\spcheck)$ for some $1\le m< n$.
We factor $m=m_0\cdot \ell$ where $(m_0,n)=1$ and $\ell|n$. Then $c$ is
an
  element of order $f=n/\ell$.
By \cite[\S 6.3]{FMAB},
the $c$-special roots are those $\alpha$ for which
  $\varpi_\alpha(c)\equiv -1/f \pmod \Zee$. There are exactly $n/f$ such
  roots. Suppose $\alpha=\alpha_k$.
Then $-k\cdot m/n\equiv -1/f \pmod
\Zee$. In
  particular, $(k,f)=1$.
Since   $\alpha$ is a special root for $SL_n(\Cee)$,
  $d_1(\alpha)=n$, and hence
  $d_1(\alpha)/f=n/f=\ell$.
On the other hand, since $\langle c\rangle$ acts freely on the Dynkin
  diagram for $G$ we see that $r_c+1=\ell$.
Let us consider the group $ L^\alpha$. By Lemma~\ref{Lalpha2},
$\widetilde L^\alpha$ is isomorphic to the subgroup of $GL_k(\Cee)\times
GL_{n-k}(\Cee)$  matrices of   equal determinant.
Hence $\ad( L^\alpha)=\ad(\widetilde L^\alpha)=PGL_k(\Cee)\times PGL_{n-k}(\Cee)$.
The map $\det\colon \widetilde L^\alpha\to \Cee^*$ induces an
  identification $\pi_1(\widetilde L^\alpha)=\Zee$ and the projection
  $\widetilde L^\alpha
  \to \ad(L^\alpha)$ sends $1\in \Zee$ to  the element $(a_k,a_{n-k})\in
  \pi_1(PGL_k(\Cee))\times \pi_1(PGL_{n-k}(\Cee))$ where $a_k$, resp. $a_{n-k}$
  generates $\pi_1(PGL_k(\Cee))$, resp. $\pi_1(PGL_{n-k}(\Cee))$.
Direct computation
  shows that $\pi_1( L^\alpha)=\Zee$ and that the natural map
  $\pi_1(\widetilde L^\alpha) \to\pi_1( L^\alpha)$ is multiplication by $f$.
Thus, we have an identification $\pi_1( L^\alpha)=\Zee[1/f]$.
Under this identification
 the element $\gamma\in \pi_1( L^\alpha)$ is $-1/f$.
  $\pi_1( L^\alpha)$ and projects to
  $\widehat \gamma=(f^{-1}a_k,f^{-1}a_{n-k})\in \pi_1(PGL_k(\Cee))\times
  \pi_1(PGL_{n-k}(\Cee))$.
Since $(k,f)=1$, the projection of this element to either factor
generates that factor.

\medskip
\noindent $ G=SO(2n+1)$, $n\ge 3$.

The $c$-special root is the unique short root $\alpha$ in the Dynkin
diagram. Direct inspection shows that $ L^\alpha=GL_n(\Cee)$, and
that $\gamma\in \pi_1( L^\alpha)$ is a generator of the fundamental
group. Thus, $\ad( L^\alpha)=PGL_n(\Cee)$ and $\widehat\gamma$ generates
$\pi_1(PGL_n(\Cee))$.
Furthermore, $d_1(\alpha)=2n$ so that $d_1(\alpha)/o(c)=n$.
Since $c$ acts on the extended Dynkin diagram for $G$
with one free orbit and $n-1$ fixed points, we see that
$r_c+1=n=d_1(\alpha)/o(c)$.

\medskip
\noindent $ G=Sp(2n)/\langle c\rangle$, $n\ge 2$:

Suppose first that $n$ is odd.
 Then there is a unique $c$-special root, the unique
long root $\alpha$. In this case, $ L^\alpha=GL_n(\Cee)/(\Zee/2\Zee)$
and
$\gamma$ generates the fundamental group of $ L^\alpha$.
Hence, $\ad( L^\alpha)=PGL_n(\Cee)$ and $\widehat\gamma$ is the square of a
 generator for this group. Since $n$ is odd, $\widehat \gamma$ is a
generator of $\pi_1(PGL_n(\Cee))$.
Since $\alpha$ is special for the simply connected form of the
 group, $d_1(\alpha)=n+1$. In this case the element $c$ acts freely on
 the nodes of the extended Dynkin diagram so that
 $r_c+1=(n+1)/2=d_1(\alpha)/o(c)$.

Now suppose that $n$ is even. Then there is a unique
$c$-special root,  the unique short
simple root $\alpha$ which is not orthogonal to the unique long simple
root. Direct computation shows that  $\widetilde L^\alpha$
is isomorphic to $GL_{n-1}(\Cee)\times
SL_2(\Cee)$ and that $c$ is the diagonal element $(-1,-1)$. Thus,
$\pi_1( L^\alpha)=\Zee$ and $\gamma$ is a generator of this group.
Furthermore, $\ad( L^\alpha)=PGL_{n-1}(\Cee)\times PGL_2(\Cee)$ and the map
$\pi_1( L^\alpha)\to \pi_1(\ad( L^\alpha))$ is onto. Thus, the
image $\widehat \gamma$ of $\gamma$ generates $\pi_1(\ad( L^\alpha))$,
and hence its projection to each factor generates the fundamental
group of that factor.
Lastly, direct computation shows that $d_1(\alpha)=n+2$. Since $c$
acts on the extended Dynkin diagram for $G$ with one fixed point and
$n/2$ free orbits, we see that
$r_c+1=(n+2)/2=d_1(\alpha)/o(c)$.

\medskip
\noindent $ G=SO(2n)$, $n\ge 4$:

For $\widetilde G=Spin(2n)$ we identify $\Lambda$ with the even integral lattice
inside $\Ar^n$. Let $e_i$ be the standard unit vector in the
$i^{\rm th}$-coordinate direction.  Then
$\Delta=\{\alpha_1,\cdots,\alpha_{n-1},\alpha_n\}$ where
$\alpha_i=e_i-e_{i+1}$ for $1\le i<n$ and $\alpha_n=e_{n-1}+e_n$.

There are two $c$-special roots $\alpha_{n-1} = e_{n-1}-e_n$
and $\alpha_n = e_{n-1}+e_n$. (Of course, these elements are
interchanged by
an outer automorphism of
$SO(2n)$.) Let $\alpha$ be one of the $c$-special roots.
Then $ L^\alpha=GL_n(\Cee)$ and $\gamma$ is a generator of $\pi_1(
L^\alpha)\cong \Zee$. Thus, $\ad( L^\alpha)=PGL_n(\Cee)$ and $\widehat
\gamma$ is a generator of this group.
Direct computation shows that $d_1(\alpha)=2(n-1)$. Since $c$ acts on
the Dynkin diagram for $G$ with two free orbits and $n-1$ fixed
points, we see that $r_c+1=n-1=d_1(\alpha)/o(c)$.

\medskip
\noindent $\widetilde G=Spin(4n+2)$, $n\ge 2$  and $c$ is an element of
  order $4$:

There is  one $c$-special root. It is the simple
root $\alpha$ corresponding to the  ``ear'' of the Dynkin diagram (i.e.
either $\alpha_{n-1}$ or $\alpha_n$) with the property that
$\varpi_\alpha(c)=-1/4\pmod\Zee$. In this case
$ L^\alpha=GL_{2n+1}(\Cee)/(\Zee/2\Zee)$. Hence $\pi_1(
L^\alpha)=\Zee$ and  $\gamma$ is a generator. Under the projection
to $\ad( L^\alpha)=PGL_{2n+1}(\Cee)$ the image $\widehat \gamma$ of $\gamma$
is the square of  the usual generator. This is clearly still a
generator.
Lastly, as above $d_1(\alpha)=2((2n+1)-1)=4n$ whereas $r_c+1=n$.
Thus, $d_1(\alpha)/o(c)=r_c+1$.

\medskip
\noindent $\widetilde G=Spin(4n)$, $n\ge 2$ and $c$ is  an element of
  order two not contained in  $\pi_1(SO(4n))$:

There is  one $c$-special root. It is the simple
root $\alpha_{n-3}$ corresponding to the node of the ``long'' arm of the
Dynkin diagram next to the trivalent node. Thus,
$\widetilde L^\alpha$  is
isomorphic to $(SL_{2n-3}(\Cee)\times
SL_4(\Cee))\times_{(\Zee/(4n-6)\Zee)}\Cee^*$
where the cyclic group is embedded in the standard way in $\Cee^*$ and
the
usual generator maps to the standard generator of $Z(SL_{2n-3}(\Cee))$
and to
the element of order $2$ in $Z(SL_4(\Cee))$. Thus, we can identify
$\widetilde L^\alpha$
with
$GL_{2n-3}(\Cee)\times_{(\Zee/2\Zee)}SL_4(\Cee)$.
The element $c$ is the image of the element
$(a,b)$ where $a$ and $b$ are central elements of order $4$ in
$GL_{2n-3}(\Cee)$ and $SL_4(\Cee)$   under the inclusion of
$\widetilde L^\alpha\subset
\widetilde G$. Thus,
$ L^\alpha$ is isomorphic
$GL_{2n-3}(\Cee)\times_{(\Zee/4\Zee)}SL_4(\Cee)$.
Hence $\gamma$ is a generator of $\pi_1( L^\alpha)=\Zee$, the
image  $\widehat\gamma$ of $\gamma$ is a generator for
$\pi_1(PGL_{2n-3}(\Cee))\times\pi_1(PGL_4(\Cee))$, and hence the projection of
$\widehat\gamma$ into either factor generates the fundamental group
of that factor.
Direct computation shows that $d_1(\alpha)=2(n+1)$. Since the action
of $c$ on the extended Dynkin diagram for $G$ has one fixed point and
$n$ free orbits, we see that $r_c+1=n+1=d_1(\alpha)/o(c)$.

\medskip
\noindent $ G=\ad(E_6)$:

There is  one $c$-special root. It is a simple
root $\alpha$ corresponding to the node next to the trivalent node on
one of
the arms of length $3$ with the property $\varpi_\alpha(c)=-1/3$.
In this case $\widetilde L^\alpha$ is isomorphic to $(SL_5(\Cee)\times
SL_2(\Cee))\times_{(\Zee/10\Zee)}\Cee^*$ where the element in
$\Zee/10\Zee$
that maps to ${\exp }(2\pi i/10)$ maps to the generator in
$Z(SL_2(\Cee))$
and to the square of the usual generator in $Z(SL_5(\Cee))$. Hence, $
L^\alpha$ is isomorphic to
$(SL_5(\Cee)\times SL_2(\Cee))\times_{(\Zee/10\Zee)}\Cee^*$
where the element in $\Zee/10\Zee$ that maps to ${\exp }(2\pi i/10)$
maps to the generator in $Z(SL_2(\Cee))$ and to the usual
generator in $Z(SL_5(\Cee))$. Thus, $\ad( L^\alpha)=PGL_2(\Cee)\times
PGL_5(\Cee)$,
$\pi_1( L^\alpha)$ is cyclic and
$\gamma$ is a generator of this group.
It follows that   $\widehat\gamma\in \pi_1(PGL_5(\Cee)\times
PGL_2(\Cee))$ generates and hence the image of $\widehat\gamma$ under
projection to either factor is  a generator of the fundamental group
of that  factor.
Direct computation shows that $d_1(\alpha)=9$. Since the action of
$\langle c\rangle$ on the extended Dynkin diagram of $E_6$ has two
free orbits and one fixed point, we see that
$r_c+1=3=d_1(\alpha)/o(c)$.

\medskip
\noindent $ G=\ad(E_7)$:

There is one $c$-special root. It corresponds to
the node of the Dynkin diagram next to the trivalent node on the arm
of length $4$. In this case  $\widetilde L^\alpha$ is isomorphic to
$(SL_3(\Cee)\times
SL_5(\Cee))\times_{(\Zee/15\Zee)} \Cee^*$, where the element in
$\Zee/15\Zee$
that maps to ${\exp }(2\pi i/15)\in \Cee^*$ maps to the usual
generator of $Z(SL_3(\Cee))$ and the square of the usual generator of
$Z(SL_5(\Cee))$.  Thus,
$ L^\alpha=(SL_3(\Cee)\times
SL_5(\Cee))\times_{(\Zee/15\Zee)} \Cee^*$ where the element in
$\Zee/15\Zee$
that maps to ${\exp }(2\pi i/15)\in \Cee^*$ maps to the inverse of
the usual
generator of $Z(SL_3(\Cee))$ and the inverse of the usual generator of
$Z(SL_5(\Cee))$.
Thus, $\ad( L^\alpha)=PGL_3(\Cee)\times PGL_5(\Cee)$, $\pi_1( L^\alpha)$ is
isomorphic to $\Zee$ and $\gamma$ is a
generator. Consequently, $\widehat\gamma$ is a generator of $\pi_1(\ad(
L^\alpha))$.

Direct computation shows that $d_1(\alpha)=10$. The action of $c$ on
the extended Dynkin diagram of $E_7$ has two fixed points and $3$ free
orbits so that $r_c+1=5=d_1(\alpha)/o(c)$.
\end{proof}

Next we compute the integer $n_{c, \alpha}$ defined in Lemma~\ref{ncalpha}:

\begin{lemma} If $\alpha$ is $c$-special and $c$ is nontrivial, then
$n_{c,\alpha} = 1$ except in the following cases:
\begin{enumerate}
\item[\rm (i)] If $\widetilde G = SL_n(\Cee)$, $\alpha$ corresponds to the $k^{\rm
th}$
vertex in the usual ordering,  and
$o(c) = d$, then $n_{c, \alpha} = n/d\cdot \gcd(k,n)$.
\item[\rm (ii)] If $\widetilde G = Spin (2n)$ and $c$ is of order $2$, then
  $n_{c, \alpha} = 2$.
\end{enumerate}
\end{lemma}
\begin{proof} If $Z(\widetilde G)$ is cyclic and $c$ is a generator, then $n_{c,
\alpha} = 1$. The remaining cases are $\widetilde G=SL_n(\Cee)$ and
$\widetilde G = Spin (2n)$, and these can
be checked directly.
\end{proof}

As in the simply connected case, we have:

\begin{lemma}\label{ceta0} There is a unique semistable $
L^\alpha$-bundle $\eta_0$ with the following properties:
\begin{enumerate}
\item[\rm (i)] $c_1(\eta_0\times _{L^\alpha} G) = c$;
\item[\rm (ii)] The Atiyah-Bott point of $\eta_0$ is $\mu_{c, \alpha}$;
\item[\rm (iii)] $\det \eta_0 = \scrO_E(-p_0)$. \endproof
\end{enumerate}
\end{lemma}

As before, a bundle $\eta$ satisfying (i) and (ii) above
is the pullback of $\eta_0$ via a translation of $E$, and we will call such an
$\eta$ a {\sl translate\/} of $\eta_0$. The following is then proved via
arguments
similar to those used in the proof of Theorem~\ref{mindim}.

\begin{theorem}\label{cmindim} Let $\eta_0$ be the bundle described in
Lemma~\ref{ceta0}. Then
$$\dim H^0(E; \ad (\eta_0\times_{ L^{\alpha}} G)) = \dim H^1(E;
\ad
(\eta_0\times_{ L^{\alpha}} G)) = r_c+2.$$
If $\xi$ is any unstable $ G$-bundle which is $C^\infty$ isomorphic
to
$\eta_0\times_{ L^{\alpha}} G$, then
$\dim H^1(E; \ad \xi)
\geq r_c+2$, with equality if and only if $\xi$ is isomorphic to
$\eta\times
_{ L^{\alpha}} G$ for some   translate
$\eta$  of $\eta_0$.\endproof
\end{theorem}

Finally, we must determine the  weights for the action of $\Cee^*$ on
$H^1(E;\frak u(\eta_0))$:

\begin{proposition}\label{cweights} Suppose that $\alpha$ is
$c$-special. Then the
$\Cee^*$ weights for the action of $\ov \varphi_\alpha(\Cee^*)$ on
$H^1(E; \frak u
(\eta_0))$, with multiplicity, are the integers $n_{c,\alpha}   g_{\ov
\beta}/n_0,
\ov \beta \in
\widetilde \Delta/w_c$.
\end{proposition}
\begin{proof} The weight for the action of $\Cee^*$ on $H^1(E; \frak
u^k(\eta_0))$ is $kn_{c, \alpha}$. By Lemma~\ref{weights}, the dimension
of
$H^1(E; \frak u^k(\eta_0))$ is $i(\alpha, k)/o_{c, \alpha} = i(\alpha,
k)/o(c)$. By Proposition~\ref{circular}, the integers
$d_k(\alpha)/o(c)$ have the circular symmetry property with respect
to
$h_\alpha$ and
$gh_\alpha/g_{\alpha} o(c)$. By Theorem~\ref{cthm}, $d_1(\alpha)/o(c) =
r_c+1$. If we define
\begin{align*}
i_c(k) &= \#\{\beta \in
\widetilde\Delta: g_{\ov  \beta }= kn_0\};\\
d_c(k) &= \sum _{k|x}i_c(x) = \sum _{\ell \geq 1}i_c(\ell k),
\end{align*}
then the integers $d_c(k)$ satisfy: $d_c (1) = r_c+1$, and the $d_c(k)$
have
the circular symmetry property with respect to $N/n_0$ and $g/n_0$,
where $N$
is the maximum value of  the $g_{\ov \beta}$. It follows by inspection
or
from \cite[Proposition 10.1.8]{BFM} that $n_0 =
o(c)g_\alpha/h_\alpha$, i.e. $h_\alpha/g_{\alpha} o(c) = 1/n_0$. By
inspection, $h_\alpha=N/n_0$. Thus $i(\alpha, k)/o(c) = i_c(k)$, and the
proof follows.
\end{proof}

\section{The nonabelian cohomology space}

\subsection{The affine space and the universal bundle}

Let $\alpha$ be an arbitrary simple root. We abbreviate $P^\alpha= P$, $L^\alpha =
L$, and $U^\alpha = U$.   Let
$\eta$ be an unliftable  semistable principal
$L$-bundle of type $c$.
There is the associated sheaf of (not necessarily
abelian) groups
$U(\eta)$. The cohomology set $H^1(E; U(\eta))$ classifies
pairs $(\xi, \varphi)$, where $\xi$ is a $ P$-bundle and $\varphi$ is an
isomorphism from the induced bundle $\xi/U$ on $ L$ to $\eta$. There is a
marked point
$0\in H^1(E; U(\eta))$, corresponding to the pair
$(\eta\times _{ L} P, I)$, where $I$ is the canonical identification of the
bundle 
$(\eta\times _{ L} P)/U$ with $\eta$. There is a corresponding functor $\mathbf{F}$
from schemes to sets defined as follows: for a scheme of finite type over $\Cee$,
$\mathbf{F}(S)$ is the set of isomorphism classes of pairs $(\Xi, \Phi)$, where
$\Xi$ is a
$ P$-bundle over $E\times S$ and $\Phi$ is an isomorphism from
$\Xi/U$ to $\pi_1^*\eta$. 

\begin{lemma} The functor $\mathbf{F}$ is represented by an affine space.
\end{lemma}
\begin{proof} Let $U_i$ be the closed subgroup of $U$ whose Lie algebra is
$\bigoplus _{k\geq i}\frak u^k$.  Then the filtration $\{U_i\}$ is a decreasing
filtration of $U$ by normal, $L$-invariant subgroups such that $U_i/U_{i+1}$ is
in the center of $U/U_{i+1}$ for every $i$, and $U_i/U_{i+1} \cong \frak u^i$. By
Theorem~\ref{A6} of the appendix, it suffices to check that 
$H^0(E;(U_i/U_{i+1})(\eta)) = H^2(E;(U_i/U_{i+1})(\eta))  =0$. The second statement
is clear since $\dim E =1$, and the first follows from Lemma~\ref{uplus}, which
implies that $H^0(E;\frak u^k(\eta) ) =0$ for every $k>0$.
\end{proof}

Thus, there is a structure of an
affine space on
$H^1(E; U(\eta))$ and a universal pair
$(\Xi_0, \Phi_0)$ over the scheme $E\times H^1(E; U(\eta))$ which
represents the functor $\mathbf{F}$. We will somewhat carelessly identify $\Xi_0$
with the associated $G$-bundle $\Xi_0\times _{ P} G$.

We now identify $\varphi _\alpha(\Cee^*)$ with $\Cee^*$. Thus we have fixed
the embedding of $\Cee^*$ in $L$. Since $L$ acts on $U$, there are induced actions
of $\Cee^*$ on
$U(\eta)$ and $\frak u(\eta)$, and hence on $H^1(E; U(\eta))$ and on $H^1(E;
\frak u(\eta))$.   Viewing $H^1(E;
U(\eta))$ as the set of pairs $(\xi, \varphi)$ as above, the action of $\Cee^*$
is via the action of $\Aut L$ on the isomorphism
$\varphi$, and this action  fixes the origin in $H^1(E; U(\eta))$, i.e.\ the
bundle $\eta\times_LP$. By Theorem~\ref{A6}, the action of $\Cee^*$ lifts to an
action on the universal principal bundle
$\Xi_0$ over $E\times H^1(E; U(\eta))$. The first goal of this section is to prove
that the action of $\varphi _\alpha(\Cee^*)$ on $H^1(E; U(\eta))$ is linearizable,
and in fact there is a $\Cee^*$-equivariant isomorphism from $H^1(E; U(\eta))$ to
$H^1(E; \frak u (\eta))$. Thus the quotient $(H^1(E; U(\eta))-\{0\})/\Cee^*$ is a
weighted projective space $\WP(\eta)$. In \S\ref{univbund}, we give a sufficient
condition for the existence of universal bundles over $E\times \WP(\eta)$.  Next we
show that, in the case where
$\alpha$ is
$c$-special and $\eta_0$ is the bundle of Proposition~\ref{eta0} or
Lemma~\ref{ceta0}, the points of
$\WP(\eta_0)$ correspond to semistable bundles whose automorphism groups have
minimal possible dimensions. In
\S\ref{sKS}, we analyze the Kodaira-Spencer homomorphism. The results of
\S\ref{sKS} will not however be used in this paper. Finally, we discuss the
singular locus of the weighted projective space and relate it to moduli spaces of
bundles with a non-simply connected structure group.

\subsection{Linearization of the action}

  We analyze the $\Cee^*$-action
on the affine space $H^1(E; U(\eta))$ more closely. Our goal will be to show that
this action can be linearized and to calculate the $\Cee^*$-weights.

\begin{lemma}\label{4.1} Every $\Cee^*$-orbit in $H^1(E; U(\eta))$ contains the
origin in its closure.
\end{lemma}
\begin{proof} Let $x\in H^1(E; U(\eta))$ be represented by the $1$-cocycle
$\{u_{ij}\}$, where $\{\Omega _i\}$ is an open cover of $E$ and   $u_{ij}\colon
\Omega _i\cap \Omega _j \to U$ is a morphism. Then $\lambda\in \Cee^*\subseteq L$
acts on the cocycle $\{u_{ij}\}$.    Define morphisms
$$\tilde u_{ij}(e, \lambda)\colon (\Omega _i\cap \Omega _j) \times  \Cee\to
U$$ as follows:
$$\tilde u_{ij}(e, \lambda) = \begin{cases} \lambda \cdot u_{ij}(e), &\text{if
$\lambda
\neq 0$;}\\
1, &\text{if $\lambda = 0$.}
\end{cases}$$
There is a
$\Cee^*$-equivariant morphism from   the unipotent
subgroup $U$ to the affine space $\frak u$    (see for example \cite[Remark, p.\
183]{Borel}). Using this $\Cee^*$-equivariant isomorphism, and
the fact that all of the $\Cee^*$-weights on $\frak u$ are positive, it is easy
to check that the $\tilde u_{ij}$ are morphisms and so define a $1$-cocycle for the
sheaf $U(\pi_1^*\eta)$ over $E\times \Cee$. Thus, they define a bundle
$\Xi$ over $E\times \Cee$, reducing to $\pi _1^*\eta$ mod $U$ and such
that $\Xi|E\times \{0\} = \eta \times _{ L} P$. By the functorial property of
$H^1(E; U(\eta))$, there is a morphism from
$\Cee$ to $H^1(E; U(\eta))$ corresponding to $\Xi$. Clearly, the image
of $0\in \Cee$ is the origin of $H^1(E; U(\eta))$, and the image of
$\Cee^*$ is exactly the $\Cee^*$-orbit of the cocycle $\{u_{ij}\}$. This proves
Lemma~\ref{4.1}.
\end{proof}

\begin{lemma}\label{4.2} Let $T$ be the tangent space of $H^1(E; U(\eta))$ at the
fixed point $0$, with the natural $\Cee^*$-action. Then  the Kodaira-Spencer
homomorphism from $T$ to $H^1(E;
\ad (\eta \times _{ L} G))$ induced by the bundle $\Xi_0$ is given by
a $\Cee^*$-equivariant isomorphism  $T\to H^1(E; \frak{u}(\eta))$ followed by the
inclusion of
$H^1(E; \frak{u}(\eta))$ as a direct summand in
$$H^1(E;\ad(\eta\times_{ L} G))\cong H^1(E;\ad_{ L}\eta)\oplus
H^1(E;{\frak u}_-(\eta))\oplus H^1(E;{\frak u}(\eta)).$$
\end{lemma}
\begin{proof} Let $\Cee[\varepsilon]$ denote the dual numbers. The space $T$ is the
set of maps from $\Spec \Cee[\varepsilon]$ to $H^1(E; U(\eta))$ such that
the closed point is mapped to the origin. By the functorial interpretation of
$H^1(E; U(\eta))$, such a morphism corresponds to a $ P$-bundle
$\Xi$ over
$E\times \Spec \Cee[\varepsilon]$, which is the pullback of $\Xi_0$, together with
an isomorphism from
$\Xi/U$ to
$\pi _1^*\eta$, and such that $\Xi$ restricts to $\eta \times _{ L} P$ over
the closed point. The second condition says that $\Xi$ is a first order
deformation of the
$ P$-bundle $\eta \times _{ L} P$. Such deformations are classified by
$$H^1(E; \ad_P (\eta \times _{ L} P)) = H^1(E; \ad_L \eta) \oplus H^1(E; \frak
u(\eta)).$$
The first condition says that the corresponding first order deformation of
the $ L$-bundle $\eta$ is trivial, or equivalently that the projection of the
Kodaira-Spencer class of $\Xi$ to $H^1(E; \ad_L \eta)$ is zero. Thus we have
defined a canonical map from $T$ to $H^1(E; \frak{u}(\eta))$. Conversely, by
reversing this construction, every element of $H^1(E; \frak{u}(\eta))$ defines a
first order deformation of $\eta
\times _{ L} P$ which reduces to $\pi _1^*\eta$ mod
$U$, so that in fact the map from $T$ to $H^1(E; \frak{u}(\eta))$ is an
isomorphism. Since this isomorphism is canonical, it is easily seen to be
$\Cee^*$-equivariant. The last statement is clear by construction.
\end{proof}

Using the previous two lemmas, we show that the $\Cee^*$-action on $H^1(E;
U(\eta))$ can be linearized.  There is the following general result about
$\Cee^*$-actions on an affine space.

\begin{lemma}\label{4.3} Let $\mathbb{A}^n$ be an affine space with a 
$\Cee^*$-action, and suppose that $0\in \mathbb{A}^n$ is a fixed point for the
action. Let $T$ be the tangent space of $\mathbb{A}^n$ at the origin, together
with the induced linear $\Cee^*$-action on $T$. Further suppose that
\begin{enumerate}
\item[\rm (i)] Every $\Cee^*$-orbit in $\mathbb{A}^n$ contains $0$ in its closure.
\item[\rm (ii)] All of the weights in the $\Cee^*$-action on $T$ are
strictly positive.
\end{enumerate}
Then there is a $\Cee^*$-equivariant isomorphism from $\mathbb{A}^n$ to $T$. Hence,
the $\Cee^*$-action on
$\mathbb{A}^n$ is linearizable, and the $\Cee^*$-weights for this action are those
for the action on $T$.
\end{lemma}
\begin{proof} Let $A =\Cee[z_1, \dots, z_n]$ be  the affine coordinate ring of
$\mathbb{A}^n$, where $0$ is defined by $z_1=\cdots = z_n=0$, and let $x_1, \dots,
x_n$ be a basis for the linear functions on
$T$. The finite-dimensional subspace of   $A$ spanned by the $z_i$ is contained in
a finite-dimensional
$\Cee^*$-invariant subspace
$V$ of
$A$, by the Cartier lemma
\cite[p.\ 25]{GIT} (or by using the grading on $A$ induced by the $\Cee^*$-action).
The map
$p\in A
\mapsto (dp)_0$ is a
$\Cee^*$-equivariant map from $A$ to $T^*$, and hence restricts to a
$\Cee^*$-equivariant map from $V$ to $T^*$. Choosing a
$\Cee^*$-equivariant splitting of the map
$V\to T^*$ defines  a $\Cee^*$-equivariant map $T^*\to A$ and thus a
$\Cee^*$-equivariant homomorphism $\Cee[x_1, \dots, x_n] \to A$. Let $f: 
\mathbb{A}^n \to T$ be the corresponding morphism. By construction, $f$ has an
invertible differential at the origin  and is
$\Cee^*$-equivariant. Thus, $f$ is injective in a neighborhood
$\Omega$ of the origin, and the image of
$f$ contains an open set
$\Omega'$ about the origin. Since the weights on $T$ are positive,  every point of
$T$ lies in the $\Cee^*$-orbit of some point of $\Omega'$. Thus $f$ is surjective.
Likewise, $f$ is injective: if $f(x_1) = f(x_2)$, then since the closures of the
$\Cee^*$-orbits of $x_1$ and $x_2$ contain the origin, and the weights of the
action on the tangent space at $0$ are all positive, it follows that there is a
$\lambda\in \Cee^*$ such that $\lambda \cdot x_1$ and $\lambda \cdot x_2$ both lie
in $\Omega$. By assumption $f(\lambda \cdot x_1) = \lambda \cdot f(x_1) = \lambda
\cdot f(x_2) = f(\lambda \cdot x_2)$. But since $f$ is injective on $\Omega$,
$\lambda \cdot x_1 = \lambda \cdot x_2$, and hence $x_1 = x_2$. It follows that
$f$ is a $\Cee^*$-equivariant bijection from $\mathbb{A}^n$ to $T$ and thus it is 
an isomorphism.
\end{proof}

\begin{corollary}\label{4.4} The 
$\Cee^*$-equivariant  morphism from $H^1(E; U(\eta))$ to $H^1(E;
\frak{u}(\eta))$ defined in Lemma~\ref{4.2} is an isomorphism. Hence, the
$\Cee^*$-action on 
$H^1(E; U(\eta))$ can be  linearized, and 
$(H^1(E; U(\eta)) -\{0\})/\Cee^*$ is a weighted projective space $\WP(\eta)$.
\endproof
\end{corollary} 

In case $\alpha$ is $c$-special and $\eta =\eta_0$, we
have calculated the corresponding weights of the weighted projective space in
Proposition~\ref{1weights} and Proposition~\ref{cweights}.

\subsection{Existence of universal bundles on the weighted projective
space}\label{univbund}

Next we discuss the existence of universal bundles over the $\Cee^*$-quotient. It
is easy to see that such bundles cannot exist at the orbifold singular points of
the weighted projective space, essentially because  there are no
local sections from the weighted projective space back to the affine space at such
points. We shall show that, away from such points, we can almost  find a universal
bundle.  In particular, there is a universal adjoint bundle away from the orbifold
singular points of the weighted projective space. 

Recall that $Z(G)\cap \varphi _\alpha(\Cee^*)$ is a finite cyclic group which we
have denoted 
$\Zee/n_{c,\alpha}\Zee$.

\begin{lemma} The subgroup $\Zee/n_{c,\alpha}\Zee$ of $\Cee^*$ acts trivially
on $H^1(E; U(\eta))$, and the quotient group
$\Cee^*/(\Zee/n_{c,\alpha}\Zee)$ acts faithfully on 
$H^1(E; U(\eta))$.
\end{lemma}
\begin{proof} The $\Cee^*$-weights are of the form $kn_{c, \alpha}$ for $1\leq
k\leq h_\alpha$, and so the lemma is clear.
\end{proof}

The fact that $\Zee/n_{c,\alpha}\Zee$  acts trivially
on $H^1(E; U(\eta))$ also follows from the fact that it is contained in $Z(G)$.
Note however that, if
$n_{c,\alpha} >1$, then the associated action of $\Zee/n_{c,\alpha}\Zee$ on
$\Xi_0$ is not in general trivial, and in fact is just multiplication by the
corresponding subgroup of the center of $ G$.

\begin{proposition}\label{4.5} Let $\WP(\eta)$ be the weighted
projective space $(H^1(E; U(\eta))
-\{0\})/\Cee^*$,
and let $\WP^0(\eta)$ denote the open subset of
$\WP(\eta)$ which is the $\Cee^*$-quotient of the set of points of
$H^1(E; U(\eta))$ where $\Cee^*/(\Zee/n_{c,\alpha}\Zee)$ acts freely. Let
$\widehat G$ be the quotient of $G$ by the subgroup $\Zee/n_{c,\alpha}\Zee$. Then
the  universal bundle
$\Xi_0$ over $E\times H^1(E; U(\eta))$ induces a principal $\widehat
G$-bundle over $E\times \WP^0(\eta)$.
\end{proposition}
\begin{proof} Let $\left(H^1(E; U(\eta))\right)^0$ be the set of points of
$H^1(E; U(\eta))$ where $\Cee^*/(\Zee/n_{c,\alpha}\Zee)$ acts freely and
effectively. We have seen that there is a lifted action of 
$\Cee^*$ on $\Xi_0|E\times \left(H^1(E; U(\eta))\right)^0$, which in
fact is free. The action of the isotropy group $\Zee/n_{c,\alpha}\Zee$ of a point
in the base on the fiber is via multiplication by elements of the center of $G$,
and thus there is an induced $\widehat G$-bundle on the $\Cee^*$-quotient of
$\left(H^1(E; U(\eta))\right)^0$ with the desired properties.
\end{proof}

It is easy to see that we could replace $\widehat
G$-bundles with bundles over an appropriate conformal form $G\times
_{\Zee/n_{c,\alpha}\Zee}\Cee^*$ of the group. 

If $n_{c,\alpha} =1$, then there is an induced $G$-bundle over $E\times
\WP^0(\eta)$. If $n_{c,\alpha} > 1$, then it is easy to see that the corresponding
$\widehat G$-bundle does not lift to a $G$-bundle. For example, in case
$G=SL_n(\Cee)$, the vector bundle over $E\times \Pee^{n-1}$ constructed by taking
the $k^{\text{th}}$ vertex of the Dynkin diagram is given as an extension
$$0 \to \pi_1^*W_k^*\otimes \pi_2^*\scrO_{\Pee^{n-1}}(1) \to \mathbf{U} \to
\pi_1^*W_{n-k} \to 0,$$
and no twist of this bundle by a line bundle will have trivial determinant. In
a future paper, we shall discuss methods for constructing universal bundles in case
$n_{c,\alpha} > 1$ via spectral covers.

\subsection{The case of a $c$-special root}\label{s4.3}

\begin{defn}\label{regdef} Let $\zeta$ be a  semistable $G$-bundle with
$c_1(\zeta) = c$. Then $\zeta $   is {\sl regular\/} if $\dim \Aut \zeta = r_c$. By
\cite[Corollary 6.3]{FM}, every $\zeta$ is S-equivalent to a regular
semistable $G$-bundle, which is unique up to isomorphism. Moreover, by
\cite[Corollary 6.3]{FM}, if $\zeta$ is not regular, then $\dim \Aut \zeta
\geq r_c+2$.
\end{defn}

 Let $\eta_0$ be the distinguished bundle defined in
Proposition~\ref{eta0} and Lemma~\ref{ceta0}. We next show that the nonzero
elements of
$H^1(E; U(\eta_0))$ correspond to regular semistable $ G$-bundles.

\begin{proposition}\label{4.6} With $\eta_0$ as above, for every $x\in H^1(E;
U(\eta_0)) -\{0\}$, let $\xi _x$ denote the principal $P$-bundle
$\Xi|E\times \{x\}$ induced by the restriction of $\Xi$ to the
slice over $x$. Then $\xi _x\times _GP$ is a regular semistable $ G$-bundle.
\end{proposition}
\begin{proof}   By Theorem~\ref{mindim} and Theorem~\ref{cmindim}, if
$\zeta$ is unstable and $c_1(\zeta) = c$, then $\dim \Aut \zeta \geq r_c+2$. As we
have noted above, the same holds for a semistable
$G$-bundle $\zeta$ which is not regular. Thus, a
$ G$-bundle
$\zeta$ is semistable and regular if and only if $\dim \Aut \zeta \leq r_c+1$. To
prove Proposition~\ref{4.6}, we shall show that, for all
$x\neq 0$,
$\dim\Aut _{ G}(\xi_x\times _PG)
\leq r_c+1$.

Let $\xi=\xi_x$. We have the inclusion of the Lie algebra $\frak p$ in
$\frak g$. Clearly, viewing $\frak g$ as a representation of $ P$, the vector
bundle $\frak g(\xi)$ is the same as $\frak g(\xi\times_PG)=\ad _{ G}(\xi\times
_PG)$.  Moreover
$\frak p(\xi) = \ad _{ P}\xi$. Thus there is an exact sequence of vector
bundles
$$0 \to \ad _{ P}\xi \to \ad _{ G}(\xi\times_PG) \to (\frak g/\frak p)(\xi ) \to
0.$$ Now replacing $\xi$ by $\eta_0 \times _{ L} P$ gives the corresponding exact
sequence
$$0 \to \ad _{ P}(\eta_0 \times _{ L} P) \to \ad _{ G}(\eta_0 \times
_{ L} G) \to \frak u_-(\eta_0) \to 0,$$ 
since $(\frak g/\frak p)(\eta_0 \times _{ L} P) = \frak u_-(\eta_0)$. Furthermore,
by Corollary~\ref{dimad} and Lemma~\ref{d2} in the simply connected case and
Theorem~\ref{cthm} in the non-simply connected case, $H^0(E;
\frak u_-(\eta_0))$ has dimension
$r_c+1$. By semicontinuity, there is a neighborhood
$\Omega$ of the origin in $H^1(E; U(\eta_0))$ such that, if $\xi =\xi _x$
corresponds to an $x\in
\Omega$, then
$\dim H^0(E; (\frak g/\frak p)(\xi )) \leq r_c+1$ as well. As every point of
$H^1(E; U(\eta_0))$ is
$\Cee^*$-equivalent to such an $x$, we must have $\dim H^0(E; (\frak g/\frak
p)(\xi)) \leq r_c+1$ for all possible $\xi$. 

Next consider the exact sequence of Lie algebras
$$0 \to \frak u\to \frak p \to \frak l \to 0.$$
There is the associated bundle sequence
$$0 \to \frak u(\xi) \to \ad _{ P}\xi\to \ad _{ L}\eta _0 \to 0.$$
Since $\frak u(\eta_0)$ is a direct sum of semistable bundles of negative degrees,
$H^0(E;
\frak u(\eta_0)) =0$. It follows as before from   semicontinuity and
$\Cee^*$-equivariance  that
$H^0(E;
\frak u(\xi)) =0$ for all $\xi$. So  $H^0(E; \ad _{ P}\xi) \subseteq H^0(E;
\ad_{ L}\eta _0)\cong \Cee$. Thus
$$h^0(E; \ad _{ G}(\xi\times_PG))\leq h^0(E; (\frak g/\frak
p)(\xi)) + h^0(E; \ad _{ P}\xi) \leq r_c+2,$$
with equality holding if and only if the map $H^0(E; \ad _{ P}\xi)\to H^0(E;
\ad_{ L}\eta _0)$ is surjective. This can only happen if the natural
homomorphism
$\Aut _{ P}\xi \to \Aut _{ L}\eta_0$ is surjective on the connected
component of the identity, which would say that every $\lambda \in\Cee^*$ in $\Aut
_{ L}\eta_0$ lifts to an element of
$\Aut _{ P}\xi$. But then $x$ must be a fixed point for the
$\Cee^*$-action, and hence $x$ is the origin. Conversely, if $x\neq 0$, then
$H^0(E; \ad _{ P}\xi) = 0$ and $\dim
\Aut _{ G}(\xi\times_PG) \leq r_c+1$, and as we have seen above, this statement
implies Proposition~\ref{4.6}.
\end{proof}

\subsection{The Kodaira-Spencer homomorphism}\label{sKS}

For the moment, we return to the case of an arbitrary root $\alpha$.
Let $x\in H^1(E; U(\eta))$, and let $\xi$ be the corresponding
$P$-bundle. The bundle
$\Xi_0$ induces a Kodaira-Spencer homomorphism from the tangent space
$T_x$ of $H^1(E; U(\eta))$ at $x$ to $H^1(E; \ad (\xi\times _PG))$,
and we wish to find some general circumstances where this map is surjective.

\begin{theorem} Suppose that the differential of the action of  $\Aut_L\eta$
on  $H^1(E; U(\eta))$ at $1\in \Aut_L\eta$ and $x\in H^1(E; U(\eta))$ is
an injective homomorphism
$H^0(E;
\ad _L\eta) \to T_{x}$. Then the Kodaira-Spencer homomorphism $T_{x}
\to H^1(E; \ad (\xi\times _PG))$ is surjective.
\end{theorem}
\begin{proof} As in
the proof of Lemma~\ref{4.2}, the tangent space $T_x$ can be identified
with bundles $\xi_\varepsilon$ over $E\times \Spec \Cee[\varepsilon]$ which
restrict to $\xi$ over the closed fiber and reduce mod $U$ to $\pi_1^*\eta$. If
$\eta$ is given by the
$1$-cocycle $\{\ell_{ij}\}$, where the $\ell_{ij}$ take values in $L$,  and $\xi$
by the
$1$-cocycle
$\{\ell_{ij}u_{ij}\}$, where the $u_{ij}$ take values in $U$, then it is easy to
see that
$\xi_\varepsilon$ is given by a
$1$-cocycle
$\{\ell_{ij}(u_{ij}+\varepsilon v_{ij})\}$, where the $v_{ij}$ are also
$U$-valued. Moreover $w_{ij} = u_{ij}^{-1}v_{ij}$ defines an element of $H^1(E;
\frak u(\xi))$. In this way, we identify 
$T_{x}$ with
$H^1(E; \frak u(\xi))$.

There is a long exact sequence 
$$0 \to \frak u(\xi) \to \ad _P\xi \to \ad _L\eta \to 0.$$
The natural map $H^1(E; \frak u(\xi)) \to H^1(E; \ad
_P\xi)$ is the Kodaira-Spencer map for deformations of the $P$-bundle $\xi$,
and its kernel is the image of the  coboundary map
$\delta \colon H^0(E; \ad _L\eta)
\to H^1(E; \frak u(\xi))$. This kernel   also contains the image of the
differential of the action of
$\Aut_L\eta$ on  $H^1(E; U(\eta))$ at $x$, which has dimension equal to $\dim
H^0(E; \ad _L\eta)$. Thus it follows by hypothesis   that
$\delta$ is injective.

The Killing form identifies the vector bundle $(\frak g/\frak p)(\xi)$ with
the dual of $\frak u(\xi)$. In particular, $(\frak g/\frak p)(\xi)$ has a 
filtration whose successive quotients are stable bundles of positive degrees.
Hence $H^1(E; (\frak g/\frak p)(\xi)) =0$ and the natural map $H^1(E;
\ad_P\xi) \to H^1(E; \ad _G(\xi\times _PG))$ is surjective. Now consider the
commutative diagram
$$\begin{CD}
@. H^1(E; \frak u(\xi)) @>>>  H^1(E; \ad _G(\xi\times _PG))@.\\
@. @VVV @| @.\\
H^0(E; (\frak g/\frak p)(\xi)) @>>> H^1(E; \ad_P\xi) @>>>  H^1(E; \ad
_G(\xi\times _PG))@>>> 0\\
@VVV @VVV @. @.\\
H^1(E; \ad _L\eta) @= H^1(E; \ad _L\eta) @. @.\\
@. @VVV @. @.\\
@. 0 @. @.
\end{CD}$$
where the middle row  and column are exact. A diagram chase shows that, if
the map $H^0(E; (\frak g/\frak p)(\xi)) \to H^1(E; \ad _L\eta)$ is
surjective, then so is the map $$H^1(E; \frak u(\xi)) \to  H^1(E; \ad
_G(\xi\times _PG)),$$ which is the statement of the theorem. The map
$H^0(E; (\frak g/\frak p)(\xi)) \to H^1(E; \ad _L\eta)$ is given by the
composition
$$H^0(E; (\frak g/\frak p)(\xi)) \to H^1(E; \ad_P\xi)\to H^1(E; \ad _L\eta),$$
and the above sequence is Serre dual to the sequence
$$H^1(E; \frak u(\xi)) \leftarrow H^0(E; (\frak g/\frak u)(\xi)) \leftarrow 
H^0(E; \ad _L\eta).$$
By the naturality of the connecting homomorphisms associated to the
following commutative diagram of short exact sequences
$$\begin{CD}
0 @>>> \frak u(\xi)  @>>> \ad _P\xi @>>> \ad _L\xi @>>> 0\\
@. @| @VVV @VVV @.\\
0 @>>> \frak u(\xi)  @>>> \ad _G(\xi\times_PG) @>>>(\frak g/\frak u)(\xi)
@>>>0,
\end{CD}$$
the composition
$H^0(E; \ad _L\eta) \to H^1(E;
\frak u(\xi))$ is just the homomorphism $\delta$, which is injective by
hypothesis. Hence by duality 
$H^0(E; (\frak g/\frak p)(\xi)) \to H^1(E; \ad _L\eta)$ is
surjective, which completes the proof.
\end{proof}

\begin{corollary}\label{local} If $\dim \Aut _L\eta =1$, and $x$ is not the origin
of $H^1(E; U(\eta))$, then the Kodaira-Spencer  homomorphism $T_{x}
\to H^1(E; \ad (\xi\times _PG))$ is surjective.  If $\alpha$ is
$c$-special, $\eta = \eta_0$, and $x$ is not the origin, then the
Kodaira-Spencer  homomorphism induces an isomorphism from $T_{x}$ modulo
the tangent space to   $\Cee^*\cdot x$  to $H^1(E;
\ad (\xi\times _PG))$.
\end{corollary}
\begin{proof} Since the $\Cee^*$-action is via strictly positive weights, the
differential of the action is injective at every nonzero point $H^1(E;
U(\eta))$, so the hypothesis of the previous theorem is satisfied. This
proves the first statement. To see the second statement, the induced map from 
$T_{x}$ modulo the tangent space to $\Cee^*\cdot x$ to $H^1(E;
\ad (\xi\times _PG))$ is surjective. The dimension of $T_{x}$ modulo the
tangent space to $\Cee^*\cdot x$ is $r_c$. Since, for every $\xi$, $\dim H^1(E;
\ad (\xi\times _PG)) \geq r_c$, equality must hold, giving a new proof of
Proposition~\ref{4.6},  and the induced map from $T_{x}$ modulo the
tangent space to $\Cee^*\cdot x$ to $H^1(E;
\ad (\xi\times _PG))$ is an isomorphism. 
\end{proof}

This shows that the map $\Psi\colon \WP(\eta_0)\to \mathcal{M}(G,c)$ is a
local diffeomorphism over the smooth points of $\mathcal{M}(G,c)$. 
Moreover, giving $\mathcal{M}(G,c)$ and $\WP(\eta_0)$ their natural
orbifold structures, it follows from Corollary~\ref{local} 
that $\Psi$ is an orbifold covering.
Let $d_k(\WP(\eta_0))$ and $d_k(\mathcal{M}(G,c))$ be the dimensions of
the subspaces of $\WP(\eta_0)$ and $\mathcal{M}(G,c)$ where $k$ divides the
order of the orbifold isotropy. The orbifold covering property implies that 
\begin{equation*}
d_k(\WP(\eta_0))\le d_k(\mathcal{M}(G,c))\ \ {\rm for\ all\ \ } k\ge 1.
\end{equation*}
The results in \cite{BFM} and \cite{FM} imply that $d_k(\mathcal{M}(G,c))= d(k)$ in
the notation of Proposition~\ref{1weights}  in the simply connected case.

In Section 5 we shall prove that $\Psi$ is an
isomorphism. In proving this result we do not appeal to Corollary~\ref{local}, but
rather use the fact that both
$d_k(\WP(\eta_0))$ and $d_k(\mathcal{M}(G,c))$ satisfy circular symmetry.
One can in fact turn this argument around. Using the above
inequality and the fact that the sums of the weights for
$\WP(\eta_0)$ and $\mathcal{M}(G,c)$ add up to $g$, one can prove directly
that  $d_k(\WP(\eta_0))=d_k(\mathcal{M}(G,c))$ for all $k\ge 1$, and hence
apply the results of \cite{BFM} to show that the $d_k(\WP(\eta_0))$ satisfy
circular symmetry. In the simply connected case, this gives a  classification-free
proof of  circular symmetry for the numbers $d_k(\alpha)$, where $\alpha$ is a
special root.

\subsection{The singular locus of the weighted projective space}

The weighted projective space $\WP(\eta)$ is naturally an orbifold. Its singular
locus (as an orbifold) corresponds to the set of points in the affine space
$H^1(E;U(\eta))$ whose isotropy group is larger than that of the generic point,
i.e.\ is larger than $\Zee/n_{c,\alpha}\Zee$. This will be the singular locus of
$\WP(\eta)$ as a variety provided that its codimension is at least two. If we
choose a linear structure and a diagonal basis of
$H^1(E;U(\eta))$, then the orbifold singular locus of $\WP(\eta)$ is a union of
weighted projective subspaces corresponding to setting all of the coordinates
equal to zero except those for which the weights are divisible by
$kn_{c,\alpha}$, where $k>1$. Our goal is to show that, when $G$ is simply
connected and $\alpha$ is special, each such subspace can be naturally identified
the moduli space for a non-simply connected subgroup of
$G$.

Recall from the proof of Proposition~\ref{highest} that the set of all $\beta
\in R$ such that $k|\beta(\varpi_\alpha\spcheck)$ is a root system $R(\alpha,
k)$. Moreover, $\Delta(\alpha,k)=(\Delta -\{\alpha\})\cup \{-\lambda_k(\alpha)\}$
is a set of simple roots for $R(\alpha, k)$. There is a semisimple subalgebra
$$\frak g(\alpha,k) = \frak h \oplus \bigoplus _{\beta \in R(\alpha, k)}\frak
g^\beta
\subseteq \frak g.$$ 
Let $G(\alpha,k)$ be the corresponding closed connected subgroup of $G$. Of course,
$\frak g(\alpha,k)$ will not be simple in general.
 Let
$\Delta(\alpha,k) = \coprod_{i\geq 1}\Delta(\alpha,k)_i$, where each subset
$\Delta(\alpha,k)_i$ corresponds to a connected component of the Dynkin diagram
of $R(\alpha,k)$, and where $-\lambda_k(\alpha) \in \Delta(\alpha,k)_1$.  Since
$G(\alpha, k)$ is semisimple, we can write
$$G(\alpha, k) =\left(\prod _iG_i\right)\Big/F,$$ where each $G_i$ is simple and
simply connected and corresponds to the subset $\Delta(\alpha,k)_i$ of
$\Delta(\alpha,k)$,  and where
$F$ is finite. If
$G$ is simply connected, then $F$ is cyclic of order $k$, generated by an element
$c_k$. 

Viewing
$-\lambda_k(\alpha)$ as an element of $\Delta(\alpha,k)$, i.e.\ a simple root for
$G(\alpha,k)$, $-\lambda_k(\alpha)$ defines a maximal parabolic subgroup
$P(\alpha,k)$ of
$G(\alpha,k)$ contained in the maximal parabolic subgroup $P^\alpha =P$ of $G$
determined by $\alpha$. The parabolic subgroup $P(\alpha,k)$ is of the form
$(P_1\times \prod_{i\geq 2}G_i)/F$, where $P_1$ is the maximal parabolic subgroup
in $G_1$ corresponding to $-\lambda_k(\alpha)$. Clearly, the  Levi factor
$(\alpha,k)$ of
$P(\alpha,k)$ is just
$L^\alpha = L= (L_1\times \prod_{i\geq 2}G_i)/F$, where $L_1$ is the Levi factor
of $P_1$. The  unipotent radical $U(\alpha, k)$ of $P(\alpha,k)$ has Lie algebra
$\bigoplus _{k|j}\frak u^j$. Let $\eta$ be a semistable $L$-bundle of negative
degree. The degree of $\eta$ is of course independent of whether we view $L$ as
the Levi factor of $P$ or of $P(\alpha,k)$, and we can define the cohomology set
$H^1(E; U(\alpha, k)(\eta))$. The inclusion of $U(\alpha, k)$ in $U$ defines a
$\Cee^*$-equivariant function $H^1(E;U(\alpha, k)(\eta))\to H^1(E;U(\eta))$ on
the level of cohomology sets, as well as a morphism between the corresponding
functors. Since the two  associated functors are both represented by affine
spaces, the function $H^1(E;U(\alpha, k)(\eta))\to H^1(E;U(\eta))$ is a
$\Cee^*$-equivariant morphism of affine
spaces. The geometric meaning of this morphism is as follows: let $F_i$ be the
projection of $F$ to the factor $G_i$, so that there is a homomorphism from
$(L_1\times \prod_{i\geq 2}G_i)/F$ to $L_1/F_1\times \prod_{i\geq 2}(G_i/F_i)$ for
$i\geq 2$. The bundle
$\eta$ thus induces an $L_1/F_1\times \prod_{i\geq 2}(G_i/F_i)$-bundle 
$\prod_i\eta_i$, where $\eta_1$ is an $L_1/F_1$-bundle
and the $\eta_i$ are
$G_i/F_i$-bundles  for $i\geq 2$, Moreover $\eta$ defines a canonical lifting
of the $(L_1/F_1)\times \prod _{i\geq 2}(G_i/F_i)$-bundle $\prod_i\eta_i$
to an $(L_1\times \prod_{i\geq 2}G_i)/F$-bundle. A class $x$ in $H^1(E;U(\alpha,
k)(\eta))$ defines a lifting of $\eta_1$ to a $P_1/F_1$-bundle $\xi_1$. The lift
$\eta$ then defines a lift of the $(P_1/F_1)\times \prod _{i\geq
2}(G_i/F_i)$-bundle defined by $\xi_1$ and the $\eta_i$, $i\geq 2$,  to a
$(P_1\times
\prod_{i\geq 2}G_i)/F$-bundle $\xi'$. Since $P(\alpha,k) =(P_1\times
\prod_{i\geq 2}G_i)/F$ is a subgroup of $P$, we can form the associated bundle
$\xi = \xi'\times _{P(\alpha,k)}P$, and this bundle is clearly the lift of $\eta$
corresponding to the image of $x$ in $H^1(E;U(\eta))$.

\begin{proposition}\label{singembed} Let $\eta$ be a semistable $L$-bundle of
negative degree
$-d$. There are compatible linear structures on
$H^1(E;U(\alpha, k)(\eta))$ and on $H^1(E;U(\eta))$ so that the morphism
$H^1(E;U(\alpha, k)(\eta))\to H^1(E;U(\eta))$ is a $\Cee^*$-equivariant embedding
of
$H^1(E;U(\alpha, k)(\eta))$ onto the linear subspace of $H^1(E;U(\eta))$ defined
by the span of all of the eigenvectors of the $\Cee^*$-action on $H^1(E;U(\eta))$
whose weights are divisible by $kn_{c,\alpha}d$.
\end{proposition}
\begin{proof}  It is   an elementary exercise to check that, if $\Cee^*$ acts
linearly and with positive weights on two affine spaces $\mathbb{A}_1$ and
$\mathbb{A}_2$ and  if $f\colon
\mathbb{A}_1\to \mathbb{A}_2$ is a $\Cee^*$-equivariant morphism whose
differential at the origin is injective, then there exist coordinates on
$\mathbb{A}_2$ for which $\Cee^*$ acts linearly and such that $f$ is a linear
embedding, and the $\Cee^*$-weights for the image of $f$ can be determined from
the differential of $f$ at the origin. The final statement then follows from the
differential computation. The differential of the morphism  $H^1(E;U(\alpha,
k)(\eta))\to H^1(E;U(\eta))$ is given by the inclusion 
$$\bigoplus _{k|j}H^1(E; \frak u^j(\eta)) \to \bigoplus _{j>0}H^1(E; \frak
u^j(\eta)).$$
Hence the image of the differential of $f$ at the origin is the span of the 
eigenvectors of the $\Cee^*$-action on $H^1(E;\frak u(\eta))$
whose weights are divisible by $kn_{c,\alpha}d$. Thus, the same ids true for $f$.
\end{proof} 

We turn now to the case of a special root $\alpha$. For simplicity, we assume
that  $G$ is simply connected, so that the finite group $F=
\pi_1(G(\alpha,k))$ is generated by an element $c_k$ of order $k$. 

\begin{proposition} Suppose that $G$ is simply connected and that $\alpha$ is
special. Let $\eta_0$ be the $L$-bundle of Proposition~\ref{eta0}.  Then:
\begin{enumerate}
\item[\rm (i)] For $i>1$, $\Delta(\alpha,k)_i$ is of type $A$ and $c_k$ projects
to a generator of the corresponding fundamental group.
\item[\rm (ii)] The  $(L_1/F_1)$-bundle $\eta_1$ induced by $\eta_0$
has the property that $\eta_1\times _{L_1/F_1}(G_1/F_1)$ is a minimally unstable
$(G_1/F_1)$-bundle, and the root
$-\lambda_k(\alpha)$ is a $c_k$-special simple root in 
$\Delta(\alpha,k)_1$. 
\end{enumerate}
\end{proposition}
\begin{proof} Part (i) follows easily from the explicit description of the
special root. To see Part (ii), it follows from Proposition~\ref{singembed} that,
if $x$ is a nonzero element of $H^1(E;U(\alpha, k)(\eta_0))$, then the image of
$x$ in $H^1(E; U(\eta_0))$ is also nonzero. In particular, if $\xi$ is the
corresponding $P$-bundle, then $\xi\times _PG$ is semistable. It is easy to check
that, in this case, the $P_1/F_1$-bundle corresponding to $x$ is again
semistable. Thus, $\eta_1\times _{P_1/F_1}(G_1/F_1)$ is a minimally unstable
$(G_1/F_1)$-bundle, and so $-\lambda_k(\alpha)$ is a $c_k$-special simple root for
$G_1/F_1$.
\end{proof}

Of course, we could also check the above proposition by a case-by-case analysis.
This  also shows that the projection $F\to G_1$ is always an embedding of $F$
into the center of $G_1$. For $i\geq 2$, the factor  $G_i$ is of type $A$ and
$F_i$ is the full center, and hence the bundle $\eta_i$ is always rigid. On the
other hand,  the image $F_1$ need not be the full center of
$G_1$. For example, if  $G$ is of type $E_8$ and $k=4$, then $G_1$
is of type $A_7$, and thus its center is isomorphic to $\Zee/8\Zee$, whereas
$F_1$ has order $4$. 

There is thus an induced morphism of weighted projective spaces. In terms of
moduli spaces, if we grant Looijenga's theorem (Theorem~\ref{LooThm}) in both the
simply connected and the non-simply connected cases, this morphism identifies the
sub-projective space of $\WP(\eta_0)$ corresponding to setting all the coordinates
in the weight spaces with weights not divisible by $kn_\alpha$ equal to zero with
the weighted projective space which is the  moduli space of unliftable
semistable $(G_1/F_1)$-bundles. Of course, on the level of $G$-bundles, this
shows that up to S-equivalence a bundle corresponding to a point of the moduli
space lying in this sub-projective space has a reduction of structure to a
semistable unliftable
$(\prod_iG_i)/F$-bundle, and conversely such a bundle defines a semistable
$G$-bundle whose moduli point lies  in this sub-projective space.   In terms of
the
$\Cee^*$-weights,  the number of such weights divisible by
$kn_\alpha$ can be related to the weights appearing for the appropriate
non-simply connected form of a subgroup. For example, if $G$ is of type $E_8$,
the root system $R(\alpha,2)$ is of type $E_7\times A_1$ and the $\Cee^*$-weights
divisible by $2$ for $G$, in other words the $g_\alpha$ such that $2|g_\alpha$,
are the weights occurring in the weighted projective space for the adjoint form of
$E_7$. These are exactly the weights appearing in the quotient diagram for
$\widetilde E_7$ modulo the action of the nontrivial element of the center,
namely twice the weights for the group of type
$F_4$:
$2,2,4,4,6$. In \cite[\S9]{BFM}, these quotient root systems appear in a different
context, unrelated to the special roots, as certain root systems $\Phi(\frak
t(k))$ constructed on certain subtori of $H$. It would be nice to understand this
somewhat mysterious connection more directly.

\section{A new proof of Looijenga's theorem}

\subsection{Statement of the theorem}

Fix a $c$-special root $\alpha$. We denote the corresponding
parabolic subgroup simply by $ P$, and similarly for $ L$. We have
defined the bundle
$\eta_0$ in Proposition~\ref{eta0} and Lemma~\ref{ceta0}. As in Section 1, let
$\mathcal{M}(G,c)$ denote the coarse moduli space of semistable $ G$-bundles
$\xi$ with $c_1(\xi) =c$, modulo S-equivalence.   We have seen that there is a
universal family of regular semistable
$G$-bundles $\Xi_0$ over $E\times (H^1(E; U(\eta_0))-\{0\})$. Thus there is an
induced morphism 
$$\widetilde\Psi\colon  H^1(E; U(\eta_0))-\{0\} \to \mathcal{M}(G,c).$$   The
morphism
$\widetilde\Psi$ is constant on
$\Cee^*$-orbits. Let $\WP(\eta_0)$ be the weighted projective space which is the
$\Cee^*$-quotient of $H^1(E; U(\eta_0))-\{0\}$ and let
$$\Psi\colon  \WP(\eta_0) \to \mathcal{M}(G,c)$$ be the morphism induced by
$\widetilde\Psi$. We can then state our version of Looijenga's theorem as follows:

\begin{theorem}\label{LooThm} Let $E$ be an elliptic curve, let $G$ be a simple
group and let $c\in \pi_1(G)$ be a generator. The morphism $\Psi\colon 
\WP(\eta_0) \to
\mathcal{M}(G,c)$ defined above is an isomorphism.
\end{theorem}

\begin{corollary} The moduli space $\mathcal{M}(G,c)$ is isomorphic to a weighted
projective space with weights $g_{\ov \beta}/n_0, \ov \beta \in \widetilde
\Delta/w_c$ as defined in Definition~\ref{defrc}.
\end{corollary}

In particular, if $G$ is simply connected, then we obtain a new proof of
Looijenga's theorem \cite{Loo}.

\subsection{The classical cases}

Let us first sketch the proof of the theorem above for the case of the classical
groups, bearing in mind the description of the bundle $\eta_0$ given in
Proposition~\ref{classical}. In case
$G=SL_n(\Cee)$, the theorem asserts that every regular semistable vector bundle
$V$ of rank $n$ and trivial determinant is S-equivalent to a  unique extension
$$0 \to W_k^*\to V \to W_{n-k} \to 0.$$ In this form, the theorem is proved in
\cite{FMW2}, Theorem 3.2(iv).

Let us next consider the case of the symplectic group. We shall show that the
morphism $\Psi$ has degree one in this case (as we shall see below, this implies
that $\Psi$ is an isomorphism). Let $V$ be a generic regular semistable symplectic
vector bundle, in other words a regular semistable vector bundle of rank $2n$ with
a nondegenerate symplectic form $A$. Here generic will mean that $V$ is a direct
sum
$$\bigoplus _{i=1}^n (\lambda _i\oplus \lambda_i^{-1}),$$ where the $\lambda_i$
are line bundles of degree zero, not of order $2$, and such that, for $i\neq j$,
$\lambda_i \neq \lambda_j^{\pm 1}$. In this case, the symplectic form on $V$ is an
orthogonal sum of symplectic forms $A_i$ on $\lambda _i\oplus \lambda_i^{-1}$. The
space of such forms which are nondegenerate corresponds to the choice of an
isomorphism from $\lambda_i^{-1}$ to itself, in other words to a nonzero multiple
of $A_i$, and the group of symplectic automorphisms of
$A_i$ is also isomorphic to $\Cee^*$. For each
$i$, the space of surjections $\varphi_i\colon  W_n^* \to \lambda _i^{\pm 1}$ is a
$\Cee^*$. Thus, the space of morphisms  $ W_n^* \to \lambda _i \oplus \lambda
_i^{-1}$ is a $\Cee^*\times \Cee^*$. It is clear that the pullback of $A_i$ to
$W_n$ under such a morphism is a nonzero symplectic form, which we denote by
$B_i$. Moreover, by varying the morphism from $ W_n^* \to \lambda _i \oplus \lambda
_i^{-1}$, we exactly get all symplectic forms on $W_n^*$ of the form $t_iB_i$,
where
$t_i\in \Cee^*$.

Now suppose we have found, for every $i$, a morphism $f_i\colon  W_n^* \to \lambda
_i 
\oplus \lambda _i^{-1}$ such that, if $t_iB_i$ is the pulled back morphism, then
$t_i\neq 0$ for every $i$ and such that $\sum _{i=1}^nt_iB_i = 0$. It follows that
the morphism $W_n ^* \to \lambda _i^{\pm 1}$ is nonzero for every $i$. By
Proposition 3.6 of \cite{FMW2}, the morphism $W_n^* \to V$ embeds $W_n^*$ as an
isotropic subbundle of $V$, and the quotient $V/W_n^*$ is isomorphic to $W_n$.
This implies that we have reduced the structure group of $V$ to a maximal
parabolic subgroup $P$ of $Sp(2n)$ corresponding to the special root, and the
corresponding $L$-bundle is $\eta_0$. Hence $(V, A)$ is in the image of
$\Psi$, and the statement that
$\Psi$ has degree $1$ is the statement that the $t_i$ above are uniquely
determined up to multiplying by a fixed nonzero scalar.

Thus, we must find $t_i$ such that $\sum _{i=1}^nt_iB_i = 0$ and show that the
$t_i$ are unique up to a scalar. To see this, note that the space of alternating
forms on
$W_n^*$ is
$H^0(E;
\bigwedge ^2W_n)$. One easily computes that $\deg \bigwedge ^2W_n = n-1$. Since
$W_n$ is stable,
$\bigwedge ^2W_n$ is a semistable vector bundle of positive degree. Thus, $h^1(E;
\bigwedge ^2W_n) = 0$, and so by Riemann-Roch, 
$h^0(E; \bigwedge ^2W_n) = n-1$. It follows that every collection of $n$ forms
$B_i$ is linearly dependent, and so some linear combination of the $B_i$ is zero.
To prove uniqueness, and also to prove that the quotient is a $W_n$, it will
suffice to show that no smaller linear combination is zero. Suppose, say, that
$\sum _{i=1}^kt_iB_i = 0$, with $k< n$ and $t_i\neq 0$ for all $i\leq k$. Consider
the induced morphism from $W_n^*$ to $\bigoplus _{i=1}^k (\lambda _i\oplus
\lambda_i^{-1})$. If
$k > n/2$, then by Proposition 3.6 of \cite{FMW2} there would exist an embedding of
$W_n^*$ as an isotropic subbundle of a vector bundle of rank $< 2n$, which is
impossible since the symplectic form on $\bigoplus _{i=1}^k (\lambda _i\oplus
\lambda_i^{-1})$ is nondegenerate. If $k\leq n/2$, then in fact the argument of
Proposition 3.6 of
\cite{FMW2} shows that the image of $W_n^*$ would be a subbundle of rank equal to
$2k$, and hence that the symplectic form on $\bigoplus _{i=1}^k (\lambda _i\oplus
\lambda_i^{-1})$ is identically zero, which is again a contradiction. Thus no
smaller linear combination of the $B_i$ is zero, proving that $\Psi$ has degree
one. A similar argument shows that, even for a regular symplectic bundle
$(V, A)$ which is not generic in the above sense, there still exists an embedding
of $W_n^*$ in $V$ as a maximal isotropic subbundle, and this embedding is unique
up to the automorphism group of $(V, A)$. Thus, we can also show directly that
$\Psi$ is a bijection instead of merely having degree one in this case.

The case of $Spin(n)$ is similar. By a generic $Spin(2n)$-bundle, we mean one
whose associated $SO(2n)$-bundle $V$ is of the form
$\bigoplus _{i=1}^n (\lambda _i\oplus \lambda_i^{-1})$, where the $\lambda_i$
are line bundles of degree zero, not of order $2$, and such that, for $i\neq j$,
$\lambda_i \neq \lambda_j^{\pm 1}$, and the quadratic form is the  
orthogonal sum of quadratic forms $A_i$ on $\lambda _i\oplus \lambda_i^{-1}$. We
consider maps $W_{n-2}^* \to V$ whose image is isotropic.   Arguments as in the
symplectic case show that there is an embedding $\iota$ of $W_{n-2}^*$ as an
isotropic subbundle of $V$, and such that the  projection to each
summand $\lambda_i^{\pm 1}$ is nonzero, and such an embedding is unique up to
orthogonal isomorphisms of
$V$. Dually, we have a map $V^*\cong V \to W_{n-2}$. Thus there is a complex
$W_{n-2}^*\xrightarrow{\iota} V \xrightarrow{\iota^*} W_{n-2}$.  The symmetric form
identifies
$V/(W_{n-2}^*)^\perp$ with $W_{n-2}$. It is easy to check that the bundle
$\Ker(\iota^*)/\operatorname{Im}(\iota)$ is a semistable $SO(4)$-bundle which
does not lift to $Spin(4)$, and hence is of the form $W_2\otimes W_2^* \cong
Q_4$. The filtration on $V$ then reduces the structure group to the appropriate
maximal parabolic subgroup as before. A very similar argument handles the case of
$SO(2n+1)$, using instead an isotropic subbundle isomorphic to $W_{n-1}^*$, and
showing that $(W_{n-1}^*)^\perp/W_{n-1}^*$ is isomorphic to $\ad W_2$.
Similar statements can also handle the case of unliftable bundles.

\subsection{Proof  of the main theorem: determinant line bundles}

We turn to the proof of Theorem~\ref{LooThm}. 

First, it is an elementary result that a morphism from a weighted projective space
to a quasiprojective variety is either constant or finite. Indeed, since every
weighted projective space has a finite ramified cover which is an ordinary
projective space, this follows from the analogous and well-known result for
$\Pee^n$. Thus, since $\WP(\eta_0)$ and $\mathcal{M}(G,c)$ are normal, if we can
show that
$\Psi$ has degree one, then it follows from Zariski's main theorem that $\Psi$ is
an isomorphism.

To calculate the degree of $\Psi$ we shall compare determinant line bundles on the
two sides. The idea will be to find line bundles $\mathbb{D}_1$, 
$\mathbb{D}_2$ on 
$\WP(\eta_0)$ and $\mathcal{M}(G,c)$ respectively, such that
$[\mathbb{D}_1] =\Psi ^*[\mathbb{D}_2]$ in $H^2(\WP(\eta_0); \Q)$ and such that
$\int _{\WP(\eta_0)}c_1(\mathbb{D}_1)^{r_c} = \int
_{\mathcal{M}(G,c)}c_1(\mathbb{D}_2)^{r_c}
\neq 0$. Since on the other hand we have $c_1(\mathbb{D}_1)^{r_c} = \deg \Psi \cdot
c_1(\mathbb{D}_2)^{r_c}$, it will then follow that $\deg \Psi = 1$. 

To define the line bundles $\mathbb{D}_i$ (which in fact will be {\it a priori\/}
only
$\Q$-Cartier divisors), we recall the definition of the determinant line bundle on
the moduli functor. First recall the definition of the moduli functor
$\mathbf{M}$ itself: for a scheme $S$ over $\Cee$, $\mathbf{M}(S)$ is the set of
isomorphism classes of principal $G$-bundles $\Xi$ over $E\times S$ such that
$\Xi|E\times
\{x\}$ is semistable for every $x\in S$. The moduli functor is coarsely
represented by $\mathcal{M}(G,c)$. Given an element of $\mathbf{M}(S)$,
corresponding to  a principal $G$-bundle $\Xi$ over
$E\times S$, we have the associated vector bundle
$\ad \Xi$ over
$E\times S$, and thus we can form the determinant line bundle $\det R\pi _2{}_*\ad
\Xi$ over $S$ (see for example Chapter 5, Section 3 of \cite{FMbook}). Since
$\pi_2$ has relative dimension one, $H^i(\ad \Xi|E\times \{s\}) = 0$  and
likewise $R^i\pi _2{}_*\ad \Xi = 0$ for $i> 1$. The fiber of  $\det R\pi
_2{}_*\ad
\Xi$ over $s\in S$ is then the complex line 
$$\bigwedge ^{\text{\rm{top}}}H^0(\ad \Xi|E\times \{s\}) \otimes \left(\bigwedge
^{\text{\rm{top}}}H^1(\ad \Xi|E\times \{s\})\right)^{-1}.$$ Here are some of the
basic properties of this line bundle:
\begin{enumerate}
\item If $R^0\pi _2{}_*\ad \Xi$ and $R^1\pi _2{}_*\ad \Xi$ are locally free, for
example if $\Xi|E\times \{x\}$ is regular for every $x\in S$, then 
$$\det R\pi _2{}_*\ad \Xi = \left(\bigwedge ^{\text{top}}R^0\pi _2{}_*\ad
\Xi\right)\otimes \left(\bigwedge ^{\text{top}}R^1\pi _2{}_*\ad \Xi\right)^{-1}.$$
\item Suppose that $S$ is smooth, that $\lambda$ is a line bundle of degree zero
on $E$,  and that
$$D_\lambda = \{\, x\in S: h^0(E; (\ad \Xi|E\times \{x\})\otimes \lambda) \neq 0$$
is a hypersurface in $S$. Let $Z_i$ be the irreducible components of $D_\lambda$,
and let $n_i$ be the length of the torsion sheaf $R^1\pi _2{}_*(\ad \Xi\otimes \pi
_1^*\lambda)$ at a generic point of $Z_i$. Then (see e.g. \cite{FMbook}, Chapter 5,
Corollary 1.2 and Proposition 3.9) there is a canonical section
$\operatorname{div}$ of 
$\big(\det R\pi _2{}_*(\ad \Xi\otimes
\pi_1^*\lambda)\big)^{-1}$ whose divisor of zeroes is $\sum _in_i Z_i$, and hence
$$c_1(\det R\pi _2{}_*\ad \Xi) = c_1(\det R\pi _2{}_*(\ad \Xi\otimes
\pi_1^*\lambda)) = -\sum _in_i[Z_i].$$
\end{enumerate}

By general results, there is an associated
$\Q$-divisor on $\mathcal{M}(G,c)$ (i.e.\ an element in $\Pic( \mathcal{M}(G,c))
\otimes \Q$, which we shall denote by
$\mathbb{D}_2$. As we shall see, $\mathbb{D}_2$ is in
fact Cartier, in other words, a line bundle.  We will find an analogous divisor
$\mathbb{D}_1$ over $\WP(\eta_0)$ and show   that (1) $\int
_{\WP(\eta_0)}c_1(\mathbb{D}_1)^{r_c} = \int
_{\mathcal{M}(G,c)}c_1(\mathbb{D}_2)^{r_c}$ and (2)
$[\mathbb{D}_1] =\Psi ^*[\mathbb{D}_2]$. This will prove Theorem~\ref{LooThm}.

\subsection{The divisor on the weighted projective space}

Consider the universal $G$-bundle $\Xi_0$ over $E\times H^1(E; U(\eta_0))$.
The action of $\Cee^*$ on $\Xi_0$ gives a linearization of the action on $\Cee^*$
on the associated vector bundle $\ad \Xi_0$ and thus on the line bundle $\det R\pi
_2{}_*\ad \Xi_0$.  Since $H^1(E; U(\eta_0))$ is an affine space, $\det R\pi
_2{}_*\ad \Xi_0$ is the trivial line bundle. Every linearization of the
$\Cee^*$-action on the trivial bundle  is given by a character
$\chi\colon  \Cee^*\to \Cee^*$, which is multiplication by $a\in \Zee$. There is
the corresponding coherent sheaf $\scrO_{\WP(\eta_0) }(a)$, which is given by
viewing $\WP(\eta_0)$ as $\operatorname{Proj}\Cee[z_0, \dots, z_{r_c}]$ with the
appropriate grading.  By a general result on weighted projective spaces
\cite{Mori},   the coherent sheaf
$\scrO_{\WP(\eta_0) }(a)$ is a line bundle if and only if
the $\Cee^*$-weights all divide $a$.  Note that, on the open set
$\WP^0(\eta_0)$ of free $\Cee^*$-orbits, the vector bundle $\ad \Xi_0$ is defined,
and in fact
$R^0\pi _2{}_*\ad \Xi_0$ and $R^1\pi _2{}_*\ad \Xi_0$ are both locally free of
rank $r_c$. Thus  there is a well-defined line bundle  $\left(\bigwedge
^{r_c}R^0\pi _2{}_*\ad \Xi_0\right)\otimes \left(\bigwedge ^{r_c}R^1\pi _2{}_*\ad
\Xi_0\right)^{-1}$, and this line bundle clearly agrees with the restriction of
$\scrO_{\WP(\eta_0)}(a)$ to $\WP^0(\eta_0)$. We next identify the integer $a$:

\begin{lemma}\label{5.4} The natural $\Cee^*$-linearization on  $\det R\pi
_2{}_*\ad \Xi_0$ corresponds to the line bundle
$\scrO_{\WP(\eta_0) }(-2gn_{c,\alpha}/n_0)$.
\end{lemma}
\begin{proof} We must show that the $\Cee^*$-linearization on $\det R\pi_2{}_*\ad
\Xi_0$ is given by the character which is raising to the power
$-2gn_{c,\alpha}/n_0$.  To compute the
$\Cee^*$-linearization, it suffices to compute the action of
$\Cee^*$ on the fiber of
$\det R\pi_2{}_*\ad \Xi_0$ over the origin, which is a fixed point for the
$\Cee^*$-action on $H^1(E; U(\eta_0))$. The fiber over $0$ is canonically
$$\bigwedge ^{\text{top}}H^0(E; \ad _G(\eta_0\times _LG)) \otimes \left(\bigwedge
^{\text{top}}H^1(E; \ad _G(\eta_0\times _LG))\right)^{-1}.$$ Now, by
Lemma~\ref{adG} and Lemma~\ref{uplus},  
$$H^0(E; \ad _{ G}(\eta_0\times _{ L} G)) = H^0(E; \frak u_-(\eta_0))
\oplus H^0(E;
\ad _{ L}\eta_0).$$ 
Since $\Cee^*$ is contained in the center of $L$, the action of $\Cee^*$ on $H^0(E;
\ad _{ L}\eta_0)$ is trivial. By Proposition~\ref{1weights} and
Proposition~\ref{cweights},
$\Cee^*$ acts on $H^0(E; \frak u_-(\eta_0))$ with weights $-n_{c,\alpha}g_{\ov
\beta}/n_0$. Thus the action on $\bigwedge ^{\text{top}}H^0(E; \ad
_{G}(\eta_0\times _{ L} G))$ is via
$-\sum _{\ov \beta}n_{c,\alpha}g_{\ov
\beta}/n_0 = -gn_{c,\alpha}/n_0$. A similar argument (or duality)
handles the case of the
$\Cee^*$-action on $\bigwedge ^{\text{top}}H^1(E; \ad _{ G}(\eta_0\times _{L}
G))$. Putting these together, we get the power $-2gn_{c,\alpha}/n_0$.
\end{proof}

\begin{remark} In case $\WP(\eta_0) $ is the weighted projective space arising
from the moduli space of
$G$-bundles, where $G$ is simply connected, one can use the Kodaira-Spencer map to
check that
$\mathbb{D}_1 = K_{\WP(\eta_0) }^{\otimes 2}$. Now, by a standard fact about
weighted projective spaces,
$K_{\WP(\eta_0) } =  \scrO_{\WP(\eta_0) }(-g)$,   and thus
$K_{\WP(\eta_0) }^{\otimes 2} = \scrO_{\WP(\eta_0) }(-2g)$. 
\end{remark}

Define $\mathbb{D}_1=  \scrO_{\WP(\eta_0)}(-2gn_{c,\alpha}/n_0)$. We now
compute the top self-intersection of $c_1(\mathbb{D}_1)$:

\begin{lemma}\label{5.3} Let $\WP^r$ be the  weighted projective space
 which is the quotient of
$\Cee^{r+1}-\{0\}$ by the action of $\Cee^*$ acting with positive weights $w_0,
\dots, w_r$, and let $a$ be an integer such that $w_i|a$ for every $i$. Then 
$$\int_{\WP^r}c_1(\scrO_{\WP^r }(a))^r = a^rd/(w_0\cdots w_r),$$ where $d =
\gcd\{w_0,
\dots, w_r\}$.
\end{lemma}
\begin{proof}  The
morphism $\Cee^{r+1} \to \Cee^{r+1}$ defined by
$$(z_0, \dots, z_r) \mapsto (z_0^{w_0}, \dots, z_r^{w_r})$$ is
$\Cee^*$-equivariant, where $\Cee^*$ acts with all weights equal to $1$ on the
domain and with weights $w_i$ on the range. Thus there is an induced cover
$f\colon \Pee^r
\to \WP^r $, and it is easy to check that the degree of this cover is
$w_0\cdots w_r/d$. There is always a natural inclusion $f^*\scrO_{\WP^r }(a) \to
\scrO_{\Pee ^r}(a)$, and one checks that this inclusion is an isomorphism
if $w_i|a$ for every $a$. In this case,
$$\int_{\WP^r}c_1(\scrO_{\WP^r }(a))^r =
\int_{\Pee^r}c_1(\scrO_{\Pee^r}(a))^r/\deg f = a^rd/(w_0\cdots w_r).$$ This proves
the formula of Lemma~\ref{5.3}.
\end{proof}

\begin{corollary}\label{5.5} $c_1(\mathbb{D}_1)^{r_c} =
(-2g)^{r_c}n_0/\prod_{\ov\beta} g_{\ov\beta}$.
\end{corollary}
\begin{proof} This is immediate from Lemma~\ref{5.4} and Lemma~\ref{5.3}, with
$w_i=n_{c,\alpha}  g_{\ov\beta}/n_0$,   $a = -2gn_{c,\alpha} /n_0$,
and $d = n_{c,\alpha} $.
\end{proof}

\subsection{The divisor on the moduli space: the simply connected case}

We turn now to the calculation of $c_1(\mathbb{D}_2)^r$. In order to make the
argument easier to follow, we begin by working out the simply connected case.
Recall that we have the finite morphism
$E\otimes _\Zee\Lambda
\to
\mathcal{M}(G,c)$. Let $\mathbb{D}_3$ be the pullback of $\mathbb{D}_2$ to a
divisor on
$E\otimes _\Zee\Lambda$. Clearly, $c_1(\mathbb{D}_2)^r
=c_1(\mathbb{D}_3)^r/\#(W)$. Thus, we shall begin by computing
$c_1(\mathbb{D}_3)^r$. For a point $\rho\in E\otimes _\Zee\Lambda$, there is an
associated flat $G$-bundle $\xi _0$, and as we have seen in Lemma 3.1 of
\cite{FM}, 
$$\ad \xi_0 \cong 
\scrO_E^r \oplus \bigoplus _{\alpha\in R}\lambda_{\alpha(\rho)},$$ where
$\lambda_{\alpha(\rho)}$ is the line bundle of degree zero corresponding to the
point $\alpha (\rho) \in E\cong \Pic^0E$, or equivalently is the line bundle
associated to the flat $U(1)$-bundle whose holonomy is given by $\alpha(\rho)$. In
particular, we see that for $\lambda =\scrO_E$, the set $D_\lambda$ defined in
the discussion on determinant line bundles is all of
$E\otimes _\Zee\Lambda$, whereas for a nontrivial line bundle $\lambda$ of degree
zero,
$D_\lambda$ is a hypersurface in $E\otimes _\Zee\Lambda$. In fact, $D_\lambda$ is
a union of distinct hypersurfaces $D_{\lambda, \alpha}$, where if $\lambda$
corresponds to the point $e\in E$, then
$$D_{\lambda, \alpha} =\{\,\rho\in E\otimes _\Zee\Lambda: \alpha (\rho) =
-e\,\}.$$ Each $D_{\lambda,
\alpha}$ is a union of translates of abelian subvarieties of
$E\otimes _\Zee\Lambda$. In particular, the hypersurface $D_{\lambda, \alpha}$ is
smooth. The next lemma says that every component of
$D_{\lambda,
\alpha}$ counts with multiplicity one in the expression for $-c_1(\mathbb{D}_3)$.

\begin{lemma}\label{5.7} For $\lambda$ a nontrivial line bundle of degree zero, 
$$c_1(\mathbb{D}_3) = -\sum _{\alpha \in R}[D_{\lambda, \alpha}].$$
\end{lemma}  
\begin{proof} There is a universal $G$-bundle $\Xi_1$ over $E\times (E\otimes
_\Zee\Lambda)$, which in fact arises from a universal $H$-bundle, which we shall
also denote by $\Xi_1$. One can describe $\Xi_1$ as follows. An $H$-bundle over
$E\times(E\otimes_\Zee\Lambda)$ is the same thing as an element of
$\Pic(E\times(E\otimes_\Zee\Lambda))\otimes \Lambda$. The inclusion
$$\Pic (E\times E)\otimes \Lambda^* \to \Pic(E\times(E\otimes_\Zee\Lambda))$$
induces an inclusion
$$\Pic (E\times E)\otimes \Lambda^*\otimes\Lambda = \Pic (E\times E) \otimes
\Hom(\Lambda, \Lambda) \to \Pic(E\times(E\otimes_\Zee\Lambda))\otimes \Lambda,$$
and we take the image of the element $\mathcal{P}\otimes \Id$. As vector bundles
over $E\times (E\otimes _\Zee\Lambda)$,
$$\ad \Xi_1 = \scrO_{E\times (E\otimes _\Zee\Lambda)}^r \oplus \bigoplus _{\alpha
\in R}\mathcal{P}_\alpha,$$ where $\mathcal{P}_\alpha$ is the pullback to $E\times
(E\otimes _\Zee\Lambda)$ of the Poincar\'e bundle $\mathcal{P}$ over $E\times E$,
via the morphism induced from $\alpha$ from $E\otimes _\Zee\Lambda$ to $E$. Thus,
by functorial properties of determinant line bundles (cf. \cite[Chapter 5,
Proposition 3.8]{FMbook}), it will suffice to show that, over $E\times E$, $\det
R\pi _2{}_*(\mathcal{P}\otimes \pi _1^*\lambda) =\scrO_E(-e)$, where as before
$\lambda = \scrO_E(e-p_0)$. It is clear in any case that the inverse of $\det R\pi
_2{}_*(\mathcal{P}\otimes \pi _1^*\lambda)$ is represented by an effective divisor
supported at $e$, and the only question is the length of
$R^1\pi _2{}_*(\mathcal{P}\otimes \pi _1^*\lambda)$. A standard calculation
using the Grothendieck-Riemann-Roch theorem shows that this length is one (cf.
\cite[Chapter 7, Lemma 1.6]{FMbook} for the case of the trivial line bundle). 
\end{proof}

Next we identify the divisor $\sum _{\alpha \in R}D_{\lambda, \alpha}$. Using the
identifications
$$H^2(E\otimes _\Zee\Lambda; \Zee) = \bigwedge^2H^1(E\otimes _\Zee\Lambda; \Zee) = 
\bigwedge^2\left(H^1(E; \Zee)\otimes \Lambda^*\right)$$ 
there is an inclusion
$$\bigwedge ^2H^1(E; \Zee )\otimes \Sym ^2\Lambda^* \subseteq
\bigwedge^2\left(H^1(E; \Zee)\otimes \Lambda^*\right),$$ 
and hence,  since there is a canonical identification  $\bigwedge^2H^1(E; \Zee)
\cong \Zee$ there is a natural
inclusion of $\Sym^2\Lambda^*$ in $H^2(E\otimes _\Zee\Lambda;
\Zee)$. Let $Q\in \Sym^2\Lambda^*$ be the quadratic form described in Section
1 defined by
$$Q= \sum _{\alpha \in R}\langle \alpha, \cdot\rangle \langle \alpha,
\cdot\rangle,$$ 
and let $\widehat Q$ be the corresponding element of $H^2(E\otimes _\Zee\Lambda;
\Zee)$. By Lemma~\ref{Looform}, $Q = (2g)I_0$, where
$I_0$ is the unique
$W$-invariant quadratic  form on $\Lambda$ such that
$I_0(\tilde \alpha \spcheck,\tilde \alpha \spcheck) =2$.

\begin{lemma}\label{5.8} $-c_1(\mathbb{D}_3) = \widehat Q$.
\end{lemma}
\begin{proof} By Lemma~\ref{5.7}, it clearly suffices to show that, for every
$\alpha
\in R$, we have an equality (under the obvious identifications)
$$[D_{\lambda, \alpha}] = \langle \alpha, \cdot\rangle \langle \alpha,
\cdot\rangle.$$ As such, this equality is a general fact about lattices $\Lambda$:
suppose that
$\Lambda$ is a lattice and $\alpha\colon  \Lambda \to \Zee$ is a homomorphism.
There is an associated morphism $E\otimes _\Zee\Lambda \to E$ which we shall 
denote by $e_\alpha$. We can define the divisor
$D_\alpha= e_\alpha ^*(p)$ for $p\in E$, as well as the cohomology class $\langle
\alpha, 
\cdot\rangle \langle \alpha, \cdot\rangle \in \Sym^2\Lambda ^* \subset H^2(E;
\Zee)$. To prove Lemma~\ref{5.8}, it is enough to prove:

\begin{claim}\label{5.9} The class of the divisor $D_\alpha= e_\alpha ^*(p)$ is
equal to
$\langle  \alpha, 
\cdot\rangle \langle \alpha, \cdot\rangle$.
\end{claim}
\noindent\textbf{Proof of Claim~\ref{5.9}.}  First assume that $\alpha$ is
primitive. Then after a suitable choice of a basis of $\Lambda$ we can assume that
$\Lambda
\cong
\Zee^r$ and that
$\alpha$ is projection onto the last factor. In this case, $D_\alpha =
E^{r-1}\times
\{p\}$ and 
$\langle  \alpha, 
\cdot\rangle \langle \alpha, \cdot\rangle  = \pi _r^*(e\wedge f)$, where $e\wedge
f$ is a positive generator of $H^2(E; \Zee)\cong \Zee$. Clearly, equality holds in
this case. If $\alpha$ is not primitive, we can write $\alpha = n\alpha _0$, where
$n$ is a nonnegative integer and $\alpha _0$ is primitive. In this case,
$e_\alpha$ factors as the morphism $\alpha_0$ followed by multiplication by $n$ on
$E$, and so $D_\alpha$ is cohomologous to $n^2$ copies of $D_{\alpha_0}$. Likewise
$$\langle  \alpha, 
\cdot\rangle \langle \alpha, \cdot\rangle=n^2\langle  \alpha_0, 
\cdot\rangle \langle \alpha_0, \cdot\rangle,$$ and so the claim follows from the
case where $\alpha$ is primitive.
\end{proof}

\subsection{The divisor on the moduli space: the non-simply connected case}

We now redo the above calculations to handle the non-simply connected case. We
have the moduli space
$\mathcal{M}(G,c)$ and the corresponding determinant line bundle as before, and we
use the notation of \S\ref{s1.5}. To make the calculation, we can pull the
determinant line bundle up to the space
$E \otimes \Lambda ^{w_c} =T_0\times T_0$, where there is a universal flat bundle.
Let $\mathbb{D}_3$ be the class of the determinant line bundle pulled back to $E
\otimes \Lambda ^{w_c}$. As before, we have an inclusion $\Sym^2(\Lambda
^{w_c})^* \to H^2(E \otimes \Lambda ^{w_c}; \Zee)$. Let $Q_0$ be the element
$(2g)(I_0|\Lambda ^{w_c})\in \Sym^2(\Lambda
^{w_c})^*$ and let $\widehat Q_0$ be the corresponding element of $H^2(E \otimes
\Lambda ^{w_c}; \Zee)$. We have the following analogue of Lemma~\ref{5.8}:

\begin{lemma}\label{5.10} $-c_1(\mathbb{D}_3) = \widehat Q_0$.
\end{lemma} 
\begin{proof} The proof is similar to that in the simply connected case for
Lemmas~\ref{5.7} and \ref{5.8}, and we shall be a little sketchy. Suppose that
$\xi$ is a flat
$ K$-bundle corresponding to the
$c$-pair
$(x,y)$. Let $\lambda$ be a fixed, general line bundle of degree zero on $E$. We
compute when $\xi$ is in the support of $(R\det\pi_2{}_*(\ad \Xi \otimes \pi
_1^*\lambda))^{-1} = -\mathbb{D}_3$. As we have seen in Lemma 4.5 of \cite{FM}, 
$$\ad \xi \cong (\frak h^{w_{c}} \otimes \scrO_E) \oplus V_0'
\oplus \bigoplus _{\mathbf{o}}(V_{y,\mathbf{o}} \otimes L_{x, \mathbf{o}}),$$
where the $\mathbf{o}$ are the orbits for the action of $w_c$ on  $R$. Here $V_0'$
is a sum of certain torsion line bundles, $L_{x, \mathbf{o}}$ is a line bundle with
holonomy
$\alpha (x)$ for any fixed choice of $\alpha
\in \mathbf{o}$, and $V_{y,\mathbf{o}}$ is the sum of the root spaces $\frak
g^\alpha,
\alpha \in \mathbf{o}$, with the action defined by $y$. It follows that
$V_{y,\mathbf{o}} \otimes L_{x, \mathbf{o}}$ is a direct sum of \textbf{distinct}
line bundles of degree zero. Given a
$w_c$-orbit
$\mathbf{o}$, let $\alpha _{\mathbf{o}}$ be a choice of $\alpha \in \mathbf{o}$.
Next we construct a universal bundle $\Xi_1$ as in the  simply connected case,
along the lines of \cite[Lemma 5.21]{FM}. The construction outlined there shows
that
$$\ad \Xi_1 = (\frak h^{w_{c}} \otimes \scrO_{E\times(E\otimes \Lambda ^{w_c})})
\oplus \pi_1^*V_0'
\oplus \bigoplus _{\mathbf{o}}(V_{y_0,\mathbf{o}} \otimes L_{x_0,
\mathbf{o}}\otimes \mathcal{P}_{\mathbf{o}}),$$
where $\pi_1\colon E\times(E\otimes \Lambda ^{w_c}) \to E$ is the projection
onto the first factor, $(x_0, y_0)$ is a fixed
$c$-pair and, as in the simply connected case,
$\mathcal{P}_{\mathbf{o}}$ is the pullback to $E\times(E\otimes \Lambda ^{w_c})$ of
the Poincar\'e bundle $\mathcal{P}$ on $E\times E$ via the morphism $E\otimes
\Lambda ^{w_c} \to E$ induced by $\alpha_{\mathbf{o}}$.  The proof of
Lemma~\ref{5.7} then shows that the divisor $\mathbb{D}_3$ is reduced.

For a general choice of
$\lambda$, there exist   $c_x$ and $c_y$ depending only on $\lambda$ such that
$\xi$ is in the support of  $-\mathbb{D}_3$ if and only if there exists an
$\mathbf{o}$ such that
\begin{align*}
\alpha _{\mathbf{o}}(x) &= c_x;\\
\sum _{\alpha \in \mathbf{o}}\alpha (y) &= c_y.
\end{align*} Thus, in cohomology,  $c_1(\mathbb{D}_3)$ corresponds to the element 
$$\sum _{\mathbf{o}}\left(\alpha _{\mathbf{o}} \otimes \sum _{\alpha \in
\mathbf{o}}\alpha\right)\in (\Lambda ^{w_c})^*\otimes (\Lambda ^{w_c})^*.$$ Now 
every $\alpha \in \mathbf{o}$ has the same restriction to  $\Lambda ^{w_c}$ as
$\alpha _{\mathbf{o}}$. Thus the above sum become
$$\sum _{\mathbf{o}}d_{\mathbf{o}}(\alpha _{\mathbf{o}}\otimes \alpha
_{\mathbf{o}}),$$ where $d_{\mathbf{o}}$ is the order of $\mathbf{o}$. On the
other hand, we have
$(2g)I_0 = \sum _{\alpha \in R}(\alpha \otimes \alpha)$. Clearly, the restriction
of this form to $\Lambda ^{w_c}$ is $\sum _{\mathbf{o}}d_{\mathbf{o}}(\alpha
_{\mathbf{o}}\otimes \alpha _{\mathbf{o}})$, and this completes the proof of
Lemma~\ref{5.10}.
\end{proof}

\subsection{Completion of the proof of Theorem~\ref{LooThm}}

We have identified a divisor $\mathbb{D}_1$ on $\WP(\eta_0)$ and computed its top
intersection. We have identified a divisor
$\mathbb{D}_2$ on
$\mathcal{M}(G,c)$, or rather its pullback to a divisor  $\mathbb{D}_3$ on a
finite cover of $\mathcal{M}(G,c)$. To find the top power of
$\mathbb{D}_3$, we use the next lemma:

\begin{lemma}\label{5.11} Let $J$ be a quadratic form on $\Lambda^*$, and view $J$
as an element of $H^2(E\otimes _\Zee\Lambda; \Zee)$ via the inclusion
$\Sym^2\Lambda ^*
\subset H^2(E; \Zee)$. If the rank of $\Lambda$ is $r$, then the top power of $J$
is $(r!)\det J$.
\end{lemma}
\begin{proof} Let $\Omega$ be the $2$-form corresponding to $J$. First suppose
that $J$ is diagonalizable with respect to some $\Zee$-basis for $\Lambda$,
corresponding to a given isomorphism $\Lambda \cong \Zee^r$. Then
$\Omega$ is of the form 
$$\sum _{i=1}^ra_i[(e\wedge f)\otimes (\pi _i)^2,$$ where $J= \sum _{i=1}^ra_i(\pi
_i)^2$, say,  $\pi_i\colon \Lambda \cong \Zee^r\to \Zee$ is projection onto the
$i^{\rm th}$ factor, and
$e\wedge f$ is the positive generator for
$H^2(E;
\Zee)$. In this case, we can write $\Omega = \sum _{i=1}^ra_i(e_i\wedge f_i)$,
where
$e_i\wedge f_i$ is the generator on the
$i^{\text{th}}$ factor of $E\otimes _\Zee\Lambda \cong E^r$. Clearly
$$\Omega ^r = (r!) a_1\cdots a_r = (r!)\det J.$$ In general, the statement makes
sense for $\Q $-coefficients. Note that 
$$\dim _{\Q }\bigwedge ^{2r}\left(H^1(E; \Q ) \otimes (\Lambda^*\otimes \Q
)\right) = 1,$$ and a basis element $b$ is given by choosing the standard positive
generator for $H^2(E\otimes _\Zee\Lambda; \Zee)$, together with a
$\Zee$-basis for $\Lambda$. Changing the $\Zee$-basis for $\Lambda$ to some new
$\Q $-basis changes the element
$b$ by $(\det X)^2$, where $X$ is the change of basis matrix. In particular, if
$X$ has determinant $1$, then $b$ is unchanged. Now every quadratic form on
$\Lambda$ can be diagonalized via a $\Q $-basis such that the change of basis
matrix relating the new $\Q $-basis to a $\Zee$-basis has determinant $1$. Thus we
may reduce to the case where $J$ is diagonalizable, where we have already checked
the result.
\end{proof}

\begin{corollary}\label{moddegree} $\int_{\mathcal{M}(G,c)}c_1(\mathbb{D}_2)^{r_c}
= (-2g)^{r_c}n_0\Big/\prod _{\ov \beta} g_{\ov \beta}$.
\end{corollary} 
\begin{proof} Let
$e$ be the degree of the covering $T_0\times T_0\to
\mathcal{M}(G,c)$. We see by Lemma~\ref{5.8} and Lemma~\ref{5.10} that it suffices
to prove that
$$(2g)^{r_c}n_0\Big/\prod _{\ov \beta} g_{\ov \beta} = (2g)^{r_c}(r_c)!\det
(I_0|\Lambda ^{w_c})/e,$$ which we can rewrite as 
$$e = \frac{(r_c)!}{n_0}\det (I_0|\Lambda ^{w_c})\prod _{\ov \beta} g_{\ov
\beta}.$$ This is exactly the statement of Theorem~\ref{degree}.
\end{proof}

To complete the proof of Theorem~\ref{LooThm}, it suffices to show that
$\Psi^*[\mathbb{D}_2]=[\mathbb{D}_1]$ in $H^2(\WP(\eta_0), \Q)$. There is a  
Zariski open and dense subset $\mathcal{M}_0$ of $\mathcal{M}(G,c)$ consisting of
semi\-stable
$G$-bundles for which the regular representative also carries a flat connection
\cite[Corollary 6.2]{FM}. Let $\widetilde {\mathcal{M}}_0$ be the preimage of this
subset in $E\otimes \Lambda^{w_c}= \widetilde {\mathcal{M}}$. 

\begin{lemma}\label{universal} Let $x\in \widetilde
{\mathcal{M}}_0$, and let
$\xi_x$ be the corresponding $G$-bundle. Then the tautological  bundle
$\Xi_1$ constructed above identifies an analytic neighborhood of $x\in
\widetilde
{\mathcal{M}}_0$ with the local semi-universal deformation space of $\xi_x$, which
is locally universal.
\end{lemma} 
\begin{proof} By the definition of $\widetilde {\mathcal{M}}_0$, there is a
unique representative up to isomorphism for the S-equivalence class of $\xi_x$ and
it is both flat and regular. Hence the map
$H^1(E; \frak h^{w_c}\otimes \scrO_E) \to H^1(E; \ad \xi_x)$ is an isomorphism.
Then the tautological  bundle
$\Xi_1$ constructed above identifies an analytic neighborhood of $x\in
\widetilde
{\mathcal{M}}_0$ with the local semi-universal deformation space of $\xi_x$. Since
$\Lie
\Aut \xi_x = \frak h^{w_c}$, it acts trivially on $H^1(E; \ad \xi_x)$ and thus the
local semi-universal deformation of $\xi_x$ is in fact locally universal. 
\end{proof} 

(See
\cite[Theorem 6.12]{FM} for a more general result along these lines.) We now claim:

\begin{lemma}\label{pullback} Let $X$ be irreducible, and let $\ov \Xi\to E\times
X$ be an
$\ad G$-bundle which lifts to a $G$-bundle on every slice $E\times \{x\}$. Let
$f\colon X \to \mathcal{M}(G,c)$ be the corresponding morphism, and suppose that
$f(X) \cap
\mathcal{M}_0 \neq \emptyset$. Let $\mathbb{D}_X =\det R\pi_2{}_*\ad \ov
\Xi$. Then $[\mathbb{D}_X] =f^*[\mathbb{D}_2]$ in $H^2(X;\Q)$.
\end{lemma}
\begin{proof} Choose a component $\widetilde X$ of the fiber product $X\times
_{\mathcal{M}(G,c)}\widetilde {\mathcal{M}}$, and let $\tilde f\colon \widetilde
X \to \widetilde {\mathcal{M}}$ be the induced map. It suffices to show that
$\tilde f ^*[\mathbb{D}_3]= [\mathbb{D}_{\widetilde X}]$ in the obvious notation.
Since both sides are algebraic, it suffices to show the following: let $\Sigma
\subseteq \widetilde X$ be an irreducible  curve such that $\tilde f(\Sigma) \cap
\widetilde {\mathcal{M}}_0 \neq \emptyset$. Then $\tilde f_*[\Sigma]\cdot
\mathbb{D}_3 =
\Sigma \cdot \mathbb{D}_{\widetilde X}$. Now choose a line bundle $\lambda$ of
degree zero on $E$ such that, if $\operatorname{div}$ is the canonical section of
$\det R\pi_2{}_*(\ad \Xi_1\otimes \pi_1^*\lambda)^{-1}$, then $\operatorname{div}
\cap \tilde f(\Sigma) \subseteq \widetilde {\mathcal{M}}_0$ is a finite set of
points. Let $p\in \Sigma $ be a point such that $\tilde f(p) \in \operatorname{div}
\cap f(\Sigma)$. Since $\tilde f(p) \in \widetilde {\mathcal{M}}_0$, it follows
from Lemma~\ref{universal}  that, in an analytic neighborhood
$\Omega$ of
$p$, the bundle $\ad \ov \Xi$ is pulled back via $\tilde f$ from   $\ad\Xi_1$. It
follows that the  line bundle $\big(\det R\pi_2{}_*(\ad \ov \Xi\otimes
\pi_1^*\lambda)\big)^{-1}$ and its canonical section are also pulled back via
$\tilde f$. Thus 
$\tilde f_*[\Sigma]\cdot
\mathbb{D}_3 =
\Sigma \cdot \mathbb{D}_{\widetilde X}$ as claimed.
\end{proof}

We cannot apply Lemma~\ref{pullback} directly to the morphism $\Psi\colon
\WP(\eta_0) \to
\mathcal{M}(G,c)$, since there is no universal $\ad G$-bundle over $\WP(\eta_0)$.
However, as in the proof of Lemma~\ref{5.3}, there is a finite cover of
$\WP(\eta_0)$ by a projective space $\Pee^{r_c}$ and an $\ad G$-bundle over
$\Pee^{r_c}$. Let
$\widehat \Psi \colon \Pee^{r_c} \to \mathcal{M}(G,c)$ be the induced morphism. It
follows from the proof of Lemma~\ref{5.3} that $\scrO_{\WP(\eta_0)}(-2gn_{c,
\alpha}/n_0)$ pulls back to $\scrO_{\Pee^{r_c}}(-2gn_{c, \alpha}/n_0)$ and that
this line bundle is the determinant line bundle on $\Pee^{r_c}$. By
Lemma~\ref{pullback}, $\widehat
\Psi ^*[\mathbb{D}_2] = c_1(\scrO_{\Pee^{r_c}}(-2gn_{c, \alpha}/n_0))$, and hence 
$\Psi ^*[\mathbb{D}_2] = [\mathbb{D}_1]$. Together with
Corollary~\ref{moddegree} and Corollary~\ref{5.5}, this completes the proof of
Theorem~\ref{LooThm}.

\appendix
\section*{Appendix: Nonabelian cohomology}
\setcounter{section}{1}

 We discuss the general formalism, due
to Grothendieck \cite{Groth}, for deciding when  a principal bundle with structure
group a quotient group can be lifted back to the full group, and for classifying
all such liftings. In general, the set of liftings (suitably interpreted) is given
by a nonabelian cohomology set. We then go on to discuss circumstances under which
this cohomology set has a natural scheme structure, and indeed represents an
appropriate functor. The arguments here are modifications of arguments due to
Deligne and Babbitt-Varadarajan
\cite{BV}, given in a somewhat different context.

\subsection{Lifting}

We begin with a very general discussion of nonabelian cohomology and liftings of
principal bundles. First suppose that $G$ is a linear algebraic group and that
$\xi$ is a principal $G$-bundle over $X$, where $X$ is a  scheme or analytic space.
Here it is understood that there is some topology for which $\xi$ is locally
trivial, for example, Zariski,
\'etale, or classical, and cohomology will always be computed with respect to the
appropriate topology. Throughout this paper, we have always worked with
holomorphic bundles and the classical topology. One basic result is that, if
$X$ is projective, there is a natural bijection between the set of isomorphism
classes of principal holomorphic $G$-bundle over
$X$ in the classical topology and the set of principal $G$-bundles over $X$ in the
\'etale topology. This follows from the method of proof of Prop.\ 20 in GAGA
\cite{GAGA} and the fact that, if $G$ is a closed   subgroup of $GL(n)$, then the
quotient
$GL(n)/G$ is quasiprojective and admits local cross-sections in the \'etale
topology. Thus the set of isomorphism classes of holomorphic principal $G$-bundles
over $X$ is canonically identified with the set of principal $G$-bundles over
$X$ in the \'etale topology. However, when we discuss representable functors
below and try to put a scheme structure on various cohomology sets, it will be
convenient to use the \'etale topology. For most of the paper, we have only
considered the case where 
$X$ is a smooth projective curve, and the issue of the correct topology is not
important. Indeed, it follows from \cite{Steinberg}
that, if $X$ is a smooth curve and $G$ is linear, then a  locally trivial
$G$-bundle in the
\'etale topology is actually Zariski locally trivial. However, we will not use this
fact. One fact about cohomology which we shall need is the following: if $X$ is a
scheme and
$V$ is a coherent sheaf on $X$, then $H^i(X; V)$ computed for the \'etale
topology is the usual sheaf cohomology computed in the Zariski topology
\cite{Milne}, III (3.8). Of course, by GAGA, a similar statement holds in the
classical topology for the analytic sheaf associated to $V$ provided that $X$ is
projective. 

If $S$ is a scheme on which $G$ acts, we can form the associated scheme 
$\xi\times_GS$ (not to be confused with fiber product). It is fibered over $X$ and
the fibers are isomorphic to
$S$. Denote the sheaf (of sets) of cross sections (regular, holomorphic, or
\'etale, depending on the context) of $\xi\times_GS$ by
$S(\xi)$. We will usually be interested in the case where $S$ is itself an
algebraic group and where $G$ acts on $S$ by homomorphisms. In this case,
$\xi\times_GS$ is a group scheme and $S(\xi)$ is a sheaf of (not necessarily
abelian) groups. For example, if $S=G$ and the action is by conjugation, then
$G(\xi)$ is the automorphism group scheme of $\xi$ and its global
sections are the group
$\Aut _G\xi$. If $S$ is a vector space and $G$ acts on $S$ linearly, then
$S(\xi)$ is the vector bundle associated to the corresponding
representation on $G$. Given an algebraic group $G$ and a space $X$, the sheaf of
morphisms from $X$ to $G$ will be denoted $\underline{G}$. Here $\underline{G}$ is
a sheaf in the Zariski, \'etale, or classical topology, depending on the context.

Suppose that
$G$ is an algebraic group and that
$N$ is a closed normal subgroup. Let $H = G/N$, with $\pi\colon G \to H$ the
induced morphism.  Let $X$ be a scheme, and let $\xi _0$ be a principal $H$-bundle
over $X$. Suppose that $\xi$ is a principal
$G$-bundle lifting $\xi_0$. Note that $G$ acts on $N$ by conjugation, so that
$N(\xi)$ is defined. If moreover $N$ is abelian, then this action of $G$ on $N$
factors through an action of $H$ on $N$, and so $N(\xi _0)$ is defined.

The group $H^0(X;H(\xi_0))$
acts on the cohomology set
$H^1(X; N(\xi))$. We have the following general result
\cite{Groth} or \cite{Serre}:

\begin{lemma}\label{A1} With the above notation, the set of all principal
$G$-bundles lifting $\xi_0$, or in other words the fiber of the class $[\xi_0]\in
H^1(X; \underline{H})$ under the natural map $H^1(X; \underline{G}) \to H^1(X;
\underline{H})$, may be identified with $H^1(X; N(\xi))/H^0(X;
H(\xi_0))$. \endproof
\end{lemma}

We will also want a slight variant of the above:

\begin{lemma}\label{A2} In the notation of Lemma~\ref{A1}, the set of all
isomorphism classes of pairs
$(\eta, \varphi)$, where
$\eta$ is a principal $G$-bundle and $\varphi$ is an isomorphism from $ \eta/N$  to
$\xi _0$, can be identified with $H^1(X; N(\xi))$. \endproof
\end{lemma}

Note that $H^0(X; H(\xi_0))$ is the group of global automorphisms of
$\xi _0$, and this group acts naturally on the set of pairs $(\eta, \varphi)$ as
above. In fact, this action is the same as the coboundary action of $H^0(X;
H(\xi_0))$ on $H^1(X; N(\xi))$.

Next, we ask the bundle $\xi_0$ lifts to an $G$-bundle. For example, suppose that
$G = N\rtimes H$ is a semidirect product of
$N$ and $H$. Then there is a natural lift of $\xi _0$ given by the choice
of an inclusion of
$H$ in $G$. In particular, the map $H^1(X; 
\underline{G}) \to H^1(X; \underline{H})$ is surjective. In this case, we can see
the identification of Lemma~\ref{A1} quite explicitly: suppose that $\xi _0$ is
defined by the cocycle $\{h_{ij}\}$ with respect to some open cover $\{U_i\}$ of
$X$. Viewing the
$h_{ij}$ as taking values in
$G$ via the inclusion, it is easy to see that, if $\xi$ is a $G$-bundle lifting
$\xi_0$ on
$H$, then we can assume that $\xi$ is given by transition functions of the form
$h_{ij}n_{ij}$. In order to be a $1$-cocycle, the $n_{ij}$ must satisfy
$$\left(h_{jk}^{-1}n_{ij}h_{jk}\right)n_{jk} = n_{ik},$$ 
which says that $\{n_{ij}\}$ defines an element of $H^1(X; N(\xi))$.
If two such cocycles, say $\{h_{ij}n_{ij}\}$ and $\{h_{ij}'n_{ij}'\}$ define
isomorphic $G$-bundles, then we can first arrange by a $1$-coboundary that
$h_{ij}=h_{ij}'$. In this case, if $\{h_{ij}n_{ij}\}$ and $\{h_{ij}n_{ij}'\}$ are
cohomologous, then there exist $h_i$ such that $h_i^{-1}h_{ij}h_j = h_{ij}$, so
that $\{h_i\}\in H^0(X; H(\xi_0))$, and moreover $n_{ij}' = h_j^{-1}n_{ij}h_j$, so
that the cocycles $\{n_{ij}\}$ and $\{n_{ij}'\}$ differ by the action of  $H^0(X;
H(\xi_0))$ on $H^1(X; N(\xi))$. Conversely, if $\{n_{ij}\}$ and $\{n_{ij}'\}$
differ by an element of $H^0(X;
H(\xi_0))$, then reversing the above argument shows that the corresponding
$G$-bundles are isomorphic.

For another example of the surjectivity of the map $H^1(X; \underline{G}) \to 
H^1(X; \underline{H})$, we
have:

\begin{lemma}\label{A3} Suppose that $N$  is abelian and that $H^2(X;
N(\xi_0)) = 0$, for example suppose that $N$ is a vector space and that $\dim X =
1$. Then the bundle $\xi_0$ lifts to a $G$-bundle.
\end{lemma}
\begin{proof} This follows from  \cite{Serre}, Corollary to Prop.\ 41, p.\ 70 or
\cite{Groth}.
\end{proof}

One trivial observation which we shall often use is the following:

\begin{lemma}\label{A4} Suppose  that we are given an exact sequence
$$1 \to U \to G \to H \to 1,$$
and that $\xi_0$ is a principal $H$-bundle over
$X$. Suppose that $U_0$ is a closed subgroup of the center of $U$ which is normal
in $G$. Finally suppose that $\widetilde \xi _0$ is a lift of $\xi_0$ to a
principal
$G/U_0$-bundle. Then the sheaf $U_0(\xi _0)$ defined by the  natural
action of $H$ on $U_0$ is isomorphic to $U_0(\widetilde \xi _0)$.
\end{lemma}
\begin{proof}  The groups $H$ and $G/U_0$ act on
$U_0$ by conjugation, and the action of $G/U_0$ factors through the projection to
$H$. Thus, the sheaves  $U_0(\xi _0)$ and $U_0(\widetilde \xi
_0)$ are identified as well.
\end{proof}

\subsection{Representability}

Let $G=LU$, where $U$ is a normal subgroup of $G$ and $G$ is a semidirect product
of $U$ and $L$.  Let $\xi_0$ be a principal $L$-bundle. We  want to find some
circumstances under which the points of the  cohomology set $H^1(X;
U(\xi_0))$ can be identified with the points in an affine space. In fact,  it is
important to prove a much stronger statement, that a certain functor corresponding
to the cohomology set is representable by an affine space. For example, if $U$ is a
vector space on which
$L$ acts linearly, then  $H^1(X; U(\xi_0))$ is an ordinary sheaf cohomology group
and thus is itself a vector space, and the corresponding affine space represents a
functor. We will encounter nonabelian groups $U$ which are unipotent. Thus,
$U$ has a filtration $U_N
\subset \cdots \subset U_1 = U$ by normal, $L$-invariant subgroups $U_i$
with the property that $U_i$ is contained in the inverse image in $U$ of the center
of $U/U_{i+1}$. The idea then, following the general
lines of \cite{BV}, will be to work inductively, starting with the case where $U$
is a vector group. The
inductive step depends on the following (\cite{BV}, Lemma 2.5.3):

\begin{lemma}\label{A5} Let $R$ be a ring. Suppose that $\mathbf{F}$ and
$\mathbf{G}$ are two covariant functors from the category of $R$-algebras to sets
and that
$\varphi\colon \mathbf{F}
\to \mathbf{G}$ is a morphism of functors, with the following property:
\begin{enumerate}
\item[\rm  (i)] $\mathbf{G}$ is represented by a polynomial algebra $R[x_1, \dots,
x_n]$ over the ring $R$.
\item[\rm  (ii)] For every $R$-algebra $S$ and for every object $\xi \in
\mathbf{G}(S)$, the
functor $\mathbf{F}_{\varphi, \xi}$ from $S$-algebras to sets defined by
$$T \mapsto \varphi (T)^{-1}(\xi'), $$
where $\xi'$ is the element of $\mathbf{G}(T)$ induced by $\xi$,
is represented by  a polynomial algebra $S[y_1, \dots, y_m]$ over $S$.
\end{enumerate}
Then the functor $\mathbf{F}$ is represented by $R[x_1, \dots, x_n, y_1, \dots,
y_m]$.
\end{lemma}
\begin{proof} Let $S= R[x_1, \dots, x_n]$ and let $\xi\in \mathbf{G}(S)$ correspond
to the identity in $\Hom _R(S,S)$. For an
$R$-algebra $T$, if $\eta \in \mathbf{F}(T)$, let $\xi' = \varphi (T)(\eta)$. Then
there exists a unique homomorphism $f\colon S \to T$ corresponding to $\xi' \in
\mathbf{G}(T)$, so that $T$ is an $S$-algebra. Now $\xi'$ is the image of $\xi\in
\mathbf{G}(S)$ under
$f_*$, and the element $\eta \in
\varphi (T)^{-1}(\xi')$ defines a unique homomorphism $S[y_1, \dots, y_m] \to T$.
Thus $\eta$ defines a unique homomorphism
$S[y_1, \dots, y_m] = R[x_1, \dots, x_n, y_1, \dots, y_m] \to T$. Conversely,
suppose given a homomorphism 
$$f\colon R[x_1, \dots, x_n, y_1, \dots, y_m] =S[y_1, \dots, y_m] \to T.$$
Then
$f$ induces a homomorphism $R[x_1, \dots, x_n]=S \to T$ and thus an element
$\xi'\in \mathbf{G}(T)$ induced by $\xi$, and the homomorphism $f$ then gives an
element
$\eta\in \mathbf{F}(T)$ mapping to
$\xi'$. Clearly these are inverse constructions. It follows that  $\mathbf{F}$ is
represented by 
$R[x_1, \dots, x_n, y_1, \dots, y_m]$.
\end{proof}

The proof shows more generally that, if $\mathbf{G}$ is represented by some
$R$-algebra $S$, and, for $\xi$ the object of $\mathbf{G}(S)$ corresponding to the
identity, if the functor from $S$-algebras to sets defined by
$$T \mapsto \varphi (T)^{-1}(\xi'), $$
where $\xi'$ is the element of $\mathbf{G}(T)$ induced by $\xi$,
is represented by an $S$-algebra $\widetilde S$, then $\widetilde S$ represents
$\mathbf{F}$.

We shall apply Lemma~\ref{A5} as follows. First let us define the functor
$\mathbf{F}$ from  
$\Cee$-algebras to sets corresponding to the group  $H^1(X; U(\xi_0))$. For a
$\Cee$-algebra $S$, let
$\mathbf{F}(S)$  be the set of isomorphism classes of pairs $(\Xi, \Phi)$, where
$\Xi$ is a principal $LU$-bundle over $X\times \Spec S$ and $\Phi:  \Xi/U \to \pi
_1^*\xi_0$ is an isomorphism from the principal $L$-bundle over $X\times \Spec S$
induced by
$\Xi$ to the pulled back bundle $\pi _1^*\xi_0$. Thus,
$$\mathbf{F}(S) = H^1(X\times \Spec S; U(\pi_1^*(\xi_0)).$$

\begin{theorem}\label{A6} Let $G= LU$ be an algebraic group over $\Cee$, where $U$
is a closed normal unipotent subgroup of $G$ and $G$ is isomorphic to the
semidirect product of $L$ and $U$. Let $X$ be a projective scheme, let $\xi _0$ be
a principal
$L$-bundle over $X$, and let
$U(\xi_0)$ be the corresponding sheaf of unipotent groups.  Let  $\{U_i\}_{i=1}^N$
be a decreasing  filtration of $U$ by normal $L$-invariant subgroups such that, for
every
$i$, 
$U_i/U_{i+1}$ is contained in the center of $U/U_{i+1}$.
Suppose that, for every $i$, 
$$H^0(X; (U_i/U_{i+1})(\xi _0)) = H^2(X; (U_i/U_{i+1})(\xi _0)) =0.$$
Then: 
\begin{enumerate}
\item[\rm  (i)] The cohomology set
$H^1(X; U(\xi_0))$ has the structure of  affine
$n$-space $\mathbb{A}^n$. More precisely,  there is a $G$-bundle $\Xi_0$
over
$X\times
\mathbb{A}^n$ and an isomorphism $\Phi_0\colon  \Xi_0/U \to \pi
_1^*\xi_0$ such that the pair $(\Xi_0, \Phi_0)$ represents the functor $\mathbf{F}$
defined above. 
\item[\rm  (ii)] There is a natural action of the  algebraic group
$\Aut_L \xi _0$ on the  affine
$n$-space $\mathbb{A}^n$ representing $H^1(X; U(\xi_0))$. This action lifts
to an action on $\Xi_0$.
\end{enumerate}
\end{theorem}
\begin{proof}  We claim that the functor $\mathbf{F}$ is
representable by an affine space. The proof is by induction on the length of the
filtration $\{U_i\}$. If this length is zero, then $U=\{0\}$ and there
is nothing to prove. Suppose that the claim has been verified for every group 
and filtration  of length less than
$N$, and let  $\{U_i\}$ be a filtration of $U$ length exactly $N$ satisfying the
hypotheses of  Theorem~\ref{A6}. If
$U_N$ is the first term in $\{U_i\}$, then the filtration $\{U_i/U_N\}$ of $U/U_N$
has length
$N-1$. By induction, the functor $\mathbf{G}$ corresponding to $LU/U_N$ is
representable. Moreover, there is an obvious morphism of functors $\varphi\colon
\mathbf{F}
\to \mathbf{G}$. Now suppose that $S$ is a $\Cee$-algebra and that we are given an
object
$\xi$ of
$\mathbf{G}(S)$, in other words a pair $(\Xi, \Phi)$, where $\Xi$
is a principal $LU/U_N$-bundle over $X\times \Spec S$ and $\Phi:   \Xi/U \to \pi
_1^*\xi_0$ is an isomorphism from the principal $L$-bundle over $X\times \Spec S$
induced by $\Xi$ to  $\pi _1^*\xi_0$. We define the  functor $\mathbf{F}_{\varphi,
\xi}$ on
$S$-algebras $T$ as follows: $\mathbf{F}_{\varphi, \xi}(T)$ consists of isomorphism
classes of pairs
$(P,\Psi)$, such that $P$ is a principal $LU$-bundle over $X\times \Spec T$
and $\Psi$ is an isomorphism from the principal
$LU/U_N$-bundle over $X\times \Spec T$ induced by $P$ to the pullback 
$\widetilde \Xi= (\Id \times f)^*\Xi$ of $\Xi$ to $X\times \Spec T$, where
$f\colon \Spec T \to \Spec S$ is the morphism corresponding to the homomorphism
$S\to T$. There is a natural map
$\mathbf{F}_{\varphi, \xi}(T) \to \varphi (T)^{-1}(\widetilde\Xi, (\Id \times
f)^*\Phi)\subseteq
\mathbf{F}(T)$. First we claim:

\begin{claim} There exists a a principal $LU$-bundle $P$ over $X\times \Spec T$
lifting the principal $LU/U_N$-bundle $\widetilde \Xi =(\Id \times f)^*\Xi$ over
$X\times
\Spec T$. In other words, for all $T$, 
$\mathbf{F}_{\varphi, \xi}(T) \neq \emptyset$.
\end{claim}
\begin{proof} By Lemma~\ref{A3}, the obstruction to finding such a lift
lives in the group $H^2(X\times \Spec T; U_N(\widetilde \Xi))$. By 
Lemma~\ref{A4}, we can identify
$U_N(\widetilde \Xi)$ with the sheaf $U_N(\pi_1^*\xi _0)$,
which is the pullback via $\pi_1^*$ of the vector bundle $V =
U_N(\xi_0)$ on $X$. Now since $\Spec T$ is affine,
$$H^2(X\times \Spec T; U_N(\widetilde \Xi)) = H^2(X\times \Spec T;
\pi _1^*V) = H^2(X; V)\otimes _\Cee T.$$
Since $H^2(X; V) =0$ by hypothesis, we can lift $\widetilde \Xi$ to a bundle $P$.
\end{proof}

Once we know that there exists one lift $P$ as in the claim, it follows from
 Lemma~\ref{A2} that the set of all such pairs $(P, \Psi)$ is classified by
$H^1(X\times \Spec T; U_N(P))$. Next we  claim that the map $\mathbf{F}_{\varphi,
\xi}(T)
\to \mathbf{F}(T)$ is injective, and thus identifies $\mathbf{F}_{\varphi, \xi}(T)$
with the fiber
$\varphi (T)^{-1}(\widetilde\Xi, (\Id \times f)^*\Phi)$. To see this, it follows
from the general formalism of nonabelian cohomology that there is a transitive
action of
$H^0(X\times \Spec T; U/U_N(P))$ on the fibers of the map from
$\mathbf{F}_{\varphi,
\xi}(T)=H^1(X\times \Spec T; U_N(P))$ to $\mathbf{F}(T) = H^1(X\times \Spec T;
U(\pi_1^*\xi_0))$, which identifies the fibers with the coset space $H^0(X\times
\Spec T; U/U_N(P))/\operatorname{Im}H^0(X\times \Spec T; U(P))$. In our case, the
unipotent group $H^0(X\times \Spec T; U/U_N(P))$ is filtered, with successive
quotients contained in $H^0( X\times \Spec T; (U_i/U_{i-1})(\pi_1^*\xi_0)) = 0$.
Hence  $\mathbf{F}_{\varphi, \xi}(T)
\to \mathbf{F}(T)$ is injective.

Using   Lemma~\ref{A4} and the fact
that
$\Spec T$ is affine,
$$H^1(X\times \Spec T; U_N(P)) = H^1(X\times \Spec T;
\pi _1^*V) = H^1(X; V)\otimes _\Cee T.$$
If $e_1, \dots, e_n$ is a basis  for $H^1(X; V)$, with dual basis $x_1, \dots,
x_n$, this says that $\mathbf{F}_{\varphi, \xi}(T) \cong \Hom _S(S[x_1, \dots,
x_n], T)$. Thus $\mathbf{F}_{\varphi, \xi}$ is representable by an affine space
over $\Spec S$. Applying  Lemma~\ref{A5} and induction, the functor $\mathbf{F}$ is
then representable by an affine space, in other words there exists an affine space
$\mathbb{A}^n =\Spec
\Cee[x_1, \dots, x_n]$ and an object $(\Xi_0, \Phi_0)\in \mathbf{F}(\Cee[x_1,
\dots, x_n])$ such that the couple 
$(\Cee[x_1, \dots, x_n], (\Xi_0, \Phi_0))$
represents $\mathbf{F}$.

The fact that the functor $\mathbf{F}$ defined above is representable implies in
particular that there is a universal $LU$-bundle over $X\times \mathbb{A}^n$, where
$\mathbb{A}^n$ is the affine space representing $\mathbf{F}$. A formal reduction
also shows that obvious extension of $\mathbf{F}$ to a functor from schemes of
finite type over
$\Cee$ to sets is also representable by  $\mathbb{A}^n$.

Lastly we must analyze  the action of $\Aut_L \xi _0$ on $H^1(X;
U(\xi _0))$. 

\begin{proposition}\label{A7} Suppose that $X$ is proper and that $L$ is a linear
algebraic group. Then $\Aut_L \xi _0$ is also a linear algebraic group, and the
natural set-theoretic action of $\Aut_L \xi _0$ on  $H^1(X; U(\xi _0))$ is an
algebraic action. Moreover, this action lifts to an action on  $\Xi_0$. 
\end{proposition}
\begin{proof} Let $\mathbf{A}$ be the functor from $\Cee$-algebras to groups
corresponding to
$\Aut_L
\xi _0$: for a $\Cee$-algebra $S$, $\mathbf{A}(S)$ is the group of automorphisms of
the pullback of $\xi _0$ to $X\times \Spec S$. Since $X$ is proper and $L$ is
affine,
$\mathbf{A}$ is representable by a linear algebraic group scheme $\Aut_L \xi _0$
over
$\Cee$. To see this, first assume that $L = GL_n$. If $V$ is the vector
bundle corresponding to $\xi _0$, $\Aut_L \xi _0$ is an affine open subset of the
affine space 
$H^0(X; End \,V)$ and  we claim that this linear algebraic group represents the
associated functor. Indeed, an automorphism of $\pi _1^*V$ is the same thing as a
section of $\pi _1^*End \, V$, in other words a  $\Spec S$-valued point $\varphi$
of
$H^0(X; End \,V)$, such that the determinant of $\varphi$ is an invertible
element of $S$, and this is the same thing as a morphism from $S$ to the Zariski
open subset of $H^0(X; End \,V)$ consisting of elements with nonzero determinant.
In general, choose an embedding of
$L$ as a closed subgroup of
$GL_n$ for some $n$, defined by polynomials $f_i$.  It is then straightforward to
verify that $\mathbf{A}$ is representable by a closed subgroup of the
corresponding     group scheme for $GL_n$.

If $\mathbf{F}$ is the functor associated to
$H^1(X; U(\xi _0))$, there is an obvious morphism of functors $\mathbf{A}\times
\mathbf{F} \to \mathbf{F}$: if the points of
$\mathbf{F}(S)$ correspond to pairs $(\Xi, \Phi)$, where $\Phi$ is an isomorphism
from the principal $L$-bundle over $X\times \Spec S$ induced by
$\Xi$ to  $\pi _1^*\xi_0$, then the automorphisms of $\pi _1^*\xi_0$ act by
composition with $\Phi$. Since $\mathbf{A}$ and $\mathbf{F}$ are representable by
$\Aut_L
\xi _0$ and  by the affine coordinate ring of $H^1(X; U(\xi _0))$ respectively,
there is a corresponding morphism
$\Aut_L \xi _0\times  H^1(X;U(\xi _0))\to H^1(X; U(\xi
_0))$, which is easily checked to give an action. It again follows formally by
representability that this action lifts to an action on $\Xi_0$.
\end{proof}

This concludes the proof of Theorem~\ref{A6}.
\end{proof}

\begin{remark} In the hypotheses of Theorem~\ref{A6}, suppose we only assume
that, for all $i$,  
$H^2(X; (U_i/U_{i+1})(\xi _0)) =0$ and that, for all $i> 1$, $H^1(X;
(U_i/U_{i+1})(\xi _0)) =0$. Then, in the inductive construction of the proof, the
fibers $\mathbf{F}_{\varphi, \xi}(T)$ are all a single point
and thus the map $\mathbf{F}_{\varphi,
\xi}(T)
\to \mathbf{F}(T)$ is automatically injective. Thus the proof goes through in this
case as well.
\end{remark}

There is also a relative version of Theorem~\ref{A6}. 

\begin{theorem}
Let $G= LU$ and the filtration $\{U_i\}$ be as in Theorem~\ref{A6}.
Suppose that $\pi\colon Z \to \Spec R$ is a flat proper morphism, and that  $\xi
_0$ is a principal $L$-bundle over $Z$, with  $V_i=
(U_i/U_{i+1})(\xi _0)$  the vector bundle associated to the action
of $L$ on $U_i/U_{i+1}$ and the principal $L$-bundle $\xi _0$. Suppose that 
$H^2(\pi^{-1}(t); V_i|\pi^{-1}(t)) = 0$ for every point $t\in X$, that
$H^0(\pi^{-1}(t); V_i|\pi^{-1}(t)) = 0$ for every point $t\in X$, and that the
$R$-module $H^1(Z; V_i)$ is locally free and compatible with base change, in the
sense that, for every homomorphism $R\to S$ of $\Cee$-algebras, with corresponding
morphism $f\colon \Spec S \to \Spec R$, we have 
$$H^1(Z\times _{\Spec R}\Spec S;
f^*V_i) \cong H^1(Z; V_i)\otimes _RS.$$
For example, if $\pi$ has relative
dimension one and, for every $i$, $\dim H^1(\pi^{-1}(t);
V_i|\pi^{-1}(t))$ is independent of $t\in \Spec R$, then $H^1(Z; V_i)$ is locally
free and compatible with base change in the above sense. Then 
\begin{enumerate}
\item[\rm (i)] There
exists a locally trivial bundle of affine spaces $\mathbb{A}$ over $\Spec R$, such
that the set of sections of $\mathbb{A}$ is isomorphic to the set 
$H^1(Z; U(\xi _0))$. 
\item[\rm (ii)] There exists a universal bundle $\Xi$ over $Z\times _{\Spec
R}\mathbb{A}$ in the obvious sense.
\item[\rm (iii)] The automorphism group scheme $\mathcal{A}$ of
$\xi_0$  acts on the bundle $\mathbb{A}$ of affine spaces over $\Spec R$, and
this action lifts to an action on $\Xi$. 
\end{enumerate}
\end{theorem}

\begin{remark} One
can also  replace $\Spec R$ in the above statements by a scheme of finite type
over $\Cee$. Moreover, in case $R =\Cee[t_1, \dots, t_j]$, or more generally if
every vector bundle over $\Spec R$ is trivial, then the inductive proof of
Theorem~\ref{A6} shows that we can take $\mathbb{A} =\Spec R[x_1, \dots, x_n]$.
\end{remark}

\bigskip
\noindent
Department of Mathematics \\
Columbia University \\
New York, NY 10027 \\
USA

\bigskip
\noindent
{\tt rf@math.columbia.edu, jm@math.columbia.edu}


\begin{thebibliography}{10}

\bibitem{Atiyah}
M. Atiyah,
\emph{Vector bundles over an elliptic curve},
Proc. London Math. Soc. \textbf{7}~(1957),
414--452.

\bibitem{AtBo} 
M. Atiyah and R. Bott,
\emph{The Yang-Mills  equations over Riemann surfaces},
Phil. Trans. Roy. Soc. London A \textbf{308}~(1982),
523--615.

\bibitem{ABS}
H. Azad, M. Barry, and G. Seitz,
\emph{On the structure of parabolic subgroups},
Comm. Algebra \textbf{18}~(1990),
551--562.

\bibitem{BV}
D. Babbitt and V.S. Varadarajan,
\emph{Local Moduli for Meromorphic Differential Equations},
Ast\'erisque \textbf{169--170}~(1989).

\bibitem{BS}
I.N. Bernshtein and O.V. Shvartsman,
\emph{Chevalley's theorem for complex crystallographic Coxeter groups},
Funct. Anal. Appl. \textbf{12}~ (1978), 
308--310.

\bibitem{Borel}
A. Borel,
\emph{Linear Algebraic Groups}, Second Enlarged Edition,
Graduate Texts in Mathematics \textbf{126},
New York, Springer-Verlag, 1991.

\bibitem{BFM}
A. Borel, R. Friedman, and J.W. Morgan,
\emph{Almost commuting elements in compact Lie groups},
Memoirs of the AMS (to appear),
math.GR/9907007.

\bibitem{Bour}
N. Bourbaki,
\emph{Groupes et Alg\`ebres de Lie}, 
Chap. 4, 5, et 6,
Masson,  Paris,
1981.

\bibitem{FMbook}
R. Friedman and J.W. Morgan,
\emph{Smooth Four-Manifolds and Complex Surfaces},
Ergebnisse der Mathematik und ihrer Grenzgebiete
3. Folge \textbf{27},  Springer-Verlag, Berlin Heidelberg New York,
1994.

\bibitem{FM}
R. Friedman and J.W. Morgan,
\emph{Holomorphic principal bundles over elliptic curves},
math.AG/9811130.

\bibitem{FMAB}
R. Friedman and J.W. Morgan,
\emph{On the converse to a theorem of Atiyah and Bott},
J. Algebraic Geometry (to appear), math.AG/0006086.

\bibitem{FMW}
R. Friedman, J.W. Morgan, and E. Witten,
\emph{Vector bundles and $F$ Theory},
Commun. Math. Phys. \textbf{187}~(1997),
679--743.

\bibitem{FMW1}
R. Friedman, J.W. Morgan, and E. Witten,
\emph{Principal $G$-bundles over elliptic curves},
Math. Research Letters \textbf{5}~(1998),
97--118.


\bibitem{FMW2}
R. Friedman, J.W. Morgan, and E. Witten,
\emph{Vector bundles over elliptic fibrations},
J. Algebraic Geometry \textbf{8}~(1999),
279--401.

\bibitem{Groth}
A. Grothendieck,
\emph{A general theory of fibre spaces with structure sheaf}, 
University of Kansas report \textbf{4}~(1955).

\bibitem{HS}
S. Helmke and P. Slodowy,
\emph{On unstable principal bundles over elliptic curves},
preprint.

\bibitem{Loo}
E. Looijenga,
\emph{Root systems and elliptic curves},
Invent. Math. \textbf{38}~(1976),
17--32.

\bibitem{Milne}
J.S. Milne,
\emph{\'Etale Cohomology},
Princeton University Press, Princeton, NJ, 1980.

\bibitem{Mori}
S. Mori,
\emph{On a generalization of complete intersections}, 
J. Math Kyoto Univ. \textbf{15}~(1975), 
619--646.

\bibitem{GIT}
D. Mumford,
\emph{Geometric Invariant Theory}, Third Edition, Ergebnisse der Mathematik und
ihrer Grenz\-gebiete Neue Folge \textbf{34}, Berlin, Springer-Verlag,
1994.

\bibitem{Ra}
A. Ramanathan,   
\emph{Stable principal bundles on a compact Riemann surface},
 Math. Ann. \textbf{213}~(1975),
129--152.

\bibitem{GAGA}
J-P. Serre, 
\emph{G\'eom\'etrie alg\'ebrique et g\'eom\'etrie analytique},
Ann. Inst. Fourier \textbf{6}~(1956),
1--42.

\bibitem{Serre}
J-P. Serre, 
\emph{Cohomologie Galoisienne},
Lecture Notes in Mathematics \textbf{5}~(1973).

\bibitem{Steinberg}
R. Steinberg,
\emph{Regular elements of semisimple algebraic groups},
Publ. Math. Inst. Hautes \'Etudes Sci. \textbf{25}~(1965),
49--80.

\end{thebibliography}
\end{document}